\documentclass{amsart}
\usepackage{epsf}
\usepackage{amssymb,latexsym, amsmath, amscd, array, graphicx}
\usepackage{fullpage}
\def\hepsffile{\leavevmode\epsffile}

\numberwithin{equation}{section}

\swapnumbers
\theoremstyle{plain}

\newtheorem{thm}{Theorem}[section]
\newtheorem{cor}[thm]{Corollary}
\newtheorem{theorem}[thm]{Theorem}
\newtheorem{lem}[thm]{Lemma}

\newtheorem{prop}[thm]{Proposition}
\newcommand\theoref{Theorem~\ref}
\newcommand\lemref{Lemma~\ref}
\newcommand\propref{Proposition~\ref}

\newcommand\defref{Definition~\ref}

\theoremstyle{definition}
\newtheorem{defin}[thm]{Definition}
\newtheorem{glos}[thm]{Glossary}

\newtheorem{rem}[thm]{Remark}
\newtheorem{remark[thm]}{Remark}
\newtheorem{rems}[thm]{Remarks}
\newtheorem{ex}[thm]{Example}
\newtheorem{example}[thm]{Example}
\newtheorem{con}[thm]{Construction}
\newtheorem{conven}[thm]{Convention}

\def\Diff{\protect\operatorname{Diff}}

\def\grad{\protect\operatorname{grad}}
\def\id{\protect\operatorname{id}}
\def\Im{\protect\operatorname{Im}}

\def\Int{\protect\operatorname{Int}}

\def\sign{\protect\operatorname{sign}}

\def\span{\protect\operatorname{span}}

\def\aug{\protect\operatorname{aug}}

\def\SKY{\mathcal{SKY}}

\def\alk{\protect\operatorname{alk}}

\def\pr{\protect\operatorname{pr}}

\def\clk{\protect\operatorname{lk}}
\def\exp{\protect\operatorname{exp}}

\def\Emb{\protect\operatorname{Emb}}
\def\det{\protect\operatorname{det}}
\def\pt{\protect\operatorname{pt}}
\def\Indet{\protect\operatorname{Indet}}
\def\Hess{\protect\operatorname{Hess}}
\def\bul{^{\scriptstyle{\scriptstyle{\bullet}}}}

\newcommand \pa[2]{\frac{\partial #1}{\partial #2}}

\def\theoref{Theorem~\ref}
\def\propref{Proposition~\ref}
\def\lemref{Lemma~\ref}

\def\eps{\varepsilon}
\def\gf{\varphi}
\def\ga{\alpha}
\def\gb{\beta}
\def\gl{\lambda}

\def\Z{{\mathbb Z}}
\def\Q{{\mathbb Q}}
\def\R{{\mathbb R}}
\def\N{{\mathbb N}}
\def\A{{\pmb A}}

\def\e{{\mathbf e}}
\def\bor{{\Omega}}
\def\cs{{\mathcal S}}

\def\ss{{X}}
\def\B{{\mathcal B}}
\def\ds{\displaystyle}

\def\1{\hbox{\rm\rlap {1}\hskip.03in{\rom I}}}
\def\Bbbone{{\rm1\mathchoice{\kern-0.25em}{\kern-0.25em}
{\kern-0.2em}{\kern-0.2em}I}}

\def\p{\partial}

\def\pp{\medskip{\parindent 0pt \it Proof.\ }}
\def\wt{\widetilde}
\def\wh{\widehat}
\def\m{\medskip}
\def\ov{\overline}
\def\un{\underline}
\long\def\forget#1\forgotten{} %
\begin{document}
\date{March 27, 2007}
\leftline{ } \centerline{ }
\title[Linking and causality in globally hyperbolic space-times]
{Linking and causality in globally hyperbolic space-times}
\author[V.~Chernov (Tchernov) and Yu.~Rudyak]{Vladimir V.
Chernov (Tchernov) and Yuli B. Rudyak}
\address{V. Chernov, Department of Mathematics,
6188 Bradley Hall, Dartmouth College, Hanover NH 03755, USA}
\email{Vladimir.Chernov@dartmouth.edu}
\address{Yu. Rudyak, Department of Mathematics, University
of Florida, 358 Little Hall, Gainesville, FL 32611-8105, USA}
\email{rudyak@math.ufl.edu}

\thanks{2000 Mathematics Subject Classification: Primary 57Q45, 53C50, Secondary 53C80, 53Z05, 57R17,
83C75}

\keywords{linking numbers, wave fronts, causality, degree of a
mapping, intersection index, globally hyperbolic space-time,
nonrefocussing space-time, Lorentz manifold, null geodesic}

\Large

\begin{abstract}
The classical linking number $\clk$ is defined when link components
are zero homologous. In~\cite{ChernovRudyakGT} we constructed the
affine linking invariant $\alk$ generalizing $\clk$ to the case of
linked submanifolds with arbitrary homology classes. Here we apply
$\alk$ to the study of causality in Lorentzian manifolds.

Let $M^m$ be a spacelike Cauchy surface in a globally hyperbolic
space-time $(\ss^{m+1}, g)$. The spherical cotangent bundle $ST^*M$
is identified with the space $\mathcal N$ of all null geodesics in
$(\ss,g).$ Hence the set of null geodesics passing through a point
$x\in X$ gives an embedded $(m-1)$-sphere $\mathfrak S_x$ in
$\mathcal N=ST^*M$ called the sky of $x.$ Low observed that if the
link $(\mathfrak S_x, \mathfrak S_y)$ is nontrivial, then $x,y\in
\ss$ are causally related. This motivated the problem (communicated
by Penrose) on the Arnold's 1998 problem list to apply the machinery
of knot theory to the study of causality. The spheres $\mathfrak
S_x$ are isotopic to the fibers of $(ST^*M)^{2m-1}\to M^m.$ They are
nonzero homologous and the classical linking number $\clk(\mathfrak
S_x, \mathfrak S_y)$ is undefined when $M$ is closed, while
$\alk(\mathfrak S_x, \mathfrak S_y)$ is well defined. Moreover,
$\alk(\mathfrak S_x, \mathfrak S_y)\in \Z$ if $M$ is not an
odd-dimensional rational homology sphere.

We give a formula for the increment of $\alk$ under passages through
Arnold dangerous tangencies. If $(\ss,g)$ is such that 
$\alk$ takes values in $\Z$ and $g$ is conformal to $\wh g$ that 
has all the timelike
sectional curvatures nonnegative, then $x, y\in \ss$ are causally related if and only if $\alk
(\mathfrak S_x, \mathfrak S_y)\neq 0$. We prove that if $\alk$ takes
values in $\Z$ and $y$ is in the causal future of $x,$ then
$\alk(\mathfrak S_x, \mathfrak S_y)$ is the intersection number of
any future directed past inextendible timelike curve to $y$ and of
the future null cone of $x.$ We show that $x,y$ in a nonrefocussing
$(\ss, g)$ are causally unrelated if and only if $(\mathfrak S_x,
\mathfrak S_y)$ can be deformed to a pair of $S^{m-1}$-fibers of
$ST^*M\to M$ by an isotopy through skies. Low showed that if $(\ss, g)$ is refocussing, then $M$ is compact. We show that the universal cover of $M$ is also compact.
\end{abstract}

\maketitle

\section{Preliminaries}\label{prelim}

We work in the $C^{\infty}$-category, and the word ``smooth'' means
$C^{\infty}$. An {\it isotopy\/} of a smooth embedding $f:P\to Q$ is 
a path in the space of smooth embeddings $P\to Q$ starting at $f.$
Given an oriented manifold $M^m$, consider its tangent
bundle $TM \to M$ and put $\mathfrak z:M\to TM$ to be the zero
section. Let $\R^+$ be the group of positive real numbers under
multiplication that acts on $TM$ as $(r,\mu)\mapsto r\mu, r\in \R^+,
\mu\in TM$. We put $STM=\bigl (TM\setminus \mathfrak
z(M)\bigr)/\R^+$ and note that the tangent bundle $TM \to M$ yields
the {\em spherical tangent bundle $\pr:STM\to M$\/} of $M.$

For the reasons discussed right before Theorem~\ref{mainthmalk}, we will assume 
that $\dim M>1.$

We denote by $T^*M \to M$ the {\em cotangent bundle} over $M$, and
we construct the {\em spherical cotangent bundle $\pr: ST^*M\to M$
\/} in a similar way. It is well known that $ST^*M$ possesses a
canonical contact structure and that the $S^{m-1}$-fibers of $ST^*M$
are Legendrian submanifolds with respect to this contact structure,
see \cite{Arnoldbook} or Appendix \ref{reviewcontact}. Note also
that the orientation of $M$ yields canonical orientations on the
fibers of spherical (co)tangent bundles. Namely, it is well known
that every spherical (co)tangent bundle is canonically oriented, and
we orient a fiber $S^{m-1}$ via the convention that the orientation
of $ST^*M$ is given by the pair (orientation of the base $M$,
orientation of the fiber $S^{m-1}$).

Given a path $\alpha:[a,b]\to M$, consider the bundle $E \to [a,b]$
induced by $\alpha$ from $ST^*M\to M$. So, we have the commutative
diagram
$$
\CD
E @ >\wh \alpha >> ST^*M\\
@VVV @VVV\\
[a,b] @>\alpha >> M.
\endCD
$$
Choose a trivialization $\iota: S^{m-1}\times [a,b] \to E$ of the
bundle $E \to [a,b]$. Define
\begin{equation}\label{epspath}
\CD \eps_{\alpha}: S^{m-1}\times [a,b] @>\iota >> E @>\wh \alpha >>
ST^*M.
\endCD
\end{equation}

Given a point $v\in M$, consider the constant path $\alpha: [0,1]\to
M, \ga(0)=v$ and define
\begin{equation}\label{epspoint}
\eps_v: S^{m-1} \longrightarrow ST^*M, \quad
\eps_v(s)=\eps_{\alpha}(s,0).
\end{equation}

Two such maps $\eps_v$ and $\eps'_v$ are isotopic via an isotopy 
such that the projection of its trace lies in a small disk containing $v.$ 
If the group $\Diff^+(S^{m-1})$ of degree one autodiffeomorphisms of
$S^{m-1}$ is connected, then such an isotopy can be chosen so that its trace 
is inside $\pr^{-1}(v).$
For example this holds for $S^3$, see~\cite{Cerf}, for $S^{2}$, 
see~\cite{Munkres, SmaleS2} and for $S^1,$ 
for trivial reasons. So since $\dim M>1$, 
any two links $(\eps_{u_1}, \eps_{v_1})$ and $(\eps_{u_2}, \eps_{v_2}),
u_1\neq v_1, u_2\neq v_2,$ are isotopic.

\def\mfootz{\footnote{Since the sequence $\pi_0(\Diff^+(D^m))\to \pi_0(\Diff^+(S^{m-1}))\to \Gamma_m\to 0$ is exact and the twisted sphere groups $\Gamma_5, \Gamma_6$ are zero, every degree one autodiffeomorphism of $S^{m-1}, m=5, 6$ extends to an  autodiffeomorphism of the unit disk $D^m\subset \R^m.$ By the results of Palais, Cerf, Milnor~\cite[Theorem 9.6]{Milnorhcobordism} every orientation preserving embeddings of $D^m$ to $\R^m$ is ambient isotopic to the identity map. Hence every degree one autodiffeomorphism of the standard unit $S^{m-1}\subset \R^m$ is isotopic to the identity map, cf~\cite[Remark page 122]{Milnorhcobordism}. Now one uses 
the exponent map $\exp_v$ and the front projection description of Legendrian knots in $ST^*M,$ see Example~\ref{contactSTM}, to get the proof.
}}

If $m=2, 3, 4$ then $\pi_0(\Diff^+(S^{m-1}))=0,$ and hence 
any two embeddings $\eps_v$ and $\eps'_v$ are Legendrian isotopic via an isotopy 
whose trace is contained in $\pr^{-1}(v).$ 
One can show~\protect\mfootz that if $m=5, 6$ then any $\eps_v$ and $\eps'_v$ are Legendrian isotopic via an isotopy such that the projection of its trace is 
inside of a small disk containing $v.$ 
For these cases any two links $(\eps_{u_1}, \eps_{v_1})$ and 
$(\eps_{u_2}, \eps_{v_2}),
u_1\neq v_1, u_2\neq v_2,$ are Legendrian isotopic.

\begin{defin}\label{df:model}
Let $f,g: S^{m-1}\to ST^*M$ be two embeddings with disjoint images
that are homotopic to a map $\eps_w$ for some $w\in M^m$. We say
that the pair $(f,g)$ is {\em unlinked} or {\em trivially linked }
if there exists a path $\gamma$ in the space of smooth embeddings
$S^{m-1} \sqcup S^{m-1}\to ST^*M$ that joins $(f,g)$ to a pair
$(\eps_u, \eps_v), u,v\in M, u\ne v$. 

If both embeddings $f, g:S^{m-1}\to ST^*M$ are Legendrian, we say that the pair
$(f,g)$ is {\em Legendrian unlinked} or {\em Legendrian trivially
linked} if there exists a path $\gamma$ as above in the space of
smooth Legendrian embeddings. (Any two trivial links are isotopic, but for $m\neq 2, 3, 4, 5, 6$ we do not 
know if it may happen that two Legendrian trivial links are not Legendrian 
isotopic.)

\end{defin}

\begin{defin}
A {\it vector field\/} on a manifold $Y$ is a smooth section of the tangent bundle $\tau_{Y}:TY\to Y$, and a {\it vector field along a $($smooth$)$ map\/} $\phi:Y_1\to Y_2$ of one manifold to another is a smooth map $\Phi:Y_1\to TY_2$ such that $\phi=\tau_{Y_2}\circ \Phi.$ 
Covector (direction, codirection, line, etc.) fields on a manifold and along a map $\phi$ are defined in a similar way. For brevity we will often write ``$\phi$ is equipped with vector field'' rather than ``$\phi$ is equipped with a vector field along it'', etc.
\end{defin}

Now we recall some basic concepts of Lorentzian geometry.

\begin{defin}\label{basic}
(a) Consider a smooth manifold $\ss^{m+1}$ equipped with a Lorentz
metric $g.$ A nonzero vector $\xi\in T\ss$ is called {\em timelike,
non-spacelike, null $($lightlike$)$, {\rm or} spacelike\/} if
$g(\xi,\xi)$ is negative, non-positive, zero, or positive,
respectively. A piecewise smooth curve is called {\em timelike,
non-spacelike, null, {\rm or} spacelike\/} if all of its velocity
vectors are respectively timelike, non-spacelike, null, or
spacelike. A smooth submanifold $M^m\subset \ss^{m+1}$ is {\em
spacelike\/} if the restriction of $g$ to $M$ is a Riemannian
metric.

(b) For each $x\in \ss$ the set of all non-spacelike vectors in
$T_x\ss$ consists of two connected components that are hemicones. A
continuous (with respect to $x\in X$) choice of a hemicone of
non-spacelike vectors in $T_x\ss$ is called the {\em time
orientation\/} of $(\ss, g).$

(c) The non-spacelike vectors from the chosen hemicones are called
{\em future pointing vectors.\/}
A piecewise smooth curve is said to be {\em future directed\/} if all of its
velocity vectors are future pointing.
\end{defin}

\begin{defin}\label{futurepast}
(a) A {\em space-time\/} $\ss=(\ss^{m+1}, g)$ is a smooth connected
time-oriented Lorentz $(m+1)$-manifold without boundary. An {\em
event\/} is a point of the space-time $\ss.$

(b) Two Lorentz metrics $g$ and $\wh g$ on $\ss^{m+1}$ are {\it conformal\/} 
if $\wh g=\Omega^2 g$ for some nowhere zero smooth function 
$\Omega:\ss\to \R.$ If $g$ and $\wh g$ are conformal, 
then a vector $\xi\in T\ss$ is timelike, nonspacelike, null, or
spacelike for $g$ if and only if it is
timelike, nonspacelike, null, or spacelike for $\wh g,$ respectively. 

(c) For two events $x,y\in \ss$ we write $x<\!\!<y$ if there is a
piecewise smooth future directed timelike curve from $x$ to $y.$ We
write $x\leq y$ if $x=y$ or if there is a piecewise smooth future
directed non-spacelike curve from $x$ to $y.$ For $x\in (\ss,g)$ we
put the {\em causal future of $x$\/} to be $J^+(x)=\{y\in \ss |
x\leq y\}$ and we put the {\em causal past of $x$\/} to be
$J^-(x)=\{y\in \ss |y\leq x\}.$ We put the {\em chronological future
of $x$\/} to be $I^+(x)=\{y\in \ss | x<\!\!< y\}$ and we put the
{\em chronological past of $x$\/} to be $I^-(x)=\{y\in \ss | y<\!\!<
x.\}$ It is easy to see that the causal and the chronological past and future of
$x$ depend only on the conformal class of the metric $g$ on $\ss.$

(d) Two events $x,y$ are {\em causally related\/} if $x\in J^+(y)$ or $y\in
J^+(x).$

(e) An open neighborhood is {\em causally convex\/} if there are no
non-spacelike curves intersecting it in a disconnected set. A
space-time is {\em strongly causal\/} if every point in it has
arbitrarily small causally convex neighborhoods. A strongly causal
space-time $(\ss, g)$ is {\em globally hyperbolic\/} if $J^+(x)\cap
J^-(y)$ is compact for all $x,y\in \ss.$

(f) A {\em Cauchy surface\/} $M$ is a subset of $\ss$ such that for
every inextendible non-spacelike curve $\gamma(t)$ in $\ss$ there
exists exactly one value $t_0$ of $t$ with $\gamma(t_0)\in M$.
\end{defin}

Clearly $y\in J^+(x)$ if and only if $x\in J^-(y);$ and $y\in
I^+(x)$ if and only if $x\in I^-(y).$ The sets $I^{\pm}(x)$ are
always open, see~\cite[Lemma 3.5]{BeemEhrlichEasley}; however, the
sets $J^{\pm}(x)$ are in general neither closed nor open,
see~\cite[pages 5--6]{BeemEhrlichEasley}. A space-time can be shown
to be globally hyperbolic if and only if it admits a Cauchy surface,
see~\cite[pages 211-212]{HawkingEllis}.

\begin{ex}\label{static}
Fix a connected oriented Riemannian manifold $(M,\ov g)$. Put
$\cdot$ to be the standard Riemannian metric on $\R^1.$ Define the
Lorentz metric $g$ on $M\times \R$ via $g\bigl( (\xi_1, \eta_1),
(\xi_2, \eta_2)\bigr)=\ov g(\xi_1, \xi_2)-\eta_1\cdot \eta_2,$ for
$(\xi_i, \eta_i)\in T_pM\times T_t\R=T_{(p,t)} (M\times \R), i=1,2.$
We denote the space-time $(M\times \R, g)$ by $(M\times \R, \ov
g\oplus -dt^2)$ and call it a {\em static space-time.\/} If $(M, \ov
g)$ is a complete Riemannian manifold, then the static space-time
$(M\times \R, \ov g\oplus -dt^2)$ is globally hyperbolic,
~see~\cite[Theorem 3.66]{BeemEhrlichEasley}, and each $M \times t$
is a Cauchy surface.
\end{ex}

\def\mfootu{\footnote{The definition of the Cauchy surface Bernal and Sanchez
use in~\cite{BernalSanchez} looks a bit weaker than the one we use.
They define Cauchy surface to be a subset of $\ss$ that is
intersected exactly once by every inextendible timelike curve,
rather than by every inextendible non-spacelike curve as we do.
However as Sanchez explained to us, their spacelike Cauchy surface
would be a Cauchy surface in our sense. This is since, every
non-spacelike curve intersects the Cauchy surface in their sense at
least once, see~\cite[Section 14, Lemma 29]{ONeill}. Moreover since
their Cauchy surface is spacelike, every non-spacelike curve would
intersect it at most once, see~\cite[Section 14, Lemma 42]{ONeill}.
}}

\m The pioneer result of Geroch~\cite{Geroch} says that every
globally hyperbolic space-time $\ss$ is {\em homeomorphic} to
$\Sigma^m\times \R$ where every $\Sigma\times t\subset \ss$ is a
(topological) Cauchy surface. The question of existence of smooth
Cauchy surfaces was considered by Seifert~\cite{Seifert} and
Dieckmann~\cite{Dieckmann}, cf.~\cite[Remark 3.77]{MinguzziSanchez}.

Bernal and Sanchez~\cite[Theorem 1]{BernalSanchez},~\cite[Theorem
1.1]{BernalSanchezMetricSplitting},~\cite[Theorem
1.2]{BernalSanchezFurther} proved the following strong result:

\begin{thm}\label{strongtheorem}
For every globally hyperbolic space-time $\ss^{m+1}$ there is an
isometry and in the same time a diffeomorphism
\begin{equation*}
h: (M^m\times \R, -\beta dt ^2+\ov g)\to (\ss, g)
\end{equation*}
where $M$ is a smooth manifold, $t:M\times \R\to \R$ is the projection,
$\beta:M\times \R\to (0, +\infty)$ is a smooth function and $\ov g$
is a smooth $2$-covariant symmetric tensor field on $M\times \R$
satisfying the following conditions$:$

\begin{description}
\item[1] for each $q\in M\times \R$ the vector $\grad t\in T_q(M\times \R)$ is
timelike and past pointing$;$
\item[2] for each $t$ the submanifold $M\times t$ of $M \times \R$ is a smooth
space-like Cauchy surface~\protect\mfootu $($i.e. it is a Cauchy
surface and the restriction of $-\beta d t ^2+\ov g$ to it is a
Riemannian metric$);$
\item[3] for each $q\in M\times \R$, the radical of $\ov g$ at $q$ is equal to
$$
\span\{\grad t\}=\span \{\partial/\partial t\}\subset T_q(M\times
\R).
$$
Here the radical of $\ov g$ at $q$ is the subspace of $T_q(M\times
\R)$ consisting of vectors that are $\ov g$ orthogonal to all
vectors in $T_q(M\times \R).$
\end{description}

In particular, the vector $\partial/\partial t$ is time-like and
future pointing, the function $t$ is increasing along all future
pointing non-spacelike curves, and the vector $\partial/\partial t$
is everywhere $(-\beta d t ^2+\ov g)$-orthogonal to the smooth
spacelike Cauchy surfaces $M\times t$ of $M \times \R$.

Moreover, for every smooth spacelike Cauchy surface $M\subset \ss$
there is an isometry $h:(M\times\R, \beta dt^2+\ov g)\to \ss $ as
above such that $h(m,0)=m,$ for all $m\in M.$

Also any two smooth spacelike Cauchy surfaces $M_1, M_2$ of $\ss$
are diffeomorphic.
\end{thm}

\begin{conven}
Throughout the paper the space-time $(\ss, g)$ is assumed to be
globally hyperbolic and oriented. The term ``Cauchy surface'' always
means ``smooth space-like Cauchy surface''. Since by our definition
every space-time is connected, Geroch Theorem~\cite{Geroch} implies
that a Cauchy surface is connected.
\end{conven}

\begin{defin}\label{df:proper}
Given a space-time $(\ss,g)$ and a Cauchy surface $M\subset \ss$,
let $h: M \times \R \to \ss$ be an isometry as in
\theoref{strongtheorem}.

We say that the isometry $h$ is {\em $M$-proper} if $h(m,0)=m$ for
all $m \in M$.

Given $t\in \R$, we let $M_t=h(M\times t)\subset \ss$ and define
$h_t: M\to M_t, h_t(m)=h(m,t)$ for all $m\in M$. Furthermore, we put
$g_t=h_t^{-1}: M_t \to M$. We put $ST^*g_t:ST^*M_t\to ST^*M,
STg_t:STM_t\to STM, ST^*h_t:ST^*M\to ST^*M_t, STh_t:STM\to STM_t$ to
be the induced maps.

Given a Cauchy surface $M\subset \ss$, choose an $M$-proper $h: M
\times \R\to \ss$ and orient $M\times \R$ by requiring $h$ to be
orientation-preserving. Now, we orient $M$ so that the pair
(orientation of $M$, orientation of $\R$) gives the orientation of
$M\times \R.$ Here the orientation of $\R$ is given by the time
orientation of $\ss$.
\end{defin}

\begin{defin}
Let $\mathcal N$ denote the space of all future directed
null-geodesics in $(\ss, g)$ modulo orientation preserving affine
reparameterizations. The {\it sky\/} $\mathfrak S_x\subset \mathcal
N$ of an event $x\in \ss$ is the subspace of all future directed
null-geodesics passing through $x.$
\end{defin}

\begin{example}\label{ex:sameskies} Let $\ov g$ be the metric on $S^m$
induced by the identification of $S^m$ with the unit sphere in
$\R^{m+1}.$ Since $(S^m, \ov g)$ is complete, the $(m+1)$-dimensional 
Einstein cylinder $(\ss,g):=(S^m\times \R, \ov g\oplus -dt^2)$ is globally 
hyperbolic, see ~\cite[Theorem 3.66]{BeemEhrlichEasley}. Given $s\in S^m$ and
$t\in \R$, put $x=(s,t)$ and $x'=(s, t+2\pi)$. Then, clearly,
$\mathfrak S_{x}=\mathfrak S_{x'}$ although $x\ne x'$.
\end{example}

\m For our future goals (see Section~\ref{introd}) we need the
concept of linking for skies. To this aim we relate skies in
$\mathcal N$ with (lifted) wave fronts in the spherical cotangent
bundle of a Cauchy surface. We explain this below.

Fix a Cauchy surface $M\subset \ss.$ An inextendible future directed
null-geodesic $\gamma$ intersects $M$ in one point $x=x(\gamma).$
Since $M$ is a spacelike surface, a $g$-orthogonal to $M$ line field
$L_y, y\in M,$ is not tangent to $M.$ Since $g$ is non-degenerate,
the lines $L_y, y\in M,$ do not contain null vectors. Thus
$T_{x}\ss=T_{x}M\oplus L_{x}$ and
\begin{equation}\label{decomposition}
\dot \gamma (x)=\xi+\eta \in T_{x}M\oplus L_{x}=T_{x}\ss,\quad
\xi\in T_xM,\,\eta\in L_x,
\end{equation}
with $\xi \neq 0$, $\eta\neq 0$, and $g(\xi, \eta)=0.$ In this way
we get a bijective map
\begin{equation}
\gf=\gf_M: \mathcal N \to STM,
\end{equation}
where $\gf(\gamma)$ is the point of $STM$ corresponding to the nonzero vector
$\xi.$

Since $M$ is a space-like surface, $g|_{M}$ is a Riemannian metric.
This allows us to identify $STM$ with the total space $ST^*M$ of the
spherical cotangent bundle, that has the natural contact structure.
Thus for a space-like Cauchy surface $M$ we get a bijective map
\begin{equation}\label{eq:psi}
\psi=\psi_{M}:\mathcal N\to ST^*M
\end{equation}
that equips $\mathcal N$ with the structure of a smooth contact
manifold.

\m Low showed~\cite{LowLegendrian} that if $M$ and $M'$ are two
smooth Cauchy surfaces, then the map
$$
f^{M'}_M=\psi_{M}\circ \psi_{M'}^{-1}:ST^*M'\to \mathcal N\to ST^*M
$$
is a contactomorphism.

(Strictly speaking the work~\cite{LowLegendrian} deals only with
$3+1$-dimensional space-times. However this result holds for
space-times of all the dimensions, see~\cite[pages
252-253]{NatarioTod}.)

\begin{defin}\label{df:fronts}
Let $M'$ be a Cauchy surface in $\ss$ and $x\in M'$. We put $\wt
W_{x, M'}=\eps_x:S^{m-1} \to ST^*M'$. Now, for an arbitrary (smooth
spacelike) Cauchy surface $M\subset \ss$, we put $\wt W_{x,
M}=f^{M'}_M\circ\eps_x :S^{m-1}\to ST^*M$. We call this embedding
{\em a lifted wave front of $x$\/} (with respect to the Cauchy
surface $M$). A {\em wave front\/} $W_{x, M}=\pr\circ\wt W_{x,
M}:S^{m-1}\to M$ is the projection of the lifted front $\wt W_{x,
M}$ to a smooth spacelike Cauchy surface $M.$
\end{defin}

The lifted wave front $\wt W_{x, M}$ is a Legendrian embedding
$S^{m-1}\to ST^*M.$ Indeed, if $M'$ is a smooth spacelike Cauchy
surface passing through $x$ that exists by \theoref{strongtheorem},
then $\wt W_{x, M'}:S^{m-1}\to ST^*M'$ is $\eps_x$ that is
Legendrian, see \cite{Arnoldbook} or Appendix \ref{reviewcontact}.
Since $f^{M'}_M$ is a contactomorphism, we get that $\wt W_{x,
M}=f^{M'}_M\circ \wt W_{x, M'}=f^{M'}_M\circ \eps_x$ is Legendrian.
(This explanation of why a lifted wave front is Legendrian was given
to us by Jose Natario, and we are grateful to him for it.)

\m Skies and lifted wave fronts are related as follows. Let
$\psi=\psi_{M'}$. Since $x\in M'$, we conclude that ${\psi(\mathfrak
S_x)}$ is the unit cotangent sphere $S_x$ that is the fiber of
$ST^*M'\to M'$ over $x$. Identifying $\mathfrak S_x$ with $S^{m-1}$,
we conclude that $\psi|_{\mathfrak S_x}: S^{m-1} \to ST^*M'$ is the
map $\eps_x=\wt W_{x, M'}$. We see that for each Cauchy surface $M$
the lifted wave front $\wt W_{x, M}=f^{M'}_M\circ \wt W_{x, M'}$ is
completely determined by the sky $\mathfrak S_x$. Moreover, if we
know $\wt W_{x, M}$ for some $M,$ then we can restore the sky
$\mathfrak S_x.$ Note that there are examples of space-times where
$\mathfrak S_x=\mathfrak S_y,$ for some $x\neq y$. Hence it is not
generally possible to restore $x$ from $\mathfrak S_x$ or from $\wt
W_{x, M}$.

\m A front $W_{x, M}$ is equipped with the natural codirection field defining the  lifting $\wt W_{x, M}$ to $ST^*M.$ In view of the identification
$ST^*M =STM$ (given by the Riemannian metric on $M$), this
codirection field yields a direction field. Since $\wt W_{x,
M}:S^{m-1}\to STM=ST^*M$ is Legendrian, this direction field
is everywhere orthogonal to the front with respect to the Riemannian
metric $g|_{M}$, cf. Example~\ref{contactSTM}.

In terms of skies, this (co)direction field can be described as
follows. For every (equivalence class of a) null-geodesic $\gamma\in
\mathfrak S_x$ the direction $\gf_M(\gamma)\in STM$ is orthogonal to
$W_{x,M}$. So, the direction and codirection fields are
$\{\gf_M(\gamma)\bigm| \gamma \in \mathfrak S_x\}$ and
$\{\psi_M(\gamma)\bigm| \gamma \in \mathfrak S_x\}$, respectively.

\m We see that there is no essential difference between the sky
$\mathfrak S_x$ and the lifted wave front $W_{x,M},$ for a Cauchy
surface $M.$ In fact, skies enable us to formulate the results in
an elegant invariant way (without making a choice of a Cauchy
surface), while the wave fronts play the role of technical tools
that are useful for proofs.

\begin{defin}\label{df:linkskies}
Assume that the events $x$ and $y$ do not lie on a common null
geodesic in a space-time $(\ss, g).$ Then the skies $\mathfrak S_x$
and $\mathfrak S_y$ are disjoint and hence for every Cauchy surface
$M\subset X$ the link $(\wt W_{x, M}, \wt W_{y, M})$ is nonsingular.
We say that the pair of skies ($\mathfrak S_x, \mathfrak S_y)$ is
{\em unlinked} or {\em trivially linked} (respectively {\em
Legendrian unlinked} or {\em Legendrian trivially linked}) if, for
every Cauchy surface $M\subset \ss$, the pair of lifted wave fronts
$(\wt W_{x,M}, \wt W_{y,M})$ is unlinked (respectively Legendrian
unlinked), as defined in \defref{df:model}.
\end{defin}

\begin{rems}\label{goodremark}
\par 1. If the events $x$ and $y$ lie on a common null geodesic, then $x$ and
$y$ are causally related for trivial reason. In this case $\mathfrak
S_x\cap \mathfrak S_y\neq \emptyset$ and hence for every Cauchy
surface $M$ the link $(\wt W_{x, M}, \wt W_{y, M})$ is singular. So,
the assumption that $x$ and $y$ do not lie on a common null geodesic
does not lead to loss of generality.

\par 2. In \theoref{welldefined} we show that the skies are unlinked provided
that the pair of lifted wave fronts $(\wt W_{x,M_0}, \wt W_{y,M_0})$
is unlinked for any {\em particular} Cauchy surface $M_0$.

\par 3. As we see from Example~\ref{ex:sameskies}, it can happen
that $\mathfrak S_x=\mathfrak S_y$ for some $x\ne y$. In this case, for each Cauchy surface $M,$ we have $\wt W_{x, M}=\wt W_{y,M}$ up to reparameterization.
\end{rems}

\begin{con}\label{con:fam}
To make the picture more familiar for some of the readers, choose a
Cauchy surface $M$ and an $M$-proper isometry $h: M \times \R \to
\ss$. Given $x\in \ss$, define
$$
\wt W_x^t=\wt W_{x,M}^t=ST^*g_t\circ \wt W_{x, M_t}: S^{m-1} \to
ST^*M_t \to ST^*M.
$$
Then the family $W_x^t:=\pr\circ \wt W_x^t: S^{m-1} \to M, t\in \R$
can be regarded as the wave front propagating in $M$. Note that for
$x\in M_t$ we have $\wt W_{x,M}^t=\eps_{g_t(x)}$.
\end{con}

\m Below we list notation and concepts that are used in the
paper consistently.

\begin{glos}\label{glos}
The basic concepts of Lorentz geometry; space-time (assumed to be
globally hyperbolic and oriented), conformal Lorentz metrics, causality, 
future and past directed curves, Cauchy surfaces (assumed to be smooth and space
like), the sets $J^{\pm}$ and $I^{\pm}$, etc, - are defined in
Definitions~\ref{basic} and \ref{futurepast}.

Differential-geometric notions for Lorentz manifolds are briefly
reminded in
\defref{difgeomLorentz}.

Given a globally hyperbolic space-time $\ss$ and a Cauchy surface
$M\subset \ss$, put $h: M \times \R \to \ss$ to be an isometry as in
\theoref{strongtheorem}. We denote by $\pi_M=\pi_{M, h}:\ss\to M$ the
composition of $h^{-1}$ and of the projection $M \times \R \to M.$
Similarly $\pi_{\R}=\pi_{\R, h}:\ss\to \R$ is the composition of
$h^{-1}$ and of the projection $M\times \R\to \R.$

Given $t\in \R$, we put $M_t=h(M\times t)\subset \ss$. We put
$g_t=\pi_M|_{M_t}: M_t \to M$ and $h_t=g_t^{-1}$, cf.
\defref{df:proper}

The maps $\eps_{\ga}$ and $\eps_v$ are described in \eqref{epspath}
and \eqref{epspoint}, respectively.

The lifted wave fronts $\wt W_{x,M}$ and wave fronts $W_{x,M}$ are
described in
\defref{df:fronts}. For the description of propagating wave fronts $\wt
W^t_{x,M}$ and $W^t_{x,M}$ see \ref{con:fam}.

Given a smooth curve $\gamma: \R \to Y$ and $t_0\in \R$, we denote
by $\dot \gamma(t_0)$ the velocity vector of $\gamma$ at $t_0$.
\end{glos}

\section{Introduction and Results}\label{introd}

Low~\cite{Low0} noticed that two events (in a globally hyperbolic
space-time) are causally related if their skies are linked in
$\mathcal N$. We explain this in greater detail in
Section~\ref{s:linkedcausal}.

The Low observation yielded the Question $8$ ``Causality in Terms of
Linking'' on V.I.~Arnold 1998 Problem List~\cite{ArnoldProblem}
which is to apply the machinery of knot theory to the study of the
relation between linking and causality. The problem was communicated
by Penrose.

\m {\bf Our paper is motivated by the above questions.} We study
relations between link theory and causality. Here we have the
following three directions of research.

\m {\bf 1. Detecting of linking.} Given two skies, how can we
recognize whether they are linked or not?

\m {\bf 2. Suitable space-times.} Low conjectured that if the Cauchy
surface is a $2$-disk with holes, then two events $x,y$ are causally
related if and only if the skies $\mathfrak S_x$ and $\mathfrak S_y$
are linked. Some special cases of this conjecture were proved by
Natario and Tod~\cite{NatarioTod}. For $(m+1)$-dimensional
space-times with $m>2$ the obvious extension of the Low conjecture
fails: Low~\cite{Low0} constructed an example of two causally
related events $x,y$ in a $(3+1)$-dimensional globally hyperbolic
space-time with Cauchy surface diffeomorphic to $\R^3$ such that the
pair $(\mathfrak S_x,\mathfrak S_y)$ is unlinked. So an interesting
question is: For which space-times the skies of every two causally
related events are linked? One of results in this direction is
\theoref{mainthmcurvature}.

\m {\bf 3. Suitable isotopies.} It is not currently known whether
the skies in the mentioned above example of Low~\cite{Low0} are
Legendrian unlinked. The modified Low conjecture posed by Natario
and Tod~\cite{NatarioTod} says that for $(3+1)$-dimensional globally
hyperbolic space-times, whose Cauchy surface $M^3$ is diffeomorphic
to a submanifold of $\R^3$, two events are causally related if and
only if their skies are Legendrian linked. 

Nevertheless the following is an example of two causally related
events in a globally-hyperbolic space-time whose skies are unlinked
even in the Legendrian sense. 

\begin{ex}\label{einstein}
Let $\ov g$ be the metric on $S^m$
induced by the identification of $S^m$ with the unit sphere in
$\R^{m+1}.$ Since $(S^m, \ov g)$ is complete, the $(m+1)$-dimensional Einstein 
cylinder $(S^m\times \R, \ov
g\oplus -dt^2)$ is globally hyperbolic, see ~\cite[Theorem
3.66]{BeemEhrlichEasley}. Given $s\in S^{m}, t\in \R$ put $x=(s, t),
x'=(s, t+2\pi)\in S^{m}\times \R.$ Put $n=-s\in S^m\subset \R^{m+1}$
and put $y=(n, t+2\pi)\in S^m\times \R.$ It is easy to see that the
events $x$ and $y$ are causally related but $(\mathfrak S_x,
\mathfrak S_y)=(\mathfrak S_{x'}, \mathfrak S_y)$ are Legendrian
unlinked. Since $S^m$ is not a submanifold of $\R^m,$ this example
contrasts but does not contradict the modified Low conjecture of
Natario and Tod~\cite{NatarioTod}.
\end{ex}

We can, however, consider a link isotopy that is even finer than
Legendrian. Namely, in Section~\ref{s:weakened} we prove that, for
so-called nonrefocussing space-times, two events $x,y$ are causally
related if and only if the link $(\mathfrak S_{x},\mathfrak S_{y})$
is isotopic through skies to the trivial link, \theoref{problem7}.

\m {\bf Now we explain the results in greater detail.} We
start with detecting that the link $(\mathfrak S_x, \mathfrak S_y)$
is nontrivial. One of the goals of the paper is to define a
generalized linking number of the pair of skies $(\mathfrak
S_{x},\mathfrak S_{y})$ that vanishes (i.e.~is equal to zero) for
unlinked pairs. In particular, the events $x$ and $y$ are causally
related if this invariant does not vanish for the pair $(\mathfrak
S_{x},\mathfrak S_{y})$. Note that for many space-times the
vanishing of our invariant implies that $(\mathfrak S_{x},\mathfrak
S_{y})$ is unlinked, see \theoref{mainthmcurvature}.

\m The (Gauss) linking number $\clk$ is the classical invariant that
often allows one to detect that the link is nontrivial. It is
defined as the intersection number of the singular chain whose
boundary is one of the linked manifolds with the other linked
manifold. In order for $\clk$ to be well-defined, the two linked
submanifolds have to be zero homologous and the sum of their
dimensions should be by one less than the dimension of the ambient
space. The link $(\wt W_{x,M}, \wt W_{y,M})$ consists of two copies
of $S^{m-1}$ in $(ST^*M)^{2m-1}$. Since $(m-1)+(m-1)+1=(2m-1)$, the
linking number $\clk(\wt W_{x,M}, \wt W_{y, M})$ would be
well-defined if $\wt W_{x, M}, \wt W_{y, M}$ were zero homologous.
Unfortunately, $\wt W_{x, M}, \wt W_{y, M}$ are homotopic to a
positively oriented $S^{m-1}$-fiber of $\pr:STM\to M$ which
generally is not zero homologous, and thus the linking number $\clk
(\wt W_{x, M}, \wt W_{y, M})$ is undefined.

\m 
When $M^m=\Int P^m,$ for some manifold $P$ with $\partial P\neq \emptyset,$ 
one can take two auxiliary negatively oriented 
fibers $S^{m-1}_{p_1}, S^{m-1}_{p_2}$ of $ST^*P\to P$ over two distinct points 
$p_1, p_2\in 
\partial P$ and define a 
(modified) linking number $\ov \clk (\wt W_{x,M}, \wt W_{y, M})=\clk \Bigl( (\wt 
W_{x, M}\sqcup 
S^{m-1}_{p_1}), (\wt W_{y, M} \sqcup S^{m-1}_{p_1})\Bigr).$ This was exactly the 
trick used by Low~\cite{Low0},~\cite{Low1},~\cite{Low2} to define 
his linking numbers of the skies in $ST^*\R^m$. Before Low this way of defining 
linking numbers for nonzero homologous circles in $ST\R^2$ was used by 
S.~Tabachnikov~\cite{Tabachnikov}. 
The general theory of linking numbers when the linked objects are zero 
homologous in the homology group of the ambient manifold modulo boundary was 
developed by 
U.~Kaiser~\cite{Kaiserbook}. When $M$ is a closed manifold, the number $\ov 
\clk$ defined 
using auxiliary negative fibers over some points is not an invariant of the link 
$(\wt W_{x, M}, \wt W_{y, M}),$ since it changes when a link component passes 
through the 
auxiliary fiber corresponding to the other link component.

\m In~\cite{ChernovRudyakGT} we constructed the {\it affine linking
invariant\/} that should be thought of as the generalization of the
linking number $\clk$ to the case of linked oriented submanifolds
realizing arbitrary homology classes.

In this paper we use this theory to define the affine linking number
$\alk(\wt W_{x,M}, \wt W_{y,M}).$ This $\alk$ invariant does not
depend on the Cauchy surface $M,$ see Theorem~\ref{welldefined}.
Hence it is an invariant of the two events $x, y$ (that do not lie
on a common null geodesic) and it can be interpreted as the affine
linking number $\alk(\mathfrak S_x,\mathfrak S_y)$ of the skies
$\mathfrak S_x,\mathfrak S_y$.

\m Actually, it was a very preliminary
version~\cite{ChernovRudyakFronts} of this paper that motivated our
work~\cite{ChernovRudyakGT}. Since~\cite{ChernovRudyakGT} was
published before this work, we rewrote it to avoid reproving results
proved in~\cite{ChernovRudyakGT}. We also changed the setup of the
work to be more familiar to people working in Lorentz geometry and
included new results about space-times for which $\alk$ completely
determines causality, about the relations between $\alk$ and the
intersection index, and about space-times for which a weakened Low
conjecture holds. The general $\alk$ invariant constructed
in~\cite{ChernovRudyakGT} takes values in a group that depends on
the ambient manifold and on the homotopy classes of the linked
submanifolds. The computation of the group is quite hard.

\m {\it In the rest of the paper we assume that space-times
have dimension $>2,$ and hence that Cauchy surfaces have dimension
$>1.$\/} The reason is the following. For a $2$-dimensional globally
hyperbolic space-time its Cauchy surface $M$ is $1$-dimensional and
the lifted wave fronts are maps of $S^0$. Since $S^0$ is not
connected,~\cite[Theorem 7.4 and Corollary 7.5]{ChernovRudyakGT}
that give a homotopy theoretical description of the range of values
of the $\alk$-invariant do not apply. Luckily in this case the Cauchy
surface $M$ is $\R$ or $S^1$ and all the links in $STM$ are
easily classified by combinatorial methods. 

\m Combining Theorem~\ref{thmfrontsunlinked}, Theorem~\ref{Main},
Proposition~\ref{usefulprop} and Theorem~\ref{welldefined} of this
work we get the following result.

\begin{theorem}\label{mainthmalk}
Let $M^m, m>1,$ be a Cauchy surface in a globally hyperbolic
space-time. Then the following holds:
\begin{enumerate}
\item If $M$ is not an odd-dimensional rational homology sphere with finite
$\pi_1(M),$ then $\alk(\mathfrak S_x,\mathfrak S_y)$ is a well
defined $\Z$-valued invariant of the link $(\mathfrak S_x,\mathfrak
S_y);$
\item If $M$ is an odd-dimensional rational homology sphere, then the invariant
$\alk(\mathfrak S_x,\mathfrak S_y)$ is well-defined if one regards
it as having values in $\Z/(\Im \deg).$ Here $\deg:\pi_m(M^m)\to \Z$
is the homomorphism that maps $[\alpha]\in \pi_m(M^m)$ to the degree
of $\alpha: S^m\to M^m;$
\item The only manifolds for which this quotient $\Z/\Im \deg$ is the trivial
group are odd-dimensional homotopy spheres;
\item $\alk(\mathfrak S_x, \mathfrak S_y)=0$ if $x,y$ are causally unrelated.
\end{enumerate}
\end{theorem}

\begin{rem}[$\alk$ as the universal order $\leq 1$ Vassiliev-Goussarov
invariant] In~\cite[Subsection 3.2]{ChernovRudyakGT} we proved that
our affine linking invariants are Vassiliev-Goussarov invariants of
order $\leq 1$ that are universal in the sense that they distinguish
all the link homotopy classes that can be distinguished using order
$\leq 1$ invariants with values in an abelian group. Since the
invariant $\alk$ constructed in this paper is a particular case of
the general construction from~\cite{ChernovRudyakGT}, we get that
$\alk(\wt W_{x, M}, \wt W_{y, M})=\alk(\mathfrak S_x,\mathfrak S_y)$
{\em is a universal Vassiliev-Goussarov link homotopy invariant of
order $\leq 1$ \/} of two linked $S^{m-1}$-spheres in $ST^*M$ that
are homotopic to a positively oriented fiber $S^{m-1}$ of
$\pr:ST^*M\to M^m, m>1.$
\end{rem}

\m A $2$-plane $E_s\subset T_s\ss$ is called {\em timelike\/} if
$g|_{E_s}$ is nondegenerate and not positive definite. A {\it timelike sectional 
curvature\/} is a sectional curvature along a timelike $2$-plane, see 
Definition~\ref{difgeomLorentz}. We prove the
following result:

\begin{theorem}[see Theorem~\ref{curvature}]\label{mainthmcurvature}
Let $(\ss, g), \dim \ss>2,$ be a globally hyperbolic space-time 
where $g$ is conformal to $\wh g$ that has all the timelike sectional curvatures 
nonnegative. Assume moreover that a Cauchy surface $M$ of $(\ss, g)$
is such that $\alk$ is a $\Z$-valued
invariant $($see Theorem~$\ref{mainthmalk})$. Then two events $x,y\in \ss$\ 
$($that do not lie on the
same null-geodesic$)$ are causally related if and only if
$\alk(\mathfrak S_x,\mathfrak S_y)\neq 0.$ In particular, they are
causally unrelated if and only if $(\mathfrak S_x,\mathfrak S_y)$
is a trivial link in $\mathcal N$.
\end{theorem}

\begin{ex}\label{exampletimelike}
Take a complete connected oriented Riemannian manifold $(M^m, \ov
g), m>1,$ of non-positive sectional curvature such that $M$ is not
an odd-dimensional rational homology sphere with finite $\pi_1(M).$
Consider a globally hyperbolic static space-time $(M\times \R, \ov
g\oplus -dt^2)$ as in Example~\ref{static}. Using~\cite[Equation
(3.21)]{BeemEhrlichEasley} one immediately gets that $(M\times \R,
\ov g\oplus -dt^2)$ has nonnegative sectional curvature on every
timelike two-plane. By Theorem~\ref{mainthmalk} $\alk(\mathfrak S_x,
\mathfrak S_y)$ is a $\Z$-valued invariant. Thus $(M\times \R, \ov
g\oplus -dt^2)$ satisfies all the conditions of
Theorem~\ref{mainthmcurvature} and two events $x,y\in (M\times \R,
\ov g\oplus -dt^2)$ (that do not lie on the same null-geodesic) are
causally related if and only if $\alk(\mathfrak S_x,\mathfrak
S_y)\neq 0.$
\end{ex}

The following Theorem shows that for $y\in J^+(x)$ the invariant
$\alk(\mathfrak S_x, \mathfrak S_y)$ gives an estimate from below on
the number of times the light rays from $x$ cross a generic past
inextendible timelike curve to $y.$

\begin{theorem}[see Theorem~\ref{alkasintersection}]
Let $(\ss^{m+1}, g), m>1$ be a globally hyperbolic space-time. Assume moreover 
that a Cauchy surface $M\subset X$ is such
that $\alk(\mathfrak S_x, \mathfrak S_y)$ is a $\Z$-valued
invariant. Let $x,y\in X$ be events that do not belong to a common
null geodesic and such that $y\in J^+(x).$ Then $\alk(\mathfrak S_x,
\mathfrak S_y)$ equals to the intersection index of the null-cone
consisting of the future directed null geodesics from the point $x$
and of a generic future directed past inextendible curve to the
point $y.$
\end{theorem}

\m In Section~\ref{mod} we develop a combinatorial method for
computing $\alk(\mathfrak S _x, \mathfrak S_y).$ This is done from the
shapes of $(W_{x,M}, W_{y,M})\subset M$ equipped with orthogonal 
to the fronts direction fields, defining their lifts to
$STM.$ This method is motivated by Arnold's~\cite{Arnold} definition
of the $J^+$-invariant of planar wave fronts. (Please, do not
confuse this $J^+$ with the causal future.)

Arnold observed that generic double points of immersed Legendrian
submanifolds in $STM=ST^*M$ correspond to the tangencies of their
cooriented projections to $M$ at which the coorienting normals to
the two immersed tangent branches point to the same direction. These
tangencies are called {\em dangerous tangencies.\/} Arnold defined
his $J^+$-invariant of a planar front by describing its increments
under passages through the dangerous self-tangencies. Thus to
compute $J^+$ one has to change the front to be ``trivial'' by a
sequence of moves that are dangerous tangencies and the
modifications corresponding to singularities of the front arising
under a generic Legendrian isotopy. Then $J^+$ of the front is the
value of $J^+$ on the trivial front plus the sum of the increments
under the dangerous tangency moves that were used. We derive a
formula for the increment of $\alk$ under the passage through the
dangerous tangency between the two fronts. (Since $\alk$ is a link
homotopy invariant, it does not change under the dangerous
self-tangency move.) When fronts are one-dimensional, our $\alk$
changes similarly to Arnold's $J^+.$

\m
Now we explain the behavior of $\alk$ under a passage through
dangerous tangency. Consider a positively oriented chart $(x_1,
\cdots, x_m)$ such that the dangerous tangency happens at the origin
where the common normal vector to the immersed branches of the two
fronts that defines their lift to $STM$ is
$-\frac{\partial}{\partial x_m}$. Locally the two fronts $W_1, W_2$
can be expressed as graphs of some functions
$$
x_m=f_i(x_1, x_2, \cdots, x_{m-1}),\, i=1, 2.
$$
Put $\sigma$ to be the number of negative eigenvalues of the Hessian
of $f_2-f_1$ at the origin. Put $\eps$ to be $+1$ if the two
oriented immersed tangent branches induce the same orientation on
the common tangent $(m-1)$-plane and put $\eps=-1$ otherwise. Put
$\alpha$ to be $+1$ (respectively $\alpha=-1$) if the
$x_m$-coordinate of the point of $W_1$ projecting to the origin in
the $(x_1, x_2, \cdots, x_{m-1})$-hyperplane after the move is
larger (respectively less) than the $x_m$-coordinate of the
corresponding point on $W_2$ after the move.

\begin{theorem}[see Theorem~\ref{hesse}]\label{dangerous}
Under a passage through a dangerous tangency $\alk$ increases by
$\eps\alpha(-1)^{\sigma}$. Recall that $\alk$ always takes values
either in $\Z$ or in $\Z_n,$ so this expression indeed makes sense.
\end{theorem}

\m We use Theorem~\ref{dangerous} to construct examples where we
can conclude that the events are causally related from the shapes of
their fronts, see Section~\ref{examplescausality}. This conclusion
can be made without the knowledge of the Lorentz metric on the
space-time, of the event points, and in many cases even without the
knowledge of topology of the globally hyperbolic space-time.

\m
In Section~\ref{s:refocussing} we discuss the refocussing phenomena,
see \defref{d:nonrefocussing}. A good property of nonrefocussing
globally hyperbolic space-times is that the map $\mu: \ss \to $
\{the space of skies\}, $x\mapsto \mathfrak S_x$ is a homeomorphism.
Low~\cite{LowNullgeodesics} introduced the concept of refocussing
spaces and noticed that a globally hyperbolic space-time with a
noncompact Cauchy surface is nonrefocussing, see
\propref{p:nonrefocussing}. We prove that a globally hyperbolic
space-time $(\ss, g)$ is nonrefocussing whenever any of its covering
space-times is, see \theoref{refocussingtheorem}. In particular, if
$\pi_1(\ss)$ is infinite, then $(\ss, g)$ is nonrefocussing.

As we discuss in Remark~\ref{Blaschke}, the question on topology of a 
refocussing space-time is related to the problems similar to the {\bf Blaschke 
conjecture\/} in Riemannian geometry.

\m
Low~\cite[Problem 7]{LowNullgeodesics} asked: ``Is there any
construction intrinsic to the space $\mathcal N$ which will enable
us to decide whether the points represented by two skies are
causally related?'' The following \theoref{problem7} gives an
affirmative answer for all nonrefocussing globally hyperbolic $(\ss,
g).$ Also, \theoref{problem7} says that a weakened version of the
Low conjecture holds for all globally hyperbolic nonrefocussing
space-times.

\begin{thm}[See Corollary~\ref{Low2true}, Definition~\ref{isotopythroughskies}]\label{problem7}
Let $(\ss, g)$ be a nonrefocussing globally hyperbolic space-time of
dimension $>2.$ Let $(x_1, x_2)$ be a pair of causally unrelated
events and let $(y_1, y_2)$ be a pair of events that do not belong
to a common null geodesic. Then the following two statements are
equivalent:
\begin{enumerate}
\item $y_1, y_2$ are causally related;
\item The link $(\mathfrak S_{y_1}, \mathfrak S_{y_2})$ is not isotopic to
$(\mathfrak S_{x_1}, \mathfrak S_{x_2})$ via an isotopy through
skies of events in $(\ss,g).$
\end{enumerate}
\end{thm}

This Theorem follows from the following more general fact that holds for all
globally hyperbolic $(\ss, g)$ and is also closely related to the
above questions.

\begin{thm}[See Theorem~\ref{Low1true}]
Let $(\ss^{m+1}, g),m>1$ be a globally hyperbolic space-time. Let $(x_1, x_2)$ 
be a pair of causally unrelated events and
let $(y_1, y_2)$ be a pair of events that do not belong to a common
null geodesic. Then the following two statements are equivalent:
\begin{enumerate}
\item $y_1, y_2$ are causally related;
\item for every pair of paths $\rho_i:[0,1]\to \ss$ 
such that $\rho_i(0)=x_i$ and $\rho_i(1)=y_i, i=1,2,$ there exists
$t\in [0,1]$ such that $\rho_1(t)$ and $\rho_2(t)$ belong to a
common null geodesic.
\end{enumerate}
\end{thm}

\section{linking and causality}\label{s:linkedcausal}

The following Theorem says that the skies of two causally unrelated
events are Legendrian unlinked. In particular, we see that for every
Cauchy surface $M$ the lifted wave fronts $\wt W_{x, M}, \wt W_{y,
M}$ of two causally unrelated events $x,y$ are Legendrian unlinked.

\begin{thm}\label{thmfrontsunlinked}
Let $(\ss, g)$ be a globally hyperbolic space-time. If $x$ and $y$
are causally unrelated events, then the pair $(\mathfrak S_x,
\mathfrak S_y)$ is Legendrian unlinked.
\end{thm}

\pp Choose a Cauchy surface $M \subset \ss$ and an $M$-proper
isometry $h: M \times \R \to \ss.$ It suffices to prove that the
lifted wave fronts $\wt W_{x,M}=\wt W_{x, M}^0$ and $\wt W_{y,M}=\wt
W_{y, M}^0$ are Legendrian unlinked in $ST^*M$.

Take $\tau_1, \tau_2\in \R$ such that $x\in M_{\tau_1}$ and $y\in
M_{\tau_2}$. Thus $h(x)=(m_1, \tau_1)$ and $h(y)=(m_2, \tau_2),$ for
some $m_1, m_2\in M.$ Without loss of generality we assume that
$\tau_1\leq \tau_2.$ There are three possible cases $\tau_1\leq
\tau_2\leq 0,$ $\tau_1\leq 0\leq \tau_2,$ and $0\leq \tau_1\leq
\tau_2.$ We prove the Theorem only for the case $\tau_1\leq 0\leq
\tau_2.$ The proof in the other two cases is similar and, in fact,
even slightly easier.

Let $S_i, i=1,2$ be a copy of $S^{m-1}$. Consider $I_1:(S_1\sqcup
S_2)\times [0, \tau_2]\to ST^*M$ defined by $I_1(s_1, t)=\wt W^t_{x,
M}(s_1)$ and $I_1(s_2, t)=\wt W^t_{y, M}(s_2),$ for $s_1\in S_1,
s_2\in S_2, t\in [0, \tau_2].$ Since $x, y$ do not lie on a common
null geodesic, we see that $I_1$ is a Legendrian isotopy between
$(\wt W_{x, M}^0, \wt W_{y, M}^0)$ and $(\wt W_{x, M}^{\tau_2}, \wt
W_{y, M}^{\tau_2})=(\wt W_{x, M}^{\tau_2}, \eps_{m_2}).$

Consider a timelike curve $\rho:[\tau_1, \tau_2]\to \ss$ given by
$\rho(t)=h(m_2, t),$ for $t\in [\tau_1, \tau_2].$ The future
directed null geodesics of the sky $\mathfrak S_x$ do not intersect
$\rho.$ Otherwise such a null geodesic followed by $\rho$ after the
intersection point is a piecewise smooth nonspacelike curve from $x$
to $y.$ This would contradict the assumption that $x$ and $y$ are
causally unrelated.

Consider $I_2:(S_1\sqcup S_2)\times [\tau_1,\tau_2]\to ST^*M$
defined by $I_2(s_1,t)=\wt W_{x, M}^t(s_1),$ $I_2(s_2,
t)=\eps_{m_2}(s_2),$ for $s_1\in S_1, s_2\in S_2, t\in [\tau_1,
\tau_2].$ Since $\rho$ does not intersect the null geodesics of the
sky $\mathfrak S_x,$ we get that $m_2\not \in \Im W_{x, M}^t$ for
all $t\in [\tau_1, \tau_2],$ and hence $I_2$ is an isotopy. Since
lifted wave fronts are Legendrian maps, we conclude that $I_2$ is a
Legendrian isotopy between $(\eps_{m_1}, \eps_{m_2})=(\wt
W^{\tau_1}_{x, M}, \eps_{m_2})$ and $(\wt W^{\tau_2}_{x, M},
\eps_{m_2}).$

Combining isotopies $I_1$ and $I_2$, we conclude that $(\wt W_{x,
M}, \wt W_{y, M})$ is Legendrian isotopic to $(\eps_{m_1},
\eps_{m_2}).$ \qed

\section{Review of the $\alk$ invariant.\/}\label{review} 
{\it In this
Section we adapt the general $\alk$ invariant constructed by us 
in~\cite{ChernovRudyakGT} to the case of linked skies.
Throughout the Section $M^m$ is a smooth connected oriented manifold of 
dimension
$m>1.$\/}

\begin{defin}[bordism group]
For a space $Y$ put $\bor_n(Y)$ {\it to be the $n$-dimensional
oriented bordism group of $Y.$\/} Recall that $\bor_n(Y)$ is the set
of the equivalence classes of (continuous) maps $g: V^n \to Y$ where
$V$ is a smooth closed oriented manifold. Here two maps $g_1: V_1
\to Y$ and $g_2: V_2 \to Y$ are equivalent if there exists a map $f:
W^{n+1}\to Y,$ where $W$ is an oriented compact smooth manifold
whose oriented boundary $\p W$ is diffeomorphic to $V_1\sqcup
(-V_2)$ and $f|_{\p W}=g_1\sqcup g_2$. Disjoint union operation
turns $\bor_n(Y)$ into an abelian group, and $\bor_n(Y)$ is
canonically isomorphic to $H_n(Y),$ for $0\le n \le 3.$ See
\cite{Rudyak, Stong, Switzer} for details.

\end{defin}

For a space $Y$, the group $\bor_0(Y)=H_0(Y)$ is the free abelian
group with the base $\pi_0(Y)$. So, every element of $\bor_0 (Y)$
can be represented as a finite formal linear combination $\sum
a_kP_k$ with $a_k\in \Z$ and $P_k\in Y$. Conversely every such
linear combination gives us an element of $\bor_0(Y)$.

\m

Put $\cs$ to be the connected component of the space of
$C^{\infty}$-mappings $S^{m-1}\to (ST^*M)^{2m-1}$ that consists of
the mappings homotopic to some (and hence to all) $\eps_v, v\in
M^m$. (Note that the mappings in $\cs$ are not assumed to be
immersions or Legendrian mappings.) Let $\cs \bul$ be the space of
pointed maps $(S^{m-1}, \star)\to (ST^*M,\star)$ such that the
corresponding maps $S^{m-1}\to ST^*M$ are in $\cs$.

\begin{lem}\label{dotsconnected}
For an oriented connected manifold $M^m, m>1,$ the space $\cs \bul$
is path connected.
\end{lem}

\pp The standard $\pi_1(M)$-action on $\cs \bul$ induces the bijection
$\cs=\cs \bul/\pi_1(M)$. Since $\cs$ is a singleton by definition,
we conclude that the above $\pi_1(M)$-action on $\cs \bul$ is
transitive. So it suffices to prove that the $\pi_1(M)$-action is
trivial.

Consider a loop $\gamma:S^1\to ST^*M$ that realizes $[\gamma]\in
\pi_1(ST^*M)$ and put $[S^{m-1}_{\star}]\in \cs \bul$ to be the
pointed homotopy class of the positively oriented fiber $S^{m-1}$ of
$\pr$ containing the base point $\star$. Consider the
$S^{m-1}$-bundle over $S^1$ induced from $\pr:ST^*M\to M$ by
$\pr\circ\gamma:S^1\to M.$ This bundle is trivial, since $\pr$ is an
oriented bundle. We choose its trivialization and obtain a bundle
map
\[
\CD
S^{m-1}\times S^1 @>>> ST^*M\\
@VVV @VVV\\
S^1 @>\pr\circ \gamma>> M.
\endCD
\]
Now we see $[\gamma][S^{m-1}_{\star}]=S^{m-1}_{\star}$ since
$\pi_1(S^1\times S^{m-1})$ acts trivially on $\pi_{m-1}(S^1\times
S^{m-1})$. Finally, $[\gamma]x=x$ for all $x\in \cs\bul,$ since the
$\pi_1(M)$-action on $\cs \bul$ is transitive. \qed

\begin{defin}[of $\B$]\label{definB}

Let $\B=\B_{\cs,\cs}$ be the space of quadruples $(\phi_1, \phi_2,
\rho_1, \rho_2)$ where $\phi_i: S^{m-1} \to ST^*M, i=1,2,$ belong to
$\cs$ and $\rho_i: \pt \to S^{m-1}$ are mappings of the
one-point-space $\pt$ such that $\phi_1\rho_1=\phi_2\rho_2$.
Clearly, $\B$ can be regarded as a subset of $\cs\times \cs \times
S^{m-1}\times S^{m-1}$, and we equip $\B$ with the subspace
topology.
\end{defin}

\begin{lem}\label{Bconnected}
For an oriented connected manifold $M^m, m>1,$ the space $\B$ is
path connected. Thus the augmentation $\aug: \bor_0(\B)\to
\bor_0(\pt)=\Z$ induced by the map $\B\to \pt$ is an isomorphism.
\end{lem}

\pp Our~\cite[Theorem 7.4]{ChernovRudyakGT} says that $\pi_0(\B)$
is the quotient of $\pi_0(\cs \bul)\times \pi_0(\cs \bul)$ by a
certain right action of $\pi_1(ST^*M)$ and a certain left action of
$\pi_1(S^{m-1})\times \pi_1(S^{m-1}).$ Now the result follows from
\lemref{dotsconnected}. \qed

\begin{defin}[of the $\mu$-pairing]\label{mupairing}
Let $\alpha_1:F_1^i\to \cs$ be a map representing $[\alpha_1]\in
\bor _i(\cs)$ and let $\alpha_2:F_2^j\to \cs$ be a map representing
$[\alpha_2]\in \bor_{j}(\cs)$. Let $\wt \alpha_l:F_l\times
S^{m-1}\to ST^*M$, $l=1,2$, be the adjoint maps i.e.~maps such that
$\wt \alpha_{l}(f,s)=(\alpha_{l}(f))(s)$. Following standard
arguments we can assume that $\wt \alpha_1$ and $\wt \alpha_2$ are
transverse. Consider the pullback diagram
\begin{equation}\label{bordism}
\CD
V@>k_1>> F_1\times S^{m-1}\\
@VVk_2V @VV\wt \alpha _1 V\\
F_2\times S^{m-1} @>\wt \alpha _2 >> ST^*M\\
\endCD
\end{equation}
of the maps $\wt \alpha_i$, $i=1,2$.

If $\wt \alpha_1$ and $\wt \alpha_2$ are transverse, then
$V=\{(f_1, s_1, f_2, s_2)| \wt \alpha_1(f_1, s_1)=\wt \alpha_2(f_2,
s_2)\}$ is a smooth closed
$\bigl(i+j+(m-1)+(m-1)-(2m-1)\bigr)=(i+j-1)$-dimensional submanifold
of $F_1\times S^{m-1}\times F_2\times S^{m-1}$. It is identified
with the transverse preimage of the diagonal in $ST^*M\times ST^*M$
under the map $\wt\alpha_1\times \wt \alpha_2:(F_1\times
S^{m-1})\times (F_2\times S^{m-1})\to ST^*M\times ST^*M,$ and hence
$V$ is canonically oriented.

Put $ \mu (\wt \alpha_1, \wt \alpha_2):V\to \mathcal B$ to be the
map sending $(f_1, s_1, f_2, s_2)\in V$ to $(\alpha_1(f_1),
\alpha_2(f_2), \rho_{s_1}, \rho_{s_2}),$ where
$\rho_{s_l}(\pt)=s_l\in S^{m-1}, l=1,2.$ As we showed
in~\cite[Theorem 2.2]{ChernovRudyakGT} the above construction yields
a well-defined pairing
\begin{equation}\label{muequation}
\begin{aligned}
\mu=\mu_{ij}: \bor _i(\cs )&\otimes \bor_j(\cs )\to \bor_{i+j-1} (\B
),\\
\mu \left([\alpha_1], [\alpha_2]\right)&=[V, \mu (\wt \alpha_1,
\wt \alpha_2)].
\end{aligned}
\end{equation}
\end{defin}

\begin{defin}\label{sigma}
Put $\Sigma$ to be the {\em discriminant\/} in $\cs \times \cs$,
i.e.~the subspace of $\cs\times \cs $ that consists of pairs $(f_1,
f_2)$ such that there exist $s_1, s_2\in S^{m-1}$ with $f_1(s_1)=f_2
(s_2)$. (We do not include into $\Sigma$ the maps that are singular
in the common sense but do not involve double points between $f_1
(S^{m-1})$ and $f_2(S^{m-1})$.)

Put $\Sigma_0$ to be the subset (stratum) of $\Sigma$ consisting of
all the pairs $(f_1, f_2)$ for which there exists precisely one pair
($s_1, s_2)$ of points $s_1, s_2\in S^{m-1}$ such that
$f_1(s_1)=f_2(s_2)$ and moreover
\begin{description}
\item[a] $s_i$ is a regular point of $f_i, i=1,2$;
\item[b] $(df_1)(T_{s_1}S^{m-1})\cap (df_2)(T_{s_2}S^{m-1})=0$.
\end{description}
\end{defin}

Note that there is a canonical map of $\Sigma_0$ into $\B.$ Namely,
we assign the commutative diagram $(f_1, f_2, \rho_{s_1},
\rho_{s_2})$ with $\rho_{s_i}:\pt \to s_i\in S^{m-1}, i=1,2,$ to the
pair $(f_1, f_2)\in \Sigma_0$ with $f_1(s_1)=f_2(s_2)$.

\begin{defin}
[of the sign of the crossing of $\Sigma_0$ and of a generic path in
$\cs\times \cs$]\label{signsigma0} Consider a singular link $(f_1,
f_2)\in \Sigma_0.$ The double point $z=f_1(s_1)=f_2(s_2)$ of it can
be resolved in two (essentially different) ways. To a resolution
$(\ov f_1, \ov f_2)$ (that is a $C^{\infty}$-small deformation of
$(f_1,f_2)$) we associate the vector $\mathbf w\in T_z ST^*M$ that
in a chart has the same direction as the vector from $\ov f_1(s_1)$
to $\ov f_2(s_2).$ We say that the resolution $(\ov f_1, \ov f_2)$
is {\em generic} if $\span\{(d\ov f_1)(T_{s_1}S^{m-1}), \mathbf w,
(d\ov f_2)(T_{s_2}S^{m-1})\}=T_{f_1(s_1)}M$.

Let $\mathfrak r_i, i=1,2$ be the positive $(m-1)$-frames in
$T_{s_i}S^{m-1}.$ Take a generic resolution of $(f_1, f_1)$ and
consider the $(2m-1)$-frame $$\{ df_1(\mathfrak r_1), {\mathbf w},
df_2(\mathfrak r_2) \} \subset T_z(ST^*M).$$ We say that the
resolution of the singular link is {\em positive\/} if this
$(2m-1)$-frame gives the canonical orientation of $ST^*M$, and we
say that the resolution is {\em negative,\/} otherwise. One checks
that the sign of the resolution does not depend on the choice of the
chart used to define $\mathbf w$.

Let $\gamma(t)$ be a path that intersects $\Sigma$ at one point
$\gamma(t_0)\in \Sigma_0$. We say that $\gamma$ intersects $\Sigma$
{\em transversally} at $\gamma(t_0)$ if $\gamma(t_0)\in \Sigma_0$
and if the resolution $(\ov f_1, \ov f_2)=\gamma(t)$ is generic for
$t$ close to $t_0$ and different from $t_0$ We put the {\em sign of
the transverse intersection of $\Sigma_0$ by $\gamma$\/} to be the
sign of the singular link resolution induced by $\gamma$, and denote
this sign by $\sigma(\gamma, t_0)=\pm 1$. Clearly if we traverse
the path $\gamma$ in the opposite direction, then the sign of the
intersection changes.

We say that a path $\gamma$ in $\cs\times \cs$ is {\em generic\/} if
it intersects $\Sigma$ at a finite number of times and these
intersections are transverse. We will also use the term ``generic
link homotopy'' for a generic path.
\end{defin}

\begin{defin}[of $\A(M)$ and of the $\alk$ invariant]\label{definalk}
Define the indeterminacy subgroup $\Indet$ of $\bor_0(\B)$ to be
the subgroup generated by the images of $\mu_{0,1}$ and $\mu_{1,0}.$
Put $\A=\A(M)=\A(ST^*M)$ to be the quotient group
$\bor_0(\B)/\Indet$ and put $q:\bor_0(\B)\to \A$ to be the quotient
homomorphism.

Our~\cite[Theorem 3.9]{ChernovRudyakGT} when applied to this work
setup says that there exists a function
\begin{equation}\label{eq:alk}
\alk: \cs\times \cs \setminus \Sigma \to \A(M)
\end{equation}
such that:

\begin{description}
\item[a] $\alk$ is constant on path connected components of
$\cs \times \cs \setminus \Sigma ;$
\item[b] if $\gamma:[a,b]\to \cs\times \cs$ is a generic
path such that $\gamma(a), \gamma(b)\not \in \Sigma$ and $t_i, i\in
I,$ are the moments when $\gamma (t_i)\in \Sigma$ $($and hence
$\gamma (t_i)\in \Sigma _0$ by the definition of the generic
path$)$, then
$$
\alk(\gamma(b))- \alk(\gamma(a))=q\Bigl (\sum _{i\in
I}\sigma(\gamma, t_i)\gamma (t_i)\Bigr ) \in \A(M).
$$
\end{description}

We showed that such $\alk$ is unique up to an additive constant.
{\em In this paper we normalize $\alk$ by the condition that
$\alk(\eps_u, \eps_v)=0$ for any two distinct $u,v\in M.$ \/}

We proved~\cite[Corollary 7.5]{ChernovRudyakGT} that for every
$\alpha\in \A(M)$ there exists a nonsingular link $(f_1, f_2)\in
\cs\times\cs\setminus \Sigma$ with $\alk(f_1, f_2)=\alpha.$ Thus
$\A(M)$ is indeed the group of values of the $\alk$-invariant.
\end{defin}

\begin{defin}[of the affine linking number of a pair of skies]\label{alkskies}
Let $(\ss, g)$ be a globally hyperbolic space-time and let $x,y\in
\ss$ be events that do not lie on a common null geodesic. Choose a
Cauchy surface $M\subset \ss.$ Since $x,y$ do not lie on a common
null geodesic the pair of lifted wave fronts $(\wt W_{x,M}, \wt
W_{y,M})$ is a point in $\cs\times \cs\setminus \Sigma$, and we put
$$
\alk_M(\mathfrak S_x, \mathfrak S_y)=\alk(\wt W_{x, M}, \wt W_{y,
M})\in \A(M)
$$
where the $\alk$ at the right-hand side means the function
\eqref{eq:alk}. \theoref{welldefined} below states that the value
$\alk_M(\mathfrak S_x, \mathfrak S_y)\in \A(M)$ does not depend on
the Cauchy surface $M$. Thus we can and shall define $\alk
(\mathfrak S_x, \mathfrak S_y):=\alk_M(\mathfrak S_x, \mathfrak
S_y),$ for any choice of $M$.
\end{defin}

\begin{thm}\label{welldefined}
Let $x, y$ be two events in a globally hyperbolic space-time $(\ss,
g)$ that do not lie on a common null geodesic. Let $M$ and $N$ be
two Cauchy surfaces in $\ss$, and let $h: M \times \R\to \ss$ and
$h': N \times \R\to \ss$ be $M$-proper and $N$-proper isometries,
respectively. Then for any $t, \tau\in \R$ the following holds:
\begin{description}
\item[1] $(\wt W_{x, M}^t, \wt W_{y, M}^t)$ is unlinked 
$($respectively Legendrian unlinked$)$ in $ST^*M$ if and only if $(\wt
W_{x, N}^{\tau}, \wt W_{y, N}^{\tau})$ is unlinked $($respectively
Legendrian unlinked$)$ in $ST^*N$.
\item[2] $\alk(\wt W_{x, M}^t, \wt W_{y, M}^t)=\alk(\wt W_{x, N}^{\tau}, \wt
W_{y, N}^{\tau})\in \A(M)=\A(N).$
\end{description}
In particular, the value
$$\alk(\mathfrak S_x,\mathfrak S_y)=\alk(\wt
W_{x, M}, \wt W_{y, M})=\alk(\wt W_{x, M}^0, \wt W_{y, M}^0)\in
\A(M)
$$
and the notion of the skies $(\mathfrak S_x, \mathfrak S_y)$
being unlinked $($respectively Legendrian unlinked$)$ are well
defined.
\end{thm}

\pp Clearly, the link $(\wt W_{x, M}^t, \wt W_{y, M}^t)$ is
Legendrian isotopic to the link $(\wt W_{x, M}^0, \wt W_{y, M}^0)$
in $ST^*M$, and the similar fact is true for the link $(\wt W_{x,
N}^{\tau}, \wt W_{y, N}^{\tau})$ in $ST^*N.$ Hence, to prove
Statement $1,$ it suffices to show that $(\wt W_{x, M}^0, \wt W_{y,
M}^0)$ is (Legendrian) unlinked if and only if $(\wt W_{x, N}^{0},
\wt W_{y, N}^{0})$ is.

Since $\alk$ invariant does not change under link isotopy, to
prove Statement $2$ it suffices to show that $\alk(\wt W_{x, M}^{0},
\wt W_{y, M}^{0})=\alk(\wt W_{x, N}^{0}, \wt W_{y, N}^{0}).$

{\em We prove the ``Legendrian unlinked'' part of statement $1.$ The
proof of the ``unlinked'' part is obtained by omitting the word
``Legendrian'' everywhere in the proof.\/} Assume that the link
$(\wt W_{x, M}^0, \wt W_{y, M}^0)$ is Legendrian unlinked in
$ST^*M$. Let $S_i, i=1,2,$ be a copy of $S^{m-1}$. Choose a
Legendrian isotopy $I_{t}:S_1\sqcup S_2\to ST^*M, t\in [0,1],$ such
that $I_0= \wt W_{x, M}^0\sqcup \wt W_{y, M}^0$ and
$I_1=\eps_{u}\sqcup \eps_{v}$ for some $u\neq v\in M.$

Put $\wt x=i(u)$ and $\wt y=i(v)$ where $i:M\to X$ is the inclusion.
Then $\wt x$ and $\wt y$ belong to the same Cauchy surface $M$ and
hence are causally unrelated. Thus by
Theorem~\ref{thmfrontsunlinked} the link $(\wt W_{\wt x, N}^{0}, \wt
W_{\wt y, N}^{0})$ is Legendrian unlinked in $ST^*N$.

Let $f^M_N:ST^*M\to ST^*N$ be the contactomorphism \eqref{eq:psi}.
Now, $f^M_N\circ I_t, t\in [0,1],$ is a Legendrian isotopy that
deforms $(\wt W_{x, N}^{0}, \wt W_{y, N}^{0})$ to the Legendrian
trivial link $(\wt W_{\wt x, N}^{0}, \wt W_{\wt y, N}^{0})$ in
$ST^*N.$

\m {\em Now we prove statement $2$ of the Theorem.\/} Let
$\gamma:[0,1]\to \cs\times \cs$ be a generic smooth link homotopy
such that $\gamma(1)= (\wt W_{x, M}^0, \wt W_{y, M}^0)$ and
$\gamma(0)=(\eps_{u}, \eps_{v})$ for some $u \neq v\in M.$ Let $t_i,
i=1, \ldots, k,$ be the time moments when $\gamma$ crosses
$\Sigma_0\subset \Sigma$ and let $\sigma(\gamma, t_i)=\pm 1$ be the
signs of these crossings, see~\ref{signsigma0}.
Lemma~\ref{Bconnected} says that $\mathcal \B$ is connected and
hence $\alk(\wt W_{x, M}^0, \wt W_{y, M}^0)=q\bigl(\sum_{i\in I}
\sigma(\gamma, t_i)\bigr),$ for the map $q:\bor_0(\B)=\Z\to \A(M).$

The smooth link homotopy $f^M_N\circ \gamma(t), t\in [0,1],$ deforms
the trivial link $(\wt W_{\wt x, N}^{0}, \wt W_{\wt y, N}^{0})$ in
$ST^*N$ to the link $(\wt W_{x, N}^{0}, \wt W_{y, N}^{0})$ in
$ST^*N.$ Clearly the link homotopy $f^M_N\circ \gamma(t), t\in
[0,1],$ is generic and it crosses $\Sigma_0\subset \Sigma$ at the
same time moments $t_i, i\in I,$ as $\gamma(t).$ Moreover since
$f^M_N$ is orientation preserving, we conclude that
$\sigma(f^M_N\circ \gamma, t_i)=\sigma(\gamma, t_i).$ Since $M$ and
$N$ are diffeomorphic, they are homotopy equivalent. Hence
$\A(M)=\A(N)$ and the maps $q:\Z\to \A(M)$ and $q:\Z\to \A(N)$ are
the same. Thus
\begin{equation*}
\begin{aligned}
\alk (\wt W_{x, N}^{0}, \wt W_{y, N}^{0})-\alk (\wt W_{\wt x,
N}^{0}, \wt W_{\wt y, N}^{0})&=q\bigl( \sum_{i\in I}
\sigma(f^M_N\circ \gamma, t_i)\bigr)=
\\
&q\bigl (\sum_{i\in I} \sigma(\gamma, t_i)\bigr)=
\alk (\wt W_{x,
M}^0, \wt W_{y, M}^0).
\end{aligned}
\end{equation*}
Since the link $(\wt W_{\wt x, N}^{0}, \wt W_{\wt y, N}^{0})$ is
trivial, $\alk (\wt W_{\wt x, N}^{0}, \wt W_{\wt y, N}^{0})=0$ and
hence $\alk (\wt W_{x, N}^{0}, \wt W_{y, N}^{0})=\alk (\wt W_{x,
M}^0, \wt W_{y, M}^0).$ \qed

\section{Computation of the group $\A(M)$.\/}\label{ComputationofA}
{\it In this section $M^m, m>1$ is a smooth connected oriented manifold.\/}

\begin{defin}\label{spec}
Given a map $\wt \alpha: S^1\times S^{m-1}\rightarrow ST^*M$, we say
that $\wt \alpha$ is {\it special} if $\wt \alpha\big |_{\star
\times S^{m-1}}$ has the form $\eps_v$ for some $v\in M$,
see~\eqref{epspoint}. Here $\star \in S^1$ is the base point.

We define $\wt \Indet\subset \Z$ to be the subgroup of $\Z$
generated by
$$
\bigl \{ \wt \alpha_*[(S^{m-1}\times S^1)]\bullet [S^{m-1}]\in
\Z=H_0(ST^*M)\text{ such that }\wt\alpha \text{ is special}\bigr\}.
$$
Here $\bullet$ is the intersection pairing of homology classes and
$[S^{m-1}]\in H_{m-1}(ST^*M)$ is the homology class of a positively
oriented fiber of $\pr:ST^*M\to M.$
\end{defin}

\begin{lem}\label{computeA}
The isomorphism $\aug:\bor_0(\B)\to \Z$ from
Lemma~$\ref{Bconnected}$ maps $\Indet$ onto $\wt\Indet$.
\end{lem}

\pp The subgroup $\Indet$ of $\bor_0(\B)$ was defined as the
subgroup generated by the images of $\mu_{0,1}$ and $\mu_{1,0}$
where $\mu_{i,j}:\bor_i(\cs)\otimes \bor_j(\cs)\to
\bor_{i+j-1}(\B),$ see~\eqref{mupairing}. In particular, the images
of $\mu_{0,1}$ and of $\mu_{1,0}$ are subgroups of $\bor_0(\B)=\Z.$
It is easy to see that $\mu_{0,1}(\alpha_0, \alpha_1)=\pm
\mu_{1,0}(\alpha_1, \alpha_0),$ for any $\alpha_0\in \bor_0(\cs),
\alpha_1\in \bor_1(\cs),$ where the sign depends on the dimension of
$M.$ Thus $\Im(\mu_{0,1})=\Indet=\Im(\mu_{1,0}).$

Take $\beta:\pt \to \eps_u\in \cs$ and $\alpha:S^1\to \cs.$ Without
loss of generality we can assume (deforming $\alpha$ if necessary)
that $\alpha(\star)=\eps_v$, for some $v\in M$ and that the adjoint
$\wt \alpha:S^1\times S^{m-1}\to ST^*M$ of $\alpha$ is transverse to
$\beta=\eps_u.$ This homotopy does not change $[S^1, \alpha]\in
\bor_1(\cs).$ Now the adjoint $\wt\alpha$ of $\alpha$ is a special
map.

From~\ref{mupairing} one verifies that the bordism class
$\mu_{1,0}\bigl ([S^1, \alpha], [\pt, \beta]\bigr)\in \bor_0(\B)$ is
represented by the set of intersection points of the maps
$\wt\alpha:S^1\times S^{m-1}\to ST^*M$ and $\wt \beta:\pt\times
S^{m-1}=S^{m-1}\to ST^*M.$ The signs at these intersection points
are equal to the signs obtained from the definition of the
intersection number of two transverse oriented submanifolds of
complimentary dimensions.

Since $\bor_0(\B)=\Z[\pi_0(\B)]=\Z$ we conclude that $\mu_{1,0}\bigl
([S^1, \alpha], [\pt, \beta]\bigr)\in \Z$ equals to the intersection
number $\wt \alpha_* [S^1\times S^{m-1}]\bullet \wt
\beta_*[S^{m-1}]\in \Z.$

Recall that $\bor_i(Y)=H_i(Y)$ for $0\leq i\leq 3$ and all spaces
$Y.$ In particular, every class in $\bor_0(\B)$ is parameterized by
a collection of oriented points and every class in $\bor_1(\B)$ is
parameterized by a collection of oriented circles. Thus
$\Indet\subset \wt \Indet\subset \Z.$

On the other hand every smooth special $\wt \alpha:S^1\times
S^{m-1}\to ST^*M$ is the adjoint of a certain map $\alpha:S^1\to
\cs.$ Thus $\wt \Indet\subset \Indet\subset \Z.$ \qed

\begin{rem} For future needs, it is convenient to regard the
augmentation isomorphism as the identification $\bor_0(\B)=\Z$. For
example, we can treat \lemref{computeA} as the equality $\wt \Indet
=\Indet$.
\end{rem}

\begin{defin}
Given oriented $m$-dimensional manifolds $N^m, M^m,$ and a
continuous map $\beta: N\rightarrow ST^*M$, we define $d(\beta)$ to
be the degree of the map
\begin{equation*}
\pr\circ\beta: N^m\rightarrow M^m.
\end{equation*}
If one of $N,M$ is not a closed manifold, then we put $d(\beta)$ and
the degree of $\pr\circ \beta$ to be zero.
\end{defin}

\begin{lem}\label{degreesphere}
A number $i\in \Z$ equals to $d(\wt\alpha)$ for some special $\wt
\alpha:S^1\times S^{m-1}\rightarrow ST^*M$ if and only if $i$ equals
to $d(\beta)$ for some $\beta:S^m\to ST^*M.$ In particular, the set
$\bigl\{ d(\wt \alpha) \text{ such that } \wt\alpha:S^1\times
S^{m-1}\to ST^*M\text{ is special}\bigr\}$ is a subgroup of $\Z$
that is the image of the homomorphism $\pi_m(ST^*M)\to \Z$ sending
the class of a map $\beta:S^m\to ST^*M$ to $d(\beta).$
\end{lem}

\pp We regard $S^{m-1}$ and $S^m$ as pointed spaces.

{\em Assume that $i=d(\wt \alpha)$ for a special $\wt \alpha.$ We
show that $i=d(\beta)$ for some $\beta:S^m\to ST^*M.$\/} Consider a
map $\ov \alpha: S^1\times S^{m-1} \rightarrow ST^*M$ such that:
\begin{description}
\item[1]$\ov \alpha\big |_{\star \times S^{m-1}}=\wt \alpha\big |_
{\star \times S^{m-1}}$,
\item[2] $\ov \alpha \big |_{S^1\times \star}=\wt \alpha \big |_
{S^1 \times \star}$,

\item[3] $\ov \alpha |_{t \times S^{m-1}}=\eps_{\ov
\alpha (t,\star)},$ for all $t\in S^1$.
\end{description}
We regard $S^1\times S^{m-1}$ as the $CW$-complex with four cells
$e^0, e^1,e^{m-1},e^m$, $\dim e^k=k$. It is easy to see that the
maps $\ov \alpha$ and $\wt \alpha$ coincide on the $(m-1)$-skeleton.
Thus, the maps $\alpha$ and $\wt \alpha$ (restricted to the
$m$-cell) together yield a map $\beta: S^m\rightarrow ST^*M$.
Clearly $d(\ov \alpha)=0$, and therefore $i=d(\beta)=d(\wt \alpha)$.

{\em Assume that $i=d(\beta)$ for some $\beta:S^m\to ST^*M.$ Let us
show that $i=d(\wt \alpha)$ for some special $\wt \alpha.$\/} Let
$\pi:S^1\times S^{m-1}\to S^{m-1}$ be the projection. Choose $v\in
M$ and put $\wh \alpha=\eps_v\circ \pi:S^1\times S^{m-1}\to ST^*M.$
Consider the maps $S^1\times S^{m-1}\to ST^*M$ that coincide with
$\wh \alpha$ on the $(m-1)$-skeleton. Up to the homotopy fixed on
the $(m-1)$-skeleton they are classified by $\pi_m(ST^*M).$ Consider
such a map $\wt \alpha:S^1\times S^{m-1}\to ST^*M$ that corresponds
to $\beta\in \pi_m(ST^*M).$ Since $d(\wh\alpha)=0$ we get that
$i=d(\beta)=d(\wt \alpha).$ \qed

\begin{prop}\label{propindetdegree}
$\aug(\Indet)=\wt \Indet\subset \Z$ is the subgroup
$\bigl\{d(\beta)|\beta: S^m\to ST^*M\bigr\}\subset \Z$.
\end{prop}

\pp By Lemma~\ref{computeA} $\aug(\Indet)=\wt \Indet\subset \Z.$
Because of \lemref{degreesphere} it suffices to show that $\wt
\Indet\subset \Z$ is the subgroup $S=\bigl\{ d(\wt \alpha) \text{
such that } \wt\alpha:S^1\times S^{m-1}\to ST^*M\text{ is
special}\bigr\}\subset \Z.$ Intersection number $\wt
\alpha_*[(S^{m-1}\times S^1)]\bullet [S^{m-1}]\in \Z=H_0(ST^*M)$ and
$d(\wt\alpha)$ do not change if we substitute a special
$\wt\alpha:S^1\times S^{m-1}\to ST^*M$ by a homotopic one. Thus in
the Definition~\ref{spec} of the generating set of $\wt \Indet$ and
in the description of $S$, it suffices to consider only special $\wt
\alpha$ that have some fixed $v\in M$ as a regular value of $\pr
\circ \wt\alpha.$

Such $\wt \alpha$ and $\eps_v, v\in M,$ are transverse at all
intersection points, and $\wt\alpha^{-1}\bigl (\Im (\wt\alpha)\cap
\Im (\eps_v)\bigr)=\bigl ( \pr\circ \wt\alpha\bigr)^{-1}(v).$
Comparing the orientations of the points in these two preimage sets
we get that $\wt \alpha_*[(S^{m-1}\times S^1)]\bullet
[S^{m-1}]=d(\wt\alpha)\in \Z.$ Thus, $S$ is the subgroup $\wt
\Indet\subset \Z$. \qed

\begin{lem}\label{oddrationalhomologysphere}
We have $\A(M)=\Z,$ unless $M$ is a closed manifold that is an
odd-dimensional rational homology sphere with finite $\pi_1(M).$ In
greater detail, if there exists a map $\beta: S^m\rightarrow ST^*M$
with $d (\beta)\neq 0,$ then $M^m$ is a closed manifold that is an
odd-dimensional rational homology sphere with finite $\pi_1(M).$
\end{lem}

\pp Set~$f=\pr\circ\beta:S^{m}\to M$ and $d=d(\beta)$. Since
$d(\beta)\ne 0,$ $M$ is closed. Let $f_!:H_*(M) \to H_*(S^{m})$ be
the transfer map, see e.g. \cite[V.2.12]{Rudyak}. Since
$f_*(f^*y\cap x)=y\cap f_*x,$ for all $x\in H_*(S^{m}) $ and $y\in
H^*(M)$, we conclude that $f_*f_!(z)=dz$, for all $z\in H_*(M)$. In
particular, since $H_i(S^m)=0,$ for $0<i<m$, we conclude that
$dH_i(M) =0,$ for $0<i<m$. Thus $H_i(M;\Q)=0,$ for $0<i<m,$ and $M$
is a rational homology sphere.

Let $e\in H^m(M)=\Z$ be the Euler class of the bundle $\pr$. Since
$f$ passes through $ST^*M$, we conclude that $f^*(e)=0$. On the
other hand, the map
$$
f^*: \Z=H^m(M^m)\to H^m(S^m)=\Z
$$
is the multiplication by $d$ with $d\ne 0$. Hence $e=0$ and
$0=e=\chi(M)=1+(-1)^m,$ where $\chi(M)$ is the Euler characteristic
of $M$. Thus $m=\dim M$ is odd.

Finally since $m>1,$ every map $f=\pr\circ \beta:S^m \to M$ passes
through the universal covering of $M$, and so $d(\beta)=0$ whenever
the fundamental group of $M$ is infinite.

\qed

\begin{defin}
We put $\deg: \pi_m(M^m) \to \Z$ to be the degree homomorphism, i.e.
the homomorphism that assigns the degree $\deg f$ to the homotopy
class of a map $f: S^{m} \to M^m$. (In fact, it coincides with the
Hurewicz homomorphism $h: \pi_m(M^m) \to H_m(M^m)$ for $M$ closed
and is zero for $M$ non-closed.)
\end{defin}

\begin{prop}\label{oddcase}
If $M$ is an odd dimensional closed manifold, then
$\Indet=\Im(\deg).$
\end{prop}

\pp Proposition~\ref{propindetdegree} implies that $\Indet \subset
\Im (\deg).$ Since $M$ is an oriented odd dimensional manifold, the
projection $\pr: ST^*M \to M$ has a section $s: M \to ST^*M,
\pr\circ s=1_M$. For closed $M$ this follows, since the Euler
characteristic of $M$ is zero and it equals to the Euler class of
$T^*M\to M.$ For non-closed oriented $M$ the bundle $T^*M\to M$ has
a section regardless of the dimension of $M.$

So, every $f:S^m\to M^m$ can be written as $f=\pr\circ \beta$ with
$\beta=sf:S^m\to ST^*M.$ Since $\deg(f)=d(\beta),$ we have $\Im
(\deg) \subset \Indet.$ Hence $\Indet=\Im (\deg).$ \qed

\m Combining Definition~\ref{definalk}, \lemref{Bconnected},
\lemref{oddrationalhomologysphere}, \propref{oddcase} and the
results of the general theory of affine linking invariants reviewed
in Section~\ref{review}, we get the following result.

\begin{thm}\label{Main}
Let $M^m, m>1,$ be a smooth connected oriented manifold. If $M$ is
not a closed manifold that is an odd dimensional rational homology
sphere with finite $\pi_1(M),$ then $\A(M)=\Z$ and the homomorphism
$q:\bor_0(\B)=\Z\to \Z$ is the identity isomorphism. Otherwise
$\A(M)=\Z/(\Im (\deg: \pi_m(M^m)\to \Z))$ and $q: \bor_0(\B)=\Z\to
\A(M)$ is the quotient homomorphism. The affine linking number
invariant $\alk$ of two component links in $ST^*M$ with components
homotopic to a positive fiber $\eps_v:S^{m-1}\to ST^*M$ of
$\pr:ST^*M\to M^m$ is a link invariant such that
\begin{enumerate}
\item $\alk$ increases by $q(1)\in \A(M)$\ $($respectively by $q(-1)\in \A(M))$
under a positive $($respectively negative$)$ transverse crossing
of $\Sigma_0$, i.e. under homotopy of the link that involves
exactly one positive $($respectively negative$)$ passage through a
transverse double point between the two link components;
\item $\alk$ is invariant under Milnor's~\cite{Milnor} link homotopy that allows
each link component to cross itself, but does not allow different
components to cross.
\end{enumerate}
This $\alk$ is uniquely defined by the normalization that it is zero
on links consisting of the positive $S^{m-1}$-fibers over two
different points of $M.$ It is a universal order $\leq 1$
Vassiliev-Goussarov link homotopy invariant of such two component
links in $ST^*M.$\qed
\end{thm}

\begin{rem}\label{alksymmetry}
One verifies that the sign of crossing of $\Sigma_0$ does not depend
on the order of link components if $m=\dim M$ is even. If $m$ is
odd, then the sign of the crossing of $\Sigma_0$ gets reversed if
one changes the order of the link components. Thus $\alk(f_1,
f_2)=(-1)^{m}\alk(f_2, f_1),$ for all $m$.
\end{rem}

The group $\A(M)$ of the values of the $\alk$-invariant appears to
be quite nontrivial even in the case when $M$ is an odd-dimensional
rational homology with finite $\pi_1(M).$

\begin{prop}\label{usefulprop}
Let $M^m, m>1,$ be a smooth connected oriented manifold. Assume
moreover that $M$ is a closed manifold that is an odd-dimensional
rational homology sphere with finite $\pi_1(M).$ $($This is the only
case when Theorem~$\ref{Main}$ does not say that $\A(M)=\Z.)$ Then
the following statements hold:
\par{\rm (i)} If $\pi_1(M)$ is a finite group of order
$k$, then $\A(M)=\Z/m\Z$ where $k$ divides $m$ {\rm (}the case $m=0$
i.e. $\A(M)=\Z$ is also possible{\rm )}.
\par{\rm (ii)} If $\A(M)=0$,
then $M$ is homeomorphic to a sphere.
\end{prop}

\pp (i) This follows because every map $S^m \to M^m$ passes through
the universal covering map $p: \wt M \to M$ which is of degree $k$.

\par(ii) If $\A(M)=0,$ then there exists a map $S^m \to M^m$
of degree 1. Since every map of degree 1 of connected closed
oriented manifolds induces an epimorphism of fundamental groups and
homology groups, we conclude that $M$ is a homotopy sphere.
Poincar\'e conjecture proved in the works of Smale~\cite{Smale} for
$m\geq 5$, Freedman~\cite{Freedman} for $m=4$, and
Perelman~\cite{Perelman1},~\cite{Perelman2} for $m=3,$ implies that
$M$ is homeomorphic to a sphere. \qed

\section{Computing $\alk$ when one of the linked spheres is 
a fiber of the spherical cotangent bundle} 

{\it In this section $M$ is a smooth connected oriented manifold of
dimension $m>1.$\/}

\begin{defin}\label{df:positive}
Let $f: U^k \to V^k$ be a smooth map of oriented manifolds, and let
$v$ be a regular value of $f$. A point $u\in f^{-1}(v)$ is called
{\em positive\/} (respectively {\em negative)\/}) if a restriction
of $f$ to a small neighborhood of $u$ is orientation preserving
(respectively orientation reversing).
\end{defin}

\begin{prop}\label{prop:constant}
Suppose that the map $f$ as in $\ref{df:positive}$ is an immersion
and that $U$ is connected. Then all the points of $U$ have the same
sign, i.e. either all of them are positive or all of them are
negative.
\end{prop}

\pp This follows since the set of all positive points is open, as
well as the set of all negative points. \qed

\m

Let $F:S^{m-1}\times [a,b]\to ST^*M$ be a smooth map such that
$F(S^{m-1}\times t)\in \cs,$ for some and then for all $t\in[a,b]$.
Let $v\in M$ be a regular value of $G=\pr \circ F:S^{m-1}\times
[a,b]\to ST^*M\to M^m$ such that $G^{-1}(v)\subset S^{m-1}\times
(a,b)$. Let $n_+(v, F)$ (respectively $n_-(v, F)$) be the number of
positive (respectively negative) points in $G^{-1}(v).$ Recall that
$\bor_0(\B)=\Z$ by Lemma~\ref{Bconnected} and that
$q:\bor_0(\B)=\Z\to \Z/\Indet=\A(M)$ is the quotient homomorphism.

\begin{lem}\label{helpful}
We have the equality $$\alk( F|_{S^{m-1}\times b}, \eps _v)-\alk(
F|_{S^{m-1}\times a}, \eps _v)=q(n_-(v, G)-n_+(v, G)).$$
\end{lem}

\pp Using a $C^{\infty}$-small perturbation of $F$ if necessary we
can and shall assume that the $[a,b]$-coordinates of all the points
in $G^{-1}(v)$ are all different. Consider the link homotopy
$H:[a,b]\to \cs\times \cs$ such that $H(t)=(F|_{S^{m-1}\times t},
\eps_v).$

Let $(s_i, t_i)\in S^{m-1}\times [a,b], i=1, \ldots, k,$ be the
points of $G^{-1}(v).$ Clearly the crossings of $\Sigma$ under the
link homotopy $H$ happen exactly at time moments $t_i.$ Since all
the values $t_i, i=1, \ldots, k,$ are distinct and $v$ is a regular
value of $G,$ we conclude that the crossings of $\Sigma$ happen
inside of $\Sigma_0\subset \Sigma$ and these crossings are
transverse. From the definition of the sign of the crossing of
$\Sigma_0$ we get that sign $\sigma(H, t_i)$ is equal to $+1$
(respectively $-1$) exactly when $(s_i, t_i)$ is a negative
(respectively positive) point of $G^{-1}(v).$

By Lemma~\ref{Bconnected}, $\bor_0(\B)=\bor_0(\pt)=\Z.$ Thus by the
definition~\ref{definalk}, $\alk$ is a link homotopy invariant that
increases by $q ( \sigma(H, t_i))$ under the crossing of $\Sigma_0$
by the link homotopy $H$ that happens at time $t_i$. \qed

\section{$\alk$-invariant and causality}\label{s:alk-causality}

Let $(\ss^{m+1}, g), m>1$ be a globally hyperbolic space-time with a Cauchy 
surface $M$. Let $x, y\in \ss$ be two events that do not lie on a common null
geodesic. If $\alk(\mathfrak S_x, \mathfrak S_y)\neq 0\in \A(M),$
then the events $x, y$ are causally related by
\theoref{thmfrontsunlinked}. The main result of this section is
Theorem~\ref{curvature} saying that for many globally hyperbolic
space-times the converse is also true, i.e. if $\alk(\mathfrak S_x,
\mathfrak S_y)=0,$ then the events $x, y$ are causally unrelated.

\begin{lem}\label{causalitythroughnullcone}
Let $(\ss^{m+1}, g), m>1$ be a globally hyperbolic space-time
and let $x, y\in \ss$ be events such that $y\in J^+(x).$ Let
$\gamma:(-\infty,\infty)\to \ss $ be an inextendible past directed
nonspacelike curve with $\gamma(0)=y.$ Then there exists a future
directed null geodesic $\nu:[0,\alpha)\to \ss$ with $\nu (0)=x$ and
$\tau_1\in [0, \alpha), \tau_2\in[0,\infty)$ such that $\nu(\tau
_1)=\gamma( \tau _2).$
\end{lem}

\pp Let $J^+(x)$ and $I^+(x)$ be the causal and the chronological
future of $x$. The set $J^+(x)$ is closed since $(\ss, g)$ is
globally hyperbolic, see~\cite[Proposition 3.16]{BeemEhrlichEasley}.
The set $I^+(x)$ is open, see~\cite[Lemma 3.5]{BeemEhrlichEasley}.

Let $M\subset \ss$ be a Cauchy surface containing $x.$ Then $\gamma(t_0)\in M$ 
for some $t_0$. We claim that $t_0\ge 0$. Otherwise the curve 
$\gamma|_{(-\infty,0)}$
followed by the past directed nonspacelike curve joining $y$ to $x$
is a past directed nonspacelike curve that intersects $M$ twice.
This is in contradiction with $M$ being a Cauchy surface.

Clearly $\gamma(t)\not\in J^+(x)$ for $t>t_0\geq 0.$ Since $J^+(x)$
is closed, $I^+(x)\subset J^+(x)$ is open, $\gamma(0)=y\in J^+(x),$
and $\gamma$ is continuous, we conclude that there exists $\tau
_1>0$ with $\gamma(\tau_1)\in J^+(x)\setminus I^+(x).$ Since $(\ss,
g)$ is globally hyperbolic, there exists a future directed null
geodesic $\nu$ from $x$ to $\gamma(\tau _1),$ see~\cite[Corollary
4.14]{BeemEhrlichEasley}. We reparameterize $\nu$ so that it is
future directed and has $\nu(0)=x.$ Now put $\tau_2$ to be such that
$\nu(\tau_2)=\gamma(\tau_1)$ and obtain the statement of the Lemma.
\qed

\m Take $x\in \ss$ and let $C=C^+(x)\subset T_x\ss$ be the hemicone
of all the future pointing null vectors. We have an obvious
$\R^+$-action on $C$. Clearly, $C/\R^+=S^{m-1}$ and in fact we have
a diffeomorphism
$$
C\cong S^{m-1}\times \R = S^{m-1}\times (0, \infty).
$$
Similarly to Riemannian manifolds, one can use geodesics to define
the exponential map $\exp=\exp_x:T_x \ss\to \ss$, cf.
\defref{difgeomLorentz}. Here the domain of $\exp$ is not the whole
$T_x\ss$ but rather a star-shaped with respect to $0\in T_x\ss$
subset $V$ of it. We put $U=V\cap C$.

\begin{lem}\label{omega}
Given a Cauchy surface $M \subset \ss$ with $x\in M$ and $U\subset
T_xX$ as above, take an $M$-proper isometry $h: M\times \R \to \ss$.
Consider the map $F: S^{m-1}\times (0,\infty)\to M \times \R,
F(s,t)=(W^t_{x,M}(s), t)$. Then there exists a diffeomorphism
$\omega: U \to S^{m-1}\times (0,\infty)$ such that the diagram
\begin{equation}\label{diagramomega}
\CD
U @>\exp_x >> \ss\\
@V\omega VV @VVh^{-1} V\\
S^{m-1}\times (0,\infty) @>F>> M\times \R
\endCD
\end{equation}
commutes. Furthermore, $U$ is an open subset of $C$.
\end{lem}

\pp Since $U\subset C\cong S^{m-1}\times (0,\infty)$, every point of
$U$ can be written as $(s,\tau)$ for some $s\in S^{m-1}, \tau \in
(0,\infty)$. Now, given $u=(s,\tau)$, there exists a unique $t=t(u)$
such that $\exp_x u\in M_{t}$. In other words,
\begin{equation}\label{tu}
t(u)=\pi_{\R}(\exp_x(u)).
\end{equation}
We put $\omega(u)=(s,t)$. It is easy to see that the above diagram
commutes and that $\omega$ is a bijection. Furthermore, $\omega$ is
smooth because of \eqref{tu}. Moreover, for each $s$ the velocity
vectors of the curve $\gamma=\gamma_s: \tau\mapsto \exp_x(s,\tau)$
are null, and hence $d\pi_{\R}(\dot\gamma(\tau))\ne 0,$ for all
$\tau$ in the domain of $\gamma$. Thus $\partial t/\partial \tau \ne
0$ everywhere. Now, since $\omega$ preserves the $s$-coordinate, we
conclude that $\omega$ is a diffeomorphism.

Now, the $m$-dimensional manifold $U$ is a subset of the
$m$-dimensional manifold $C$, and so $U$ is open because of the
Invariance of Domain Theorem. \qed

\begin{defin}\label{df:curvature}
The {\em timelike sectional curvatures\/} in a space-time $(\ss, g)$
are the sectional curvatures along the {\it timelike\/} $2$-planes
in $T\ss,$ i.e. $2$-planes $E_s\subset T_s\ss$ such that $g|_{E_s}$
is a nondegenerate form that is not positive definite.
(See~\ref{difgeomLorentz} for a more thorough definition.)
\end{defin}

\begin{prop}\label{prop:immersion}
Let $(\ss^{m+1}, g), m>1$ be a globally hyperbolic space-time where $g$ is 
conformal to $\wh g$ that has all the timelike sectional curvatures
nonnegative. Take a point $x\in \ss$ and Cauchy surface $M \subset
\ss$ with $x\in M$. Choose an $M$-proper isometry $h.$ Define
$$
G=G_g:S^{m-1}\times (-\infty, \infty) \to M, \quad G(s,t)= W^t_{x,M}(s).
$$
Then the restrictions of $G$ onto $S^{m-1}\times (-\infty, 0)$ and
onto $S^{m-1}\times (0,\infty)$ are immersions.
\end{prop}

\pp 
Let $\Omega:\ss\to \R$ be nowhere zero smooth function 
such that $\wh g=\Omega^2 g.$ 
A Cauchy surface $M$ in $(\ss, g)$ is a Cauchy Surface in $(\ss, \wh g).$
Moreover if $h:(M\times \R, \beta dt^2+\ov g)\to (\ss, g)$ is an 
$M$-proper isometry, then $h:\bigl (M\times \R, 
(\Omega^2\circ h)(\beta dt^2+\ov g)\bigr)
\to (\ss, \wh g)$ also is an $M$-proper isometry. 
The null geodesics for $g$ and 
$\wh g $ coincide up to reparameterization, see~\cite[Lemma 
9.17]{BeemEhrlichEasley}, that is not in general an affine reparameterization. 
Hence the maps $G_g, G_{\wh g}:S^{m-1}\times (-\infty, \infty) \to M$ are equal. 
Thus without loss of generality 
we assume that $g$ has all the timelike sectional curvatures nonnegative.

We prove that $G:S^{m-1}\times (0,\infty)\to M$ is an immersion;
the restriction $G:S^{m-1}\times (-\infty, 0)\to M$ can be
considered similarly. For brevity we denote $W_{x,M}^t$ by $W^t$ and
we denote $\exp_x$ by $\exp.$

First, we prove that the map $F:S^{m-1}\times (0,\infty)\to M\times
\R, F(s,t)=(G(s,t),t)$ is an immersion. Let $V\subset T_x\ss$ be the
maximal subset where $\exp$ is defined and let $U\subset V$ be as in
Lemma~\ref{omega}. Let $\rho:[0, b]\to \ss$ be a geodesic starting
at $x.$ The point $\rho(b)$ is conjugate to $x=\rho(0)$ along $\rho$
if and only if the exponential map $\exp:T_x\ss \to \ss$ is
singular at $b\dot \rho(0)\in T_x\ss,$ i.e. if and only if the
differential $(d\exp)_{b\dot\rho(0)}:T_{b\dot\rho(0)}(T_x\ss)\to
T_{\rho(b)}\ss$ is not of full rank,~see~\cite[Proposition 10,
Section 10]{ONeill}. All the non-spacelike geodesics in $(\ss, g)$
do not have any conjugate points, since all the timelike sectional
curvatures in $(\ss^{m+1}, g)$ are nonnegative,
see~\cite[Proposition 11.13]{BeemEhrlichEasley}.

Hence $(d \exp)_u:T_uU\to \ss$ is of full rank, for every $u\in U.$
Therefore $\exp|_U:U \to \ss$ is an immersion. Since the diagram
\eqref{diagramomega} is commutative, we get that $F$ is an
immersion.

For all $(s,t) \in S^{m-1}\times (0,\infty)$, let $V(s,t)$ be the
image of the linear map $ (d W^t)(s): T_sS^{m-1}\to T_a M, \quad
a=W^t(s). $ Since $F$ is the immersion, we conclude that
$W^t:S^{m-1}\to M$ is an immersion for all $t>0.$ So, in order to
show that $G:S^{m-1}\times (0, \infty)\to M$ is an immersion, is
suffices to show that, for all $(s,t)\in S^{m-1}\times (0,
\infty)$,
\begin{equation}\label{notin}
dG(s,t)(\partial/\partial t)\not \in V(s,t).
\end{equation}
So let us take a point $(s_0, t_0) \in S^{m-1}\times (0,\infty)$ and
prove that \eqref{notin} holds for $(s, t)=(s_0,t_0)$.

Let $z=W_{x,M_{t_0}}(s_0)\subset M_{t_0}\subset X$ and
$y=g_{t_0}(z)\in M .$ Put $L=\span\bigl(dh(y,
t_0)(\frac{\partial}{\partial t})\bigr)\subset T_z\ss.$

By definition of an $M$-proper isometry $h$ we have a direct sum
decomposition
\begin{equation}\label{firstdecomp}
T_z\ss=T_zM_{t_0}\oplus L
\end{equation}
and $L$ is the $g$-orthogonal compliment of $T_zM_{t_0}$ in
$T_z\ss$.

Take a null curve $\gamma(t)$ defined by $\gamma(t)=W_{x,
M_t}(s_0)\in M_t\in \ss.$ Clearly up to reparameterization $\gamma$
is a null geodesic through $x.$ Put $\xi=\dot\gamma(t_0)\in T_z\ss$
and use~\eqref{firstdecomp} to decompose $\xi$ as
\begin{equation*}
\xi=\xi_1+\xi_2,\quad \xi_1\in T_zM_{t_0},\, \xi_2\in L.
\end{equation*}
Since $L$ is $g$-orthogonal to $M_{t_0},$ the direction of $\xi_1$
is the direction that defines the lifted wave front $\wt W_{x,
M_{t_0}}$ at $s_0\in S^{m-1}$, cf. \eqref{decomposition}. Since $\wt
W_{x, M_{t_0}}$ is Legendrian, $\xi_1$ is a nonzero vector that is
$g|_{M_{t_0}}$-orthogonal to $\Im (d W_{x, M_{t_0}})(s_0).$ In
particular, $\xi_1\not\in \Im (d W_{x, M_{t_0}})(s_0).$

Since $W^t=g_t(W_{x, M_t})=\pi_{M}(W_{x, M_t})$ for all $t$ and
$L=\ker d\pi_{M}(z),$ we have
$dg_{t_0}(\xi_1)=d\pi_{M}(\xi)=dG(s_0, t_0)(\frac{\partial
}{\partial t})$ and $dg_{t_0}\bigl (\Im (d W_{x,
M_{t_0}})(s_0)\bigr)=\Im dW^{t_0}(s_0)=V(s_0, t_0).$ Since
$g_{t_0}:M_{t_0}\to M$ is a diffeomorphism, we conclude that
$dG(s_0, t_0)(\frac{\partial }{\partial t})\not \in V(s_0,
t_0).$\qed

\m

Recall that by Theorem~\ref{Main} $\A(M)=\Z$ for all smooth
connected oriented $M^m, m>1,$ unless $M$ is a closed manifold that
is an odd dimensional rational homology sphere with finite
$\pi_1(M).$

\begin{thm}\label{curvature}
Let $(\ss^{m+1}, g), m>1$ be a globally hyperbolic space-time where $g$ is 
conformal to $\wh g$ that has all the timelike sectional curvatures
nonnegative. Furthermore, assume that $\A(M)=\Z$ for a Cauchy
surface $M^m$ of $\ss.$ Let $x, y\in \ss$ be two events that do not
lie on a common null geodesic. Then the following statements $(1),$
$(2),$ and $(3)$ are equivalent$:$
\begin{enumerate}
\item $x$ and $y$ are causally related$;$
\item $\alk(\mathfrak S_x, \mathfrak S_y)\neq 0 \in \A(M)=\Z;$
\item the skies $\mathfrak S_x, \mathfrak S_y$ are nontrivially linked in
$\mathcal N$.
\end{enumerate}
\end{thm}

Many space-times satisfying all the conditions of Theorem~\ref{curvature}
are constructed in Example~\ref{exampletimelike}.

\pp By Theorem~\ref{Main} $\alk$ is a link homotopy invariant that
is normalized to be zero when the lifted wave fronts are unlinked.
Thus (2) $\Longrightarrow$ (3). Furthermore, (3) $\Longrightarrow$
(1) by \theoref{thmfrontsunlinked}.

Now we prove that (1) $\Longrightarrow$ (2).
The null geodesics for $g$ and 
$\wh g $ coincide up to reparameterization, see~\cite[Lemma 
9.17]{BeemEhrlichEasley}, that is not in general an affine reparameterization.
Clearly $M$ is a Cauchy surface with respect to both 
$g$ and $\wh g.$ Thus the spaces $\mathcal N=STM$ for $(\ss, g)$ and 
for $(\ss, \wh g)$ are naturally diffeomorphic and the links 
$(\mathfrak S_x, \mathfrak S_y)\subset STM=\mathcal N$ computed for the 
two metrics coincide. Moreover $x,y$ are causally related in $(\ss, g)$ 
if and only if they are causally related in $(\ss, \wh g).$ Thus without 
loss of generality we assume that $g$ has all the timelike sectional curvatures
nonnegative.

Remark~\ref{alksymmetry} says that $\alk(\mathfrak S_x, \mathfrak
S_y)=(-1)^{m}\alk(\mathfrak S_y, \mathfrak S_x).$ Thus it suffices
to prove that $\alk(\mathfrak S_x, \mathfrak S_y)\ne 0$ whenever
$y\in J^+(x)$. So, we assume that $y\in J^+(x)$.

Choose a Cauchy surface $M\ni x$ and an $M$-proper isometry $h: M
\times \R \to \ss$. For brevity we denote $W^t_{x,M}$ by $W^t$. Let
$\tau\in \R$ be the unique value such that $y\in M_{\tau}$. Clearly,
$\tau>0$.

Without loss of generality we can and shall assume that
$\pi_M(x)\neq \pi_M(y).$ Indeed, if $\pi_M(x)=\pi_M(y),$ then we can
construct an auxiliary event $z\in M_{\tau}$ such that $\pi_M(z)\neq
\pi_M(y),$ $z\in J^+(x),$ events $x,z$ do not lie on a common null
geodesic, and $\alk(\mathfrak S_x, \mathfrak S_y)=\alk(\mathfrak
S_x, \mathfrak S_z)$. We construct the event $z$ is follows.

Put $v=\pi_M(y)=\pi_M(x)$ so that $y=h(v, \tau), x=h( v, 0).$ Since
$y\in J^+(x)$ and $x,y$ do not lie on a common null
geodesic,~\cite[Corollary 4.14]{BeemEhrlichEasley} says that $y\in
I^+(x).$ By~\cite[Lemma 3.5]{BeemEhrlichEasley} $I^+(x)$ is open and
hence there exists an open neighborhood $\wt U\subset M$ containing
$v$ such that $h(\wt U, \tau)\subset I^+(x).$ Since $x,y$ do not lie
on a common null geodesic, $y\not \in \Im W_{x, M_{\tau}}.$ Since
$\Im W_{x, M_{\tau}}$ is compact and $y=h(v, \tau)\not \in \Im W_{x,
M_{\tau}},$ there exists an open connected $U\subset \wt U$
containing $v$ such that $h(U, \tau)\cap \Im W_{x,
M_{\tau}}=\emptyset.$ Choose $u\neq v\in U$ and put $z=h(u, \tau).$
Clearly $z\in J^+(x), \pi_{M}(z)=u\neq v=\pi_{M}(x),$ and since
$z\not \in \Im W_{x, M_{\tau}},$ the events $x,z$ do not lie on a
common null geodesic.

Let us prove that $\alk(\mathfrak S_x, \mathfrak S_y)=\alk(\mathfrak
S_x, \mathfrak S_z)$. Take a path $\beta:[0,1]\to U$ with
$\gb(0)=v$ and $\gb(1)=u$. Let $S_1$ and $S_2$ be two copies of
$S^{m-1}$. Define $I:(S_1\sqcup S_2)\times [0,1]\to ST^*M$ by
setting $I(s_1,t)=\wt W^{\tau}(s_1), I(s_2,t)=\eps_{\beta}(s_2,t),$
for $s_1\in S_1, s_2\in S_2, t\in [0,1].$ Since $h(U, \tau)\cap \Im
W_{x, M_{\tau}}=\emptyset$ we get that $U\cap \Im W^{\tau
}=\emptyset.$ Since $\Im \beta\subset U,$ we conclude that $I$ is a
link isotopy and $\alk(\wt W^{\tau}_{x, M}, \eps _v)=\alk(\wt
W^{\tau}_{x, M}, \eps _u).$

Using Theorem~\ref{welldefined} and the above identity we have
\begin{equation}
\begin{aligned}
\alk(\mathfrak S_x, \mathfrak S_y)&=\alk(\wt W_{x, M_{\tau}}, \wt
W_{y, M_{\tau}})= \alk(\wt W^{\tau}_{x, M}, \wt W^{\tau}_{y,
M})\\
&=\alk(\wt W^{\tau}_{x, M}, \eps_v)=\alk(\wt W^{\tau}_{x, M},
\eps_u)= \alk(\wt W^{\tau}_{x, M}, \wt W^{\tau}_{z,
M})\\
&=\alk(\mathfrak S_x, \mathfrak S_z).
\end{aligned}
\end{equation}
Thus, we can and shall assume that $\pi_M(y)\ne \pi_M(x)$.

\m

Let $v=\pi_M(y)$, so that $h(v, \tau)=y$, and let
$\gamma:(-\infty,+\infty)\to \ss$ be an inextendible past directed
timelike curve given by $\gamma(t)=h(v, \tau-t)\in M_{\tau
-t}\subset \ss.$ Lemma~\ref{causalitythroughnullcone} applied to $x,
y$ and $\gamma$ implies that there exists a future directed null
geodesic $\nu:[0, \alpha)\to \ss$ and $\tau_1\in [0, \alpha),
\tau_2\in [0, +\infty)$ such that $\nu(0)=x$ and
$\nu(\tau_1)=\gamma(\tau_2).$ Reparameterize $\nu$ as $\wt \nu
(t)=W_{x, M_{t}}(s), t\geq 0,$ for some $s\in S^{m-1}.$ Then there
exists $\ov \tau_1\in [0, +\infty)$ such that $M_{\ov \tau_1}\ni \wt
\nu(\ov \tau_1)=\gamma(\tau_2)\in M_{\tau -\tau _2}.$ Hence $\ov
\tau_1=\tau-\tau_2$ and since $\ov \tau_1, \tau_2\geq 0,$ we have
$\ov \tau_1\in [0, \tau]$ and $v\in \Im W^{\ov \tau_1}.$

\m
Define $G: S^{m-1}\times [0,\tau] \to M$ by setting $G(s,t)=W^t(s).$
Since $v=\pi_M(y)\neq \pi_M(x)$, we conclude that $v\not \in \Im
G|_{S^{m-1}\times 0}=\Im W^0=x.$ Since $x,y$ do not lie on a common
null geodesic, $y\not\in \Im W_{x, M_{\tau}}$ and therefore $v\not
\in \Im W^{\tau}=\Im G|_{S^{m-1}\times \tau}.$ So $G^{-1}(v)\subset
S^{m-1}\times(0, \tau)$. By \propref{prop:immersion}
$G|_{S^{m-1}\times (0, \tau]}$ is an immersion.
\propref{prop:constant} and the fact $v\in \Im W^{\ov \tau_1}=\Im
G|_{S^{m-1}\times \ov \tau_1}$ imply that all the points in
$G^{-1}(v)\neq \emptyset$ have the same sign. Thus one of $n_+(v,
G)$ and $n_-(v, G)$ is zero and the other is nonzero. By the
assumption of the Theorem $\A(M)=\Z$ and hence $q:\bor_0(\B)=\Z\to
\Z$ has zero kernel. Now \lemref{helpful} for $F(s, t)=\wt W^t_{x,
M}(s), a=0, b=\tau$ and Theorem~\ref{welldefined} imply that
$\alk(\mathfrak S_x, \mathfrak S_y)=\alk(\wt W^{\tau}_{x, M}, \wt
W^{\tau}_{y, M})= \alk(\wt W^{\tau}_{x, M}, \eps_v)=\alk(\wt
W^{0}_{x, M}, \eps_v)+q(n_-(v, G)-n_+(v, G))=\alk(\eps_x,
\eps_v)+q(n_-(v, G)-n_+(v, G))=q(n_-(v, G)-n_+(v, G))\neq 0\in
\A(M).$ \qed

\section{$\alk$-invariant and intersection numbers.}

Recall the following definition.

\begin{defin}[intersection number $f_1\bullet f_2$]
\label{intersectionumber} Let $f_1:N_1\to L^l$ and $f_2:N_2\to L^l$
be transverse mappings of oriented manifolds of complimentary
dimensions into an oriented manifold $L^l.$ Assume that the
preimages of $\Im f_1\cap \Im f_2$ under $f_1$ and $f_2$ are finite
sets. Take $(n_1, n_2)\in N_1\times N_2$ such that
$f_1(n_1)=f_2(n_2)$ and take positive orientation frames $\mathfrak
r_1\subset T_{n_1}N_1$ and $\mathfrak r_2\subset T_{n_2}N_2.$ Put
$\sigma (n_1, n_2)=+1$ if $\{df_1(\mathfrak r_1), df_2(\mathfrak
r_2)\}\subset T_{f_1(n_1)}L=T_{f_2(n_2)}L$ is a positive orientation
frame of $L,$ and put $\sigma( n_1, n_2)=-1$ otherwise. Since $f_1$
and $f_2$ are transverse and $\dim N_1+\dim N_2=\dim L,$ $f_i$ is an
immersion in a neighborhood of $n_i, i=1, 2.$ Hence $\sigma(n_1,
n_2)=\pm 1$ is well defined and does not depend on the choices of
$\mathfrak r_1, \mathfrak r_2.$

The {\em intersection number\/} $f_1\bullet f_2\in \Z$ is defined as
\begin{equation}
f_1\bullet f_2=\sum_{\{(n_1, n_2)\in N_1\times N_2|
f_1(n_1)=f_2(n_2)\}}\sigma(n_1, n_2).
\end{equation}
\end{defin}

Let $(\ss^{m+1}, g), m>1$ be a globally hyperbolic space-time,
and let $x,y$ be two events in $\ss$ that do not lie on a common
null geodesic and such that $y\in J^+(x).$ Let $\exp=\exp_x$ and $U$
be as in \lemref{omega}.

Let $\gamma$ be a future directed past inextendible timelike curve
that ends at $y$ and does not pass through $x$. We say that $\gamma$
is {\em generic} (with respect to $\exp_x|_U$) if it is transverse
to $\exp_{x}|_U$ and it does not pass through the self-intersection
points of $\exp_x|_U.$ It is possible to show that every $\gamma$ as
above can be made generic by a $C^{\infty}$-small deformation. Note
that if $\gamma$ is generic and $\exp(u)\in \Im \gamma,$ for some
$u\in U,$ then $\exp|_U$ is an immersion in a neighborhood of $u,$
since otherwise $\gamma$ and $\exp_x$ are not transverse for
dimension reasons.

\forget

\begin{prop}\label{timeliketransversality}
Let $(\ss^{m+1}, g), m>1$ be a globally hyperbolic space-time,
and let $x_1, x_2\in (\ss, g)$ be such that $x_2\in J^+(x_1)$ and
$x_1, x_2$ do not belong to a common null geodesic. Let $U\subset
C^+(x)\subset T_{x_1}\ss$ be as in Lemma~\ref{omega}.

Choose $\ov x\not \in x_1\sqcup \Im (\exp_{x_1}|_U)$ and let
$\gamma$ be a future directed smooth timelike curve from $\ov x$ to
$x_2.$ Then we can $C^{\infty}$-small perturb $\gamma$ so that
$\gamma$ is a future directed smooth timelike curve from $\ov x$ to
$x_2$ that does not pass through $x_1,$ and moreover $\gamma$ is
transversal to $\exp_{x_1}|_U,$ and it does not pass through the
selfintersection points of $\exp_{x_1}|_U.$
\end{prop}

\pp Choose a Cauchy surface $M\subset X$ containing $x_1$ and an
$M$-proper isometry $h:M\times \R\to \ss.$ For brevity we write
$W^t$ for $W^{t}_{x, M}$ and we write $\exp$ for $\exp_{x_1}|_U.$
Define $\ov F:S^{m-1}\times [0, +\infty)\to M\times \R$ via $\ov
F(s, t)=(W^t(s), t),$ and put $F=\ov F|_{S^{m-1}\times (0,
+\infty)}.$ Put $\ov U=U\sqcup {\bf 0}\in T_{x_1}\ss.$ By
Lemma~\ref{omega} $h\circ F\circ \omega=\exp: U\to \ss,$ where
$\omega:U\to S^{m-1}\times (0, +\infty)$ is a diffeomorphism.

Since $x_1, x_2$ do not belong to a common null geodesic, both $\ov
x, x_2\not \in x_1\sqcup \Im \exp=\Im \exp_{x_1}|_{\ov U}.$ Thus
$h^{-1}(\ov x), h^{-1}(x_2)\not \in \Im \ov F=\Im \bigl (h^{-1}\circ
\exp_{x_1}|_{\ov U}\bigr ).$ Since $S^{m-1}$ is compact and $\ov
F(s, t)=(W^t(s), t)$ is continuous, there are open neighborhoods
of $h^{-1}(\ov x), h^{-1}(x_2)$ that do not intersect $\Im \ov F.$
Hence there are open neighborhoods of $\ov x, x_2$ that do not
intersect $\Im \exp_{x_1}|_{\ov U}.$ For this reason everywhere
below we choose a $C^{\infty}$-small perturbation of $\gamma$ that
is trivial close to $\ov x$ and to $x_2.$ Since being a future
directed timelike curve is an open condition in the space of all
curves in $(\ss, g),$ if the perturbation of $\gamma$ is small
enough, then $\gamma$ will automatically be a future directed
timelike curve from $\ov x$ to $x_2.$

The set of curves that are transversal to $\exp$ is residual in the
space of all smooth curves in $\ss$, see~\cite[Theorem
2.4]{BiasiSaeki}. Thus by a $C^{\infty}$-small perturbation of
$\gamma$ we can assume that $\gamma$ is transversal to $\exp$ and
does not pass through $x_1.$ In particular for dimension reasons,
$\exp$ is an immersion in a neighborhood of every $u\in
\exp^{-1}(\Im \gamma\cap \Im \exp).$ Also the set $\Im \gamma\cap
\Im \exp$ is finite.

Let us show that $\gamma$ can be pushed off the selfintersection
points of $\exp$ by a small perturbation in the neighborhoods of
these points. Assume that $z\in \Im \gamma\cap \Im \exp$ and there
are several points $u_i\in U, i\in I,$ such that $\exp(u_i)=z.$ Put
$(n, t)=h^{-1}(z).$ Then $\omega(u_i)=(n_i, t)\in S^{m-1}\times (0,
+\infty),$ for some $n_i\in S^{m-1}.$ Since $\exp$ is an immersion
in a neighborhood of each $u_i,$ $W^t$ is an immersion in a
neighborhood of each $n_i\in S^{m-1}.$ Since $S^{m-1}$ is compact,
$I$ is a finite set.

Since $S^{m-1}$ is compact, there is a small neighborhood of
$h^{-1}(z)\in M_t$ that is not crossed by any branches of $W^t$
except of the immersed branches passing through $z.$ Thus there is a
small neighborhood of $z\in \ss=M\times \R$ that does not intersect
with any other branches of $\exp$ except of the immersed branches
passing through $z.$

Two immersed branches of $\exp$ can not be tangent at $z\in \ss,$
since otherwise $\wt W^t$ has a double point and is not an
embedding. Thus we showed that for a timelike curve $\gamma$
transversal to $\exp$ and for every $z\in \Im \gamma\cap \Im \exp$,
the number of points $u_i$ such that $\exp(u_i)=z$ is finite and
moreover all the branches of $\exp$ passing through $z$ are immersed
and not tangent to each other at $z.$ Moreover there exists a small
neighborhood of $z$ that is not crossed by any other branches of
$\exp$ except of the immersed branches passing through $z.$ By a
$C^{\infty}$-small perturbation of $\gamma$ in this neighborhood we
can push $\gamma$ off $z$ (that is the set of selfintersection
points of $\exp$ intersected with this neighborhood). Repeat this
procedure for all the points $z\in \Im \gamma\cap \Im \exp.$\qed

\forgotten

\begin{thm}\label{alkasintersection}
Let $(X^{m+1}, g), m>1$ be a globally hyperbolic space-time. Let $x, y, U,$ and 
$\gamma$ be as above. Then $\alk(\mathfrak
S_x, \mathfrak S_y)=q(\exp |_U\bullet \gamma)\in \A(M),$ where
$q:\Z=\bor_0(\B)\to \A(M)$ is the homomorphism from Theorem
$\ref{Main}$ and $M$ is any Cauchy surface.
\end{thm}

\pp Take a Cauchy surface $M$ and an $M$-proper isometry $h: (M
\times \R, \mathfrak g=-\beta dt^2+\ov g) \to (\ss, g).$ Without
loss of generality we assume that $x\in M=M_0$. Since $y\in J^+(x),$
we conclude that $\pi_{\R}(y)>0$. So without loss of generality we
assume that $y\in M_1.$ Since the $\R$-coordinate in $M\times \R$ is
strictly increasing along all future directed non-spacelike curves
we assume (reparameterizing $\gamma$ if necessary) that
$\pi_{\R}(\gamma(t))=t$. For brevity we denote $\exp_x |_U$ by
$\exp,$ we denote $W^t_{a, M}$ by $W^t_a,$ and we denote $\wt
W^t_{a, M}$ by $\wt W^t_a,$ for $a\in M.$

Define
$$
F:S^{m-1}\times [0, +\infty)\to M\times\R,\quad F(s,t)=(W^t_x(s),
t).
$$
By Lemma~\ref{omega} there exists an orientation preserving
diffeomorphism $\omega: U\to S^{m-1}\times (0, +\infty)$ such that
$h\circ F\circ \omega=\exp.$

Since $\gamma$ does not pass through $x=\Im h F|_{S^{m-1}\times 0}$
and $\omega, h$ are orientation preserving, we get that $\exp\bullet
\gamma=F\bullet (h^{-1}\gamma).$ Since $\pi_{\R}(\gamma(t))=t$ we
get that $h^{-1}\gamma(t)\not \in \Im F,$ for $t< 0.$ So in the
computation of $F\bullet (h^{-1}\gamma)=\exp\bullet \gamma$ we can
substitute $\gamma$ by its restriction to $[0,1].$ For brevity we
denote $\gamma |_{[0,1]}$ by $\gamma.$

We define $\ga:[0,1]\to M, \ga(t)=\pi_M(\gamma(t))$ and put
$u=\ga(0), v=\ga(1).$ Define $\wh \gamma:[0,1]\to M\times \R,$
$\wh\gamma(t)=h^{-1}\gamma (t).$ Clearly $\wh \gamma(t)=(\ga(t), t),
t\in [0,1].$ Since $\exp \bullet \gamma=F\bullet
(h^{-1}\gamma)=F\bullet \wh \gamma$ and $\alk(\mathfrak S_x,
\mathfrak S_y)=\alk(\wt W^{1}_{x}, \wt W^{1}_{y})$ by
Lemma~\ref{welldefined}, it suffices to show that $\alk(\wt
W^{1}_{x}, \wt W^{1}_{y})=q(F\bullet \wh \gamma)\in \A(M).$

\m
Let $S^{m-1}_i, i=1,2,$ be a copy of $S^{m-1}$. Let $H:
(S_1^{m-1}\sqcup S_2^{m-1})\times [0, 1]\to ST^*M$ be a link
homotopy given by $H|_{S_1^{m-1}\times t}=\wt W^t_{x}$ and
$H|_{S^{m-1}_2\times t}=\eps_{\ga(t)}, t\in [0,1].$ This link
homotopy deforms the trivial link $(\wt W_{x}^0, \eps_u)=(\eps_{x},
\eps_{u})$ to $(\wt W^{1}_{x}, \wt W^{1}_{y})=(\wt W^{1}_{x},
\eps_{v}).$ Thus to prove the Theorem it suffices to show that the
link homotopy $H$ is generic and the sum of the signs of the
crossings of $\Sigma_0$ during $H$ equals to $F \bullet \wh
\gamma\in \Z.$ {\it We do this below by showing that all the
crossings of $\Sigma$ under $H$ are in the bijective correspondence
with the points of $\Im F\cap \Im \wh \gamma,$ they happen in
$\Sigma_0,$ and are transverse. After this we show that the sign of
a crossing of $\Sigma_0$ under $H$ coincides with the sign of the
corresponding intersection point of $\Im F\cap \Im \wh \gamma.$\/}

Clearly $z\in\Im F\cap \Im \wh \gamma$ exactly when $F(s_1,
\tau)=z=\wh \gamma(\tau),$ for some $s_1\in S^{m-1},\tau\in [0,1].$
Since $\gamma$ is generic, such $s_1, \tau$ are uniquely determined
by $z.$ Hence the points $z\in \Im F\cap \Im \wh \gamma$ are in the
bijective correspondence with the pairs $(s_1,\tau), s_1\in
S^{m-1},\tau\in[0,1]$ such that $W^{\tau}_{x}(s_1)=\ga(\tau)\in M.$

The crossings of $\Sigma$ under the link homotopy $H_t=(\wt W^t_{x},
\eps_{\ga(t)})$ also happen exactly at time moments $\tau$ when
$W^\tau_{x}(s_1)=\ga(\tau),$ for some $s_1\in S^{m-1}_1.$ Since
$\gamma$ is generic, this $s_1$ is uniquely determined by $\tau.$
The preimages of the double point of the singular link are $s_1\in
S^{m-1}_1$ and the unique $s_2\in S^{m-1}_2$ such that
$\eps_{\ga(\tau)}(s_2)=\wt W^{\tau}_x(s_1).$ So we established a
bijective correspondence between the points in $\Im F\cap \Im \wh
\gamma$ and the crossings of $\Sigma$ under $H.$

Let us prove that all these crossings happen in $\Sigma_0.$ Since
$\wt W^{\tau}_{x}, \eps_{\alpha(\tau)}$ are embeddings, we conclude
that condition {\bf a} from \defref{sigma} holds for $f_1=\wt
W_x^{\tau}$ and $f_2=\eps_{\ga(\tau)}$. Note that if $\exp(a)\in
\Im \gamma$ then $\exp$ is an immersion in a neighborhood of $a,$
since $\gamma$ is transverse to $\exp$. Hence \lemref{omega} implies
that $\exp|_{\omega^{-1}(S^{m-1}\times \tau)}$ has image in $M_{\tau}$  and it is an immersion in a neighborhood of $\omega^{-1}(s_1, \tau).$ Thus $W^{\tau}_{x}=\pi_M \exp|_{\omega^{-1}
(S^{m-1}\times \tau)}=\pr \circ \wt W^{\tau}_x$ is an immersion in a
neighborhood of $s_1.$ Since $\Im
\pr\circ\eps_{\alpha(\tau)}=\alpha(\tau),$ we get that $\Im (d\wt
W^{\tau}_x) (T_{s_1}S^{m-1}_1)\cap \Im (d
\eps_{\alpha(\tau)})(T_{s_2}S^{m-1}_2)={\bf 0}$ and condition {\bf
b} of~\defref{sigma} holds. Thus all the crossings of $\Sigma$ under
link homotopy $H$ happen in $\Sigma_0.$

If we show that a tangent frame to $\wt W^\tau_{x}$ at $\wt
W^\tau_{x}(s_1)$, a tangent frame to $\eps_{\ga(\tau)}$ at
$\eps_{\ga(\tau)}(s_2)=\wt W^\tau_{x}(s_1),$ and the vector $\bf w$
from the definition of $\sigma(H, \tau)$ form a linearly independent
family, then the crossing of $\Sigma_0$ is transverse. Consider the
differential $d\pr: TST^*M \to TM$ of $\pr: ST^*M \to M$. Since
$\eps_{\ga(\tau)}$ is the inclusion of an $S^{m-1}$-fiber of
$\pr:ST^*M\to M,$ it suffices to show that the images under $d\pr$
of ${\bf w}$ and of a tangent frame to $\wt W^\tau_{x}$ at $\wt
W^\tau_{x}(s_1)$ are linearly independent in
$T_{W_{x}^{\tau}(s_1)}M$. As we remarked, $W^{\tau}_{x}$ is an
immersion in a neighborhood of $s_1.$ So to prove that the crossing
of $\Sigma_0$ is transverse is suffices to show that $d\pr({\bf
w})\notin \Im (dW^\tau_{x})(T_{s_1}S^{m-1}_1).$

Given $s\in S^{m-1}$, we define $\gb_s:[0,1]\to
M,\,\gb_s(t)=W^t_x(s).$ Clearly $h(\gb_{s_1}(t), t), t\in [0,1],$ is
(up to a reparameterization) an arc of the null geodesic whose
velocity vectors define the points $\wt W_{x, M_{t}}(s_1)\in
ST^*M_t, t\in [0,1].$ For brevity we denote the vector field
$\partial/\partial t$ on $\R$ by $\partial_t$. Put
$$\xi=\dot{\beta}_{s_1} (\tau)\in
T_{\ga(\tau)}M=T_{W^{\tau}_{x}(s_1)}M
$$
and note that $\wh \xi:=\xi +\partial_t\in T_z(M\times \R)$ is the
velocity vector of the curve $(\gb_{s_1}(t), t), t\in [0,1].$ Since
$h$ is an isometry, $\wh \xi$ is a future pointing null vector with
respect to the Lorentz metric $\mathfrak g.$

Put
$$
\eta=\dot\ga (\tau)\in T_{\ga(\tau)}(M\times \R)
$$
Note that $\wh \eta:=\eta+\partial_t$ is the velocity vector of $\wh
\gamma$ and hence it is a future pointing timelike vector with
respect to the Lorentz metric $\mathfrak g.$

It is easy to see that $d \pr({\bf w})=\xi-\eta.$

Let $\ov g=h_{\tau}^*(g)$ be the Riemannian metric on $M.$ The
direction of the vector $\xi$ is $\wt W^{\tau}_{x}(s_1),$ where we
identify $ST^*M$ and $STM$ via the metric $\ov g.$ Since $\wt
W^{\tau}_{x}$ is Legendrian, $\xi$ is $\ov g$-orthogonal to $\Im
dW^{\tau}_{x}(T_{s_1}S_1^{m-1}).$ To show that $d \pr({\bf
w})=\xi-\eta \not \in \Im dW^{\tau}_{x}(T_{s_1}S_1^{m-1})$ it
suffices to show that $\ov g\bigl( \xi-\eta, \xi)>0.$

The Lorentz product of a future pointing timelike and of a future
pointing null vector is negative. Hence
\begin{equation}\label{eqvec0}
\begin{aligned}
0<0-\mathfrak g(\wh \eta, \wh \xi)&=\mathfrak g(\wh \xi, \wh
\xi)-\mathfrak g(\wh \eta, \wh \xi)= \mathfrak g(\wh \xi-\wh \eta,
\wh \xi) \\
&=\mathfrak g \bigl(\xi-\eta, \xi\bigr)= \ov g(\xi-\eta, \xi).
\end{aligned}
\end{equation}
Thus $\ov g\bigl( \xi-\eta, \xi)>0$ and all the crossings of
$\Sigma_0$ under $H$ are transverse.

\m To finish the proof of the Theorem it suffices to show that the
intersection point $z\in \Im F\cap \Im \wh \gamma$ has the same sign
as the corresponding crossing of $\Sigma_0$ under homotopy $H.$

Let $\mathfrak r$ be a positive orientation frame in
$T_{s_1}S^{m-1}_1$. Then the sign $\sigma\bigl ( (s_1,\tau),
s_2\bigr)$ of the intersection point $z$ is the sign of the
orientation of $M\times \R$ given by the frame $\{dF (\mathfrak r),
\wh \xi, \wh \eta\}$.

The vectors of $d F (\mathfrak r)$ are tangent to $M\times
\tau\subset M\times\R$ and are spacelike with respect to $\mathfrak
g.$ The vector $\wh \xi$ is null. The straight line homotopy
$(1-\gl)\wh\xi +\gl\xi, \gl\in [0,1],$ of $\wh \xi$ to the
spacelike vector $\xi$ is $\mathfrak g$-orthogonal to $dF(\mathfrak
r)$ and induces a homotopy of the frame $\{(dF) (\mathfrak r),
\wh\xi, \wh \eta\}$ to the frame $\{(dF) (\mathfrak r), \xi, \wh
\eta\}.$

\m We claim that the frame stays nondegenerate during the homotopy,
so that the orientations of $M\times \R$ given by the initial frame
and the final frame are equal. If this is false, then since the
vectors in $dF(\mathfrak r)$ are linearly independent and have zero
$\R$-coordinate, we get that there exist a {\em spacelike} vector
$\zeta\in T_{\ga(\tau)}M\subset T_z(M\times \R)$, a value $\gl\in
[0,1]$, and $a,b\in\R$ (with at least one of $a,b$ nonzero) such
that $\zeta\in \span (dF(\mathfrak r))$ and such that
\begin{equation}\label{eqvec1}
a(\eta+\partial_t)+b\Bigl((1- \gl) (\xi+\partial_t)
+\gl\xi\Bigr)+\zeta={\bf 0}.
\end{equation}
Equating the coefficients at $\partial_t$, we see that
$a=-b(1-\gl),$ and since at least one of $a,b$ is non-zero, $b\neq
0.$ Substitute $a=-b(1-\gl)$ into~\eqref{eqvec1} to get
\begin{equation}\label{eqvec2}
-b(1-\gl)\eta+b\xi+\zeta={\bf 0}.
\end{equation}
Since $\xi$ is $\mathfrak g$-orthogonal to $dF(\mathfrak r)$ and
$\zeta\in \span (dF(\mathfrak r))$, we conclude that $\mathfrak
g\bigl(\xi, \zeta\bigr)=0.$ Thus from~\eqref{eqvec2} we have we that
$\gl\neq 1,$ and hence
$\eta=\frac{1}{(1-\gl)}\xi+\frac{1}{b(1-\gl)}\zeta.$ Thus
\begin{equation}\label{eqvec3}
\begin{aligned}
\mathfrak g (\eta,\eta)&= \frac{1}{(1-\gl)^2}\mathfrak g
(\xi,\xi)+\frac{2}{b(1-\gl)^2}\mathfrak g (\xi,\zeta
)+\frac{1}{b^2(1-\gl)^2}\mathfrak
g (\zeta,\zeta)\\
&= \frac{1}{(1-\gl)^2}\mathfrak g (\xi,\xi
)+0+\frac{1}{b^2(1-\gl)^2}\mathfrak g (\zeta, \zeta )>\mathfrak g
(\xi,\xi),
\end{aligned}
\end{equation}
since $\lambda\in [0,1]$ and $\zeta$ is a spacelike vector. Since
$\mathfrak g\bigl(\xi, \partial_t \bigr)=0=\mathfrak
g\bigl(\eta,\partial _t\bigr),$ the vectors $\wh \eta$ and
$\partial_t$ are timelike, the vector $\wh \xi$ is null,
and~\eqref{eqvec3} holds, we conclude that
\begin{equation}\label{eqvec4}
\begin{aligned}
0&>\mathfrak g(\wh \eta, \wh \eta)=\mathfrak g\bigl( \eta +\partial_
t, \eta +\partial_ t\bigr)=\mathfrak g\bigl(\eta,
\eta\bigr)+\mathfrak g(\partial_t, \partial_t)\\
&>\mathfrak g\bigl(\xi,\xi\bigr)+\mathfrak g(\partial_t, \partial
_t)= \mathfrak g(\wh \xi, \wh \xi)=0.
\end{aligned}
\end{equation}
This is a contradiction.

\m Since $M\times \R$ is a product of oriented manifolds and the two
frames above give equal orientations of it, we see that the sign
$\sigma\bigl ((s_1, \tau), \tau)$ of the intersection point $z\in
\Im F\cap \Im \wh \gamma$ is positive exactly when
$\{d(W^{\tau}_{x})(\mathfrak r), \xi\}=\{d(W^{\tau}_{x})(\mathfrak
r), \dot{\beta}_{s_1}({\tau})\}$ is a positive orientation frame of
$M.$

Recall that $ST^*M$ was oriented in such a way that an $m$-frame
projecting to a positive frame on $M^m$ followed by a positive
orientation frame of the $S^{m-1}$-fiber is a positive orientation
frame of $ST^*M.$ Since $\eps_{\ga(\tau)}$ is an inclusion of the
positively oriented fiber, we conclude that $\sigma(H, \tau )=+1$
exactly when $\{d(W^{\tau}_{x})(\mathfrak r), d \pr({\bf
w})\}=\{d(W^{\tau}_{x})(\mathfrak r), \xi-\eta\}$ is a positive
orientation frame of $M.$

Since $\xi$ is $\ov {g}$-orthogonal to the immersed branch of the
front $W^{\tau}_{x}$, and since by~\eqref{eqvec0} $\ov {g}
(\xi-\eta, \xi)>0,$ we conclude that $\xi$ and $\xi-\eta$ point to
the same half-space of $T_{W^{\tau}_{x}(s_1)}M\setminus \bigl(\Im
d(W^{\tau}_{x})(T_{s_1}S^{m-1}_1)\bigr).$ Thus the orientations of
$M$ given by the frames $\{d(W^{\tau}_{x})(\mathfrak r), \xi\}$ and
$\{d(W^{\tau}_{x})(\mathfrak r), d \pr({\bf w})\}$ are equal. Hence
the signs of the intersection points of $F$ with $\wh \gamma$ and of
the corresponding crossings of $\Sigma_0$ under $H$ coincide. \qed

\section{Computing the increment of $\alk$ under the passage through a dangerous
tangency.}\label{mod} 

Let $(\ss^{m+1},g),m>1$ be a globally hyperbolic space-time. 
Definition~\ref{definalk} and Lemma~\ref{Bconnected} imply
that $\alk(\mathfrak S_x, \mathfrak S_y)$ can be computed as
follows. Take a Cauchy surface $M^m \subset \ss$ and $t\in \R$ and
choose a generic path $\ga:[a,b]\to \cs\times \cs$ that deforms a
pair $(\eps_{u}, \eps_{v})$ to $(\wt W_{x,M}^t, \wt
W_{y,M}^t)\subset ST^*M.$ Put $t_i, i\in I,$ to be the time moments
when $\ga$ crosses $\Sigma.$ Since $\ga$ is generic these crossings
happen in $\Sigma_0,$ and we put $\sigma(\ga, t_i)=\pm 1, i\in I,$
to be the signs of these crossings, see~\ref{signsigma0}. By
\theoref{Main}
$$
\alk(\mathfrak S_x, \mathfrak S_y)=q\bigl(\sum _{i\in I} \sigma(\ga,
t_i)\bigr)\in \A(M).
$$
Such a path $\ga:[a,b]\to \cs\times \cs$ can be described as a
family of maps $ \wt \ga^{\tau}: S_1^{m-1}\sqcup S_2^{m-1}\to ST^*M,
\tau\in [a,b]. $ It also can be described as a family of maps $ \un
\ga^{\tau}=\pr\circ \wt\ga^{\tau}: S_1^{m-1}\sqcup S_2^{m-1}\to M,
\tau\in [a,b], $ equipped with a covector field $\theta^{\tau}_s\in
T^*_{\un\ga^{\tau}(s)}M, s\in S_1^{m-1}\sqcup S_2^{m-1},$ that
defines the lift of $\un \ga^{\tau}$ to $\wt \ga^{\tau}.$ In terms
of the last description the crossings of $\Sigma$ by $\ga$
correspond to the triples $(\tau, s_1, s_2)\in [a,b]\times
S^{m-1}_1\times S^{m-1}_2$ such that $\un \ga^{\tau}(s_1)=\un
\ga^{\tau}(s_2)\in M$ and the nonzero covectors
$\theta^{\tau}_{s_1}, \theta^{\tau}_{s_2}$ are positive multiples of
each other in $T_{\un\ga^{\tau}(s_1)}^*M=T_{\un\ga^{\tau}(s_2)}^*M.$
(The triples $(\tau, s_1, s_2)$ at which $\un \ga^{\tau}(s_1)=\un
\ga^{\tau}(s_2)$ and the covectors $\theta^{\tau}_{s_1},
\theta^{\tau}_{s_2}$ are negative multiples of each other do not
correspond to the double points of $\wt \ga^{\tau},$ since then $\wt
\ga^{\tau}(s_1)$ and $\wt \ga^{\tau}(s_2)$ are the opposite points
of the $S^{m-1}$-fiber of $\pr:ST^*M\to M.$)

As we know, the lifted wave fronts $\wt W_{x,M}^{t}, \wt W_{y,M}^t$
are Legendrian embeddings $S^{m-1}\to ST^*M$ that are each
Legendrian isotopic to a Legendrian embedding $\eps_w:S^{m-1}\to
ST^*M, w\in M.$ Put $\mathcal L\subset \cs$ to be the connected
component of the space of Legendrian immersions $S^{m-1}\to ST^*M$
that contains the Legendrian embeddings $\eps_w, w\in M.$ Clearly
the generic path $\ga$ joining $(\eps_{u}, \eps_{v})$ to $(\wt
W_{x,M}^t, \wt W_{y,M}^t)$ can be chosen so that $\Im (\ga)\subset
\mathcal L\times \mathcal L\subset \cs\times \cs.$ Let
$\gl:[a,b]\to\mathcal L\times \mathcal L $ be such a path. (We
changed the notation from $\ga$ to $\gl$ in order to emphasize that
$\gl$ is a path in $\mathcal L$ rather than in the whole $\cs$.)
Let $\wt \gl^{\tau}:S^{m-1}_1\sqcup S^{m-1}_2\to ST^*M, \tau\in
[a,b],$ be the corresponding family of maps and let $\un
\gl^{\tau}=\pr\circ \wt \gl^{\tau}:S^{m-1}_1\sqcup S^{m-1}_2\to M,
\tau\in [a,b],$ be the family of maps equipped with a covector field
$\theta_s^{\tau}$ that defines the lift of $\un \gl^{\tau}$ to $\wt
\gl^{\tau}.$ Since $\wt \gl^{\tau}$ are Legendrian, the covectors
$\theta^{\tau}_{s}\in T_{\gl^{\tau}(s)}^*M$ vanish on $(\un
\gl^{\tau})_*(T_{s}S_i^{m-1})$ for $s\in S_i^{m-1}, i=1,2.$ If
$\gl:[a,b]\to \mathcal L\times \mathcal L$ is generic, then the
crossings of $\Sigma$ by $\gl$ happen in $\Sigma_0$ and correspond
to the triples $(\tau, s_1, s_2)$ as above with an extra condition
that $\un \gl^{\tau}$ restricted to small neighborhoods of $s_1,
s_2$ is an immersion. Since $\theta^{\tau}_{s}$ vanishes on $(\un
\gl^{\tau})_*(T_{s}S_i^{m-1})$ for $s\in S_i^{m-1}, i=1,2,$ we get
that $\gl^{\tau}|_{S_1^{m-1}}$ and $\gl^{\tau}|_{S_2^{m-1}}$ are
tangent at $\gl^{\tau}(s_1)=\gl^{\tau}(s_2).$ Combining all this
together we see that the crossings of $\Sigma_0$ by a generic
$\gl:[a,b]\to \mathcal L\times \mathcal L$ correspond to the so
called {\em Arnold's~\cite{Arnold} dangerous tangencies\/} of $\un
\gl^{\tau}|_{S_1^{m-1}}$ and $\un \gl^{\tau}|_{S_2^{m-1}}.$ These
are the instances when the immersed branches of $\un
\gl^{\tau}|_{S_1^{m-1}}$ and $\un \gl^{\tau}|_{S_2^{m-1}}$ are
tangent at exactly one point, this tangency point has exactly one
preimage on each of $S_1^{m-1}, S_2^{m-1},$ and the covectors
defining the Legendrian lifts of $\un \gl^{\tau}|_{S_1^{m-1}}$ and
$\un \gl^{\tau}|_{S_2^{m-1}}$ at the tangency point are positive
multiples of each other. (We will see that the tangency point is of
order one, since $\lambda$ is a generic path and $\sigma(\lambda,
\tau_0)$ is well-defined.)

\m Below we give a formula for computing the sign $\sigma(\gl,
\tau_0)$ of the crossing of $\Sigma_0$ that corresponds to the
passage through Arnold's dangerous tangency of $\un
\gl^{\tau_0}|_{S_1^{m-1}}$ and $\un \gl^{\tau_0}|_{S_2^{m-1}}.$

Put $\un \gl_i^{\tau}=\un \gl^{\tau} |_{S^{m-1}_i}, i=1,2,$ and
equip them with the restrictions of the covector field
$\theta^{\tau}_s.$ Consider a positively oriented chart $\gf: U\to
\R^m, U\subset M$ with local coordinates $\{x_1, \ldots, x_m\}$
such that:

\begin{itemize}
\item $\un\gl^{\tau_0}_1(s_1)=\gf^{-1}({\bf 0})=\un\lambda^{\tau_0}_2(s_2)\in M$
is the dangerous tangency point;
\item the restriction of $\un\gl^{\tau_0}_i$ to the preimage $V_i$ of $U$ under
$\un\gl^{\tau_0}_i$ is an embedding, $i=1,2;$
\item the common tangent hyperplane to
$\gf \un\gl^{\tau_0}_i|_{V_i}, i=1,2$ at the point
$\gf \un \gl^{\tau_0}_1(s_1)={\bf 0}=\gf \un \gl^{\tau_0}_2(s_2)$
is given by the equation $x_m=0$;
\item $(\gf^{-1})^*(\theta^{\tau_0}_{s_1})$ is a positive multiple of $-dx_m;$
\item $\Im\bigl(\gf \un \gl_i^{\tau_0}|_{V_i}\bigr)$ is given by an equation
$x_m=f_i(x_1, \ldots, x_{m-1}),$ for some smooth function $f_i,
i=1,2.$
\end{itemize}

Since $\un \gl^{\tau_0}_1$ and $\un \gl^{\tau_0}_2$ are dangerously
tangent at $\gf^{-1} ({\bf 0}),$ we conclude that
$(\gf^{-1})^*(\theta^{\tau_0}_{s_2})$ is a positive multiple of
$-dx_m.$ We put $\eps$ to be $+1$ if $\un \gl_1^{\tau_0}$ and $\un
\gl_2^{\tau_0}$ induce the same orientation on the common tangent
$(m-1)$-plane at $\gf^{-1}({\bf 0})$, and we put $\eps=-1$
otherwise. Put $g=f_2-f_1$ and let $\Hess g({\bf 0})$ denote the
Hessian of $g$ at ${\bf 0}\in \R^{m-1}.$ Put $\alpha$ to be the sign
of the $m$-th coordinate of the difference $\gf(\un
\gl_2^{\tau'}(s_2))-\gf(\un \gl_1^{\tau'}(s_1))\in \R^m,$ for
$\tau'$ slightly bigger than $\tau_0.$

\begin{theorem}\label{hesse}
\[
\sigma(\gl, t_0)=(-1)^k\alpha\eps=\ga\eps\sign(\det \Hess g({\bf
0}))
\]
where $k$ is the number of negative eigenvalues of $\Hess g({\bf
0})$.
\end{theorem}

\begin{rem}
Since we consider the passage through a dangerous tangency point
that corresponds to a transverse crossing of $\Sigma_0$, we know
that $\sigma(\gl, \tau_0)$ is defined. In particular, by
Theorem~\ref{hesse} $\Hess g({\bf 0})$ is nondegenerate and hence
the tangency point is of order one. Similarly $\alpha$ is
well-defined, i.e. the difference of the $m$-th coordinates that we
used to define $\alpha$ is nonzero.

A version of this Theorem appeared in our
preprint~\cite{ChernovRudyakFronts}. Also some ingredients of this
formula appeared in the work of T.~Ekholm, J.~Etnyre, and
M.~Sullivan~\cite[Proposition 3.3 and Lemma 3.4]{EES} in a different
situation, where the authors compute the Thurston-Bennequin
invariant of a Legendrian submanifold of $\R^{2n+1}.$
\end{rem}

\pp Since $\gf:U\to \R^m$ is a positively oriented chart, we get
that the sign of the crossing of $\Sigma_0$ under the lifts of $\un
\gl_1^{\tau}$ and of $\un \gl_2^{\tau}$ to $ST^*M$ is equal to the
sign of the crossing of $\Sigma_0$ under the lift to $ST^*\R^m$ of
the branches of $\gf \un \gl_1^{\tau}$ and of $\gf \un \gl_2^{\tau}$
in $\R^m$ that are equipped by the covector field
$(\gf^{-1})^*(\theta_s^{\tau}), s\in S^{m-1}_1\sqcup S^{m-1}_2.$ We
use the flat Riemannian metric on $\R^m$ to identify $ST^*\R^m$ and
$ST\R^m.$ Under this identification the codirection of a covector
$\theta\in T^*\R^m$ corresponds to the direction of the vector
$\theta^+\in T\R^m$ that is orthogonal to $\ker \theta$ and
satisfies $\theta(\theta^+)>0.$ Thus to prove Theorem~\ref{hesse},
it suffices to show that the formula in its formulation indeed gives
the sign of the crossing of $\Sigma_0$ under the lifts of the
branches of $\gf \un \gl_1^{\tau}$ and of $\gf \un \gl_2^{\tau}$ to
$ST\R^m.$

\m If one changes orientation of one of the two $S^{m-1}$-spheres
parameterizing $\un \gl_1^{\tau}$ and $\un \gl_2^{\tau},$ then both
expressions in the statement of Theorem~\ref{hesse} change sign.
Thus without loss of generality we can assume that the
orientations induced by $\gf \un \gl_1^{\tau_0}$ and by $\gf \un
\gl_2^{\tau_0}$ on a common tangent hyperplane $\{ (x_1, \ldots,
x_m)|x_m=0 \}\subset T_{\bf 0}\R^m$ are equal to the standard
orientation of the $\R^{m-1}$-plane, and hence $\eps=1$.

\m Without loss of generality we identify $V_i, i=1,2,$ with
$\R^{m-1}$, and we put $(v_1, \ldots, v_{m-1})$ to be the
coordinates on $\R^{m-1}.$ We parameterize the branches of $\gf
\gl_i^{\tau_0}|_{V_i}, i=1,2,$ by the maps $\R^{m-1}\to \R^m$. After
an orientation preserving reparameterization, the branch $\gf \un
\gl_i^{\tau_0}|_{V_i}, i=1,2$ is given by the parametric equations
\[
x_k=v_k \text{ for } k=1, \ldots, m-1,\quad x_m=f_i(v_1, \ldots,
v_{m-1}).
\]

We consider the unit hemisphere $S^-=\{(x_1, \ldots, x_m)\bigm |
x_1^2+\cdots +x_m^2=1, x_m< 0\}\subset \R^m$, and we equip $S^-$
with local coordinates $\{y_1, \ldots, y_{m-1}\}$ by setting
$y_k(p)=x_k(p)$ for all $p\in S^-$ and $k=1, \ldots, m-1$.

Put $\wt \mu_i^{\tau}, i=1,2$ to be the lift of the branch of $\gf
\un \gl^{\tau}_i$ to $ST\R^m.$ It is obtained by mapping a point
$v\in \R^{m-1}$ to the direction of the unit vector normal to $\gf
\un \gl^{\tau}_i$ at $\gf \un \gl^{\tau}_i(v)$ on which the
corresponding covector $(\gf^{-1})^*(\theta^{\tau}_v)$ is positive.
So at the dangerous tangency point ${\bf 0}\in \R^m$ the two unit
length vector fields defining the lifts $\wt \mu_i^{\tau}, i=1,2,$
are equal to $-{\partial}/{\partial x_m}.$

Let $b$ be the unique point in $\Im \wt \mu_1^{\tau_0}\cap \Im \wt
\mu_2^{\tau_0}.$ Clearly $\pr(b)={\bf 0}\in \R^m$ and the product
$\R^m\times S^-$ can be considered as the codomain of the chart
$\psi=\{x_1, \ldots, x_m, y_1, \ldots, y_{m-1}\}$ at $b\in ST\R^m$.
The parametric equations for the lifts $\wt
\mu_i^{\tau_0}:\R^{m-1}\to ST\R^{m-1}, i=1,2$ are
\begin{equation}\label{fronts}
\begin{split}
x_k=v_k \text{\, for\, } k=1, \ldots, m-1;\quad
x_m=f_i(v_1, \ldots, v_{m-1}); \\
y_k=\frac{1}{r_i}\pa{f_i}{v_k} \text{\, for\, } k=1, \ldots, m-1,
\text{ where }
r_i=\sqrt{1+\sum_{k=1}^{m-1}\left(\pa{f_i}{v_k}\right)^2}.
\end{split}
\end{equation}

This holds since $(y_1, \ldots, y_{m-1},-1/r_i)$ is the unit normal
vector to $\Im\gf \un \gl_i^{\tau_0}$ at $\gf \un
\gl_i^{\tau_0}({\bf v})$ and this normal vector for ${\bf v}={\bf
0}\in \R^{m-1}$ coincides with $-{\partial}/{\partial x_m}.$

Let $\mathbf w$ be the vector from \defref{signsigma0}. Let $\mathbf
w=(\ga_1, \ldots, \ga_{2m-1})$ in the chart $\psi$. Clearly $\ga$
from the statement of Theorem~\ref{hesse} is equal to the sign of
$\ga_m$.

To make the notation simpler for a function $h:\R^{m-1}\to \R$ we
put
$$
\partial_k h=\pa{h}{v_k}\quad \text{ and } \quad 
\partial _{k,l}h=\frac{\partial
^2 h}{\partial v_k\partial v_l}.
$$
For $i=1,2$, the positive tangent frame to $\wt \mu_i^{\tau_0}$ is
given by vectors
$$
\xi_k^{(i)}=(\partial_k x_1, \ldots, \partial_k x_m, \partial_k
y_1,\ldots,
\partial_ky_{m-1}),\ k=1, \ldots, m-1
$$
where $x_k$ and $y_k$ are from \eqref{fronts} with the corresponding
value of $i$. So, according to \defref{signsigma0}, the sign
$\sigma(\gl, \tau_0)$ is equal to the sign of the polyvector
$$
\xi_1^{(1)} \wedge \cdots \wedge \xi_{m-1}^{(1)} \wedge \mathbf w
\wedge \xi_1^{(2)} \wedge \cdots \wedge \xi_{m-1}^{(2)},
$$
i.e. to the sign of the determinant with column vectors
$\xi^{(1)}$'s, $\mathbf w$ and $\xi^{(2)}$'s computed at ${\bf
v}={\bf 0}.$ Clearly $$\partial_l\Bigl(\frac{\partial
_k{f_i}}{r_i}\Bigr)=\frac{r_i\partial_{k,l} f_i-\partial _k
f_i\partial _l r_i}{r_i^2}.$$ Since $\partial_k f_i({\bf 0})=0, k=1,
\cdots, m-1,$ and $r_1({\bf 0})=1=r_2({\bf 0}),$ we get that
$\partial _l(y_k)=\partial_l (\frac{\partial _k{f_i}}{r_i})({\bf
0})=\partial_{k,l}f_i({\bf 0})$ for $y_k$ from~\eqref{fronts} with
the corresponding value of $i.$ Thus $\sigma(\gl, \tau_0)$ equals to
the sign of the determinant
\begin{equation}\label{hugematrix}
\begin{vmatrix}
1 & \cdots & 0 &\ga_1 & 1 &\cdots & 0 \cr
0 & \cdots & 0 &\ga_2 & 0 &\cdots & 0 \cr
\vdots &\vdots & \vdots & \vdots &\vdots & \vdots &\vdots\cr
0 & \cdots & 1 &\ga_{m-1} & 0 &\cdots & 1 \cr
\partial_1 f_1 &\cdots &
\partial_{m-1} f_1&\ga_{m} & \partial_1 f_2 
&\cdots & \partial_{m-1} f_2\cr
\partial_{11} f_1 &\cdots &
\partial_{m-1,1} f_1&\ga_{m+1} & \partial_{11} f_2 
&\cdots & \partial_{m-1,1} f_2\cr
\partial_{12} f_1 & \cdots &
\partial_{m-1,2} f_1&\ga_{m+2} & \partial_{12} f_2 
&\cdots & \partial_{m-1,2} f_2\cr
\vdots & \vdots &\vdots & \vdots &\vdots &\vdots &\vdots \cr
\partial_{1, m-1} f_1 &\cdots & \partial_{m-1,m-1} f_1&\ga_{2m-1} 
& \partial_{1, m-1} f_2 &\cdots & \partial_{m-1,m-1} f_2 \cr
\end{vmatrix}
\end{equation}
evaluated at ${\bf 0}\in \R^{m-1}=\{(v_1, \ldots, v_{m-1})\}$. Here
the up-left and up-right $(m-1)\times (m-1)$ blocks of the matrix
are identity matrices.

Subtract the $k$-th column from the $(m+k)$-th one, $k=1, \ldots , m-1$ 
to get the determinant
\begin{equation}\label{subtraction}
\begin{vmatrix}
1 & \cdots & 0 &\ga_1 & 0 &\cdots & 0 \cr 0 & \cdots & 0
&\ga_2 & 0 &\cdots & 0 \cr \vdots &\vdots & \vdots &
\vdots &\vdots &\vdots &\vdots\cr 0 & \cdots & 1
&\ga_{m-1} & 0 &\cdots & 0 \cr
\partial_1 f_1 &\cdots & \partial_{m-1} f_1&\ga_{m} 
& \partial_1 g &\cdots &
\partial_{m-1} g\cr
\partial_{11} f_1 &\cdots &
\partial_{m-1,1} f_1&\ga_{m+1} & \partial_{11} g 
&\cdots & \partial_{m-1,1} g\cr
\partial_{12} f_1 &\cdots &
\partial_{m-1,2} f_1&\ga_{m+2} & \partial_{12} g 
&\cdots & \partial_{m-1,2} g\cr \vdots & \vdots &\vdots & \vdots
&\vdots &\vdots &\vdots \cr
\partial_{1, m-1} f_1 &\cdots & \partial_{m-1,m-1} f_1&\ga_{2m-1} 
& \partial_{1, m-1} g &\cdots & \partial_{m-1,m-1} g\cr
\end{vmatrix}
\end{equation}
evaluated at ${\bf 0}.$ Since $\partial_kg({\bf 0})=0=\partial
_kf_1({\bf 0}), k=1\cdots, m-1$, this determinant equals to
$\ga_m\det \Hess g({\bf 0})$ and we proved the Theorem. \qed

\forget
\begin{rem}
Let us explain how to compute $\sign \det \Hess f(0)$ in a
coordinate-free way, i.e. without choosing the chart $\gf$. Choose a
Riemannian metric on $M$ and take the unique torsion free connection
$\nabla$ on $M$ compatible with the Riemannian metric. Let $T$ be
the common tangent space to $\un \lambda^{\tau_0}_i, i=1,2$ at $\un
\lambda^{\tau_0}_1(s_1)=\un \lambda^{\tau_0}_1(s_2).$ For $i=1,2$
let ${\mathbf n_i}$ be the normal unit vector field to $\un
\lambda^{\tau_0}_i$ near $\un \lambda^{\tau_0}_1(s_1)=\un
\lambda^{\tau_0}_1(s_2)$ on which the evaluation of the covector
field $\theta^{\tau_0}_{s}, s\in S_1^{m-1}\sqcup S_2^{m-1}$ is
positive.

For $i=1,2$ consider the Weingarten operator
\begin{equation}\label{weingarten}
A_i: T \to T, \quad A(\e)=\nabla_{\e}{\mathbf n_i}.
\end{equation}
We set $B=A_2 - A_1$. It is well known that each $A_i$ is a
self-adjoint operator, \cite [Ch. 7]{Spivak}, and therefore $B$ is.
One can prove that the sign of $\det B$ is equal to the sign of
$\det\Hess f(a)$. This follows, since in the above considered
coordinate-choosing case we have $B=\Hess f(a)$.
\end{rem}
\forgotten

\forget
\begin{rem}
We note interesting similarities between the terms in our formula
for the increment of the $\alk$ invariant and in the formulas for
the Bennequin invariant of a Legendrian submanifold of the standard
contact $\Bigl (\R^{2n+1}=(x_1, \cdots, x_n, y_1, \cdots, y_n, z);
\ker (\sum y_i dx_i-dz)\Bigr )$ derived by T.~Ekholm, J.~Etnyre and
M.~Sullivan~\cite[Proposition 3.3 and Lemma 3.4]{EES}. The formulas
in~\cite{EES} compute the Bennequin invariant through the front
projection of the Legendrian submanifold to the $(x_1, \cdots, x_n,
z)$-subspace. They are sums over pairs of points of the front
projection with equal $(x_1, \cdots, x_n)$-coordinates (but
different $z$-coordinates) at which the tangent planes to the two
front branches are parallel. The increment of such a pair contains a
term which is the sign of the determinant of the Hessian of the
difference of the functions parameterizing the front branches. The
formulas in~\cite{EES} also involve the sign determined by whether
two orientations on a hyperplane induced by a such a pair of points
coincide or not. (This sign is hidden into an expressions that
involves the number of cusps of various types met along a path on a
front connecting the two points.)

These similarities are certainly expected since the Bennequin
invariant is the selflinking number, i.e. the linking number of the
Legendrian submanifolds and of its canonical infinitesimal shift.
However, since the Bennequin invariant is defined through the
classical linking number, it is not defined when the Legendrian
submanifold is nonzero homologous. In our work the submanifolds are
not zero homologous so there is no immediate relation between our
Theorem~\ref{hesse} and the results of~\cite{EES}.
\end{rem}
\forgotten

\m
\begin{ex}[calculation of $\sigma(\gl, \tau_0)$] Consider the passage through
a dangerous tangency point in a positively oriented chart $(x_1,
\ldots, x_m)$ shown in Figure~\ref{definesign.fig}. Assume that the
tangency in the Figure happens along the $(x_1, \ldots,
x_{m-1})$-hyperplane and that the $x_m$-axis points to the right in
the Figure. Assume that $\un \gl_1^{\tau}$ is the ``left'' surface
in the Figure and that $\un \gl _2^{\tau}$ is ``right'' surface. The
vector field in the Figure is the unit vector field normal to the
branches of $\un \gl_1^{\tau}$ and $\un \gl _2^{\tau}$ on which the
evaluation of the covector field $\theta^{\tau}_s, s\in
S^{m-1}_1\sqcup S^{m-1}_2,$ is positive. Then $\alpha=-1,$ $\sign
\det \Hess f({\bf 0})=1,$ and thus $\sigma (\gl, \tau_0)=-\eps.$
That is $\sigma(\gl, \tau_0)=-1$ if the two tangent branches induce
the same orientation on the common tangent $(m-1)$-hyperplane and
$\sigma(\gl, \tau_0)=+1$ otherwise.

\begin{figure}[htbp]
\begin{center}
\epsfxsize 10cm \hepsffile{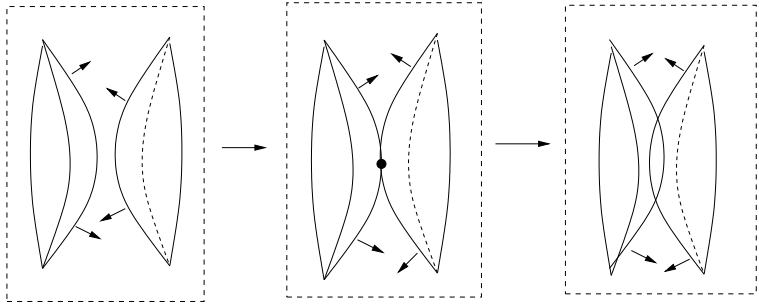}
\end{center}
\caption{Dangerous tangency}\label{definesign.fig}
\end{figure}
\end{ex}

\section{Examples}\label{examplescausality}

To illustrate the usage of the affine linking invariant consider the
following examples.

\begin{ex}\label{example1}
Let us show how one can use $\alk$ to determine that two events are
causally related. Let $(\ss^{m+1}, g), m>1$ be a globally hyperbolic
space-time and let $M^m$ be a Cauchy surface in
$X$. For brevity we denote $W_{a,M}$ by $W_{a}$ and we denote $\wt
W_{a, M}$ by $\wt W_a, a \in \ss$. Let $x, y$ be two events that do
not lie on a common null geodesic. From
Definition~\ref{df:linkskies} and Remarks~\ref{goodremark} it
follows that $x,y$ do not lie on a common null geodesic if and only
if the lifted wave fronts $(\wt W_{x}, \wt W_{y})$ form a
nonsingular link in $ST^*M.$

To compute the value of $\alk(\wt W_{x}, \wt W_{y})\in \A(M)$ we
take a generic homotopy deforming a trivial link $(\eps_{u},
\eps_{v}),$ $u\neq v\in M,$ to $(\wt W_{x}, \wt W_{y}).$ Let $p$ and
$n$ be the number of positive and negative crossings of
$\Sigma_0\subset \Sigma$ under the homotopy. Then $\alk(\mathfrak
S_x, \mathfrak S_y)= \alk(\wt W^{t}_{x}, \wt W^t_{y})=q(p-n)\in
\A(M),$ for the homomorphism $q:\bor_0(\B)=\Z\to \Z.$ If $\alk(\wt
W_{x}, \wt W_{y})=\alk(\mathfrak S_x, \mathfrak S_y)\neq 0\in
\A(M),$ then we conclude that $x$ and $y$ are causally related, see
\theoref{thmfrontsunlinked} and \theoref{welldefined}.

Observe that this computation and conclusion can be made just from
the shape of the cooriented and oriented fronts $W_{x}, W_{y}$ on a
Cauchy surface $M$, without the knowledge of the event points $x, y$
and of the Lorentz metric $g$ on $\ss$. Moreover, if $M$ is not
homeomorphic to an even dimensional sphere $S^{2k},$ then one does
not have to equip the pictures of the fronts with orientations. This
is since for such manifolds a positively oriented $S^{m-1}$-fiber of
$\pr:ST^*M\to M$ is not free homotopic to a negatively oriented
fiber $S^{m-1}$, see Theorem~\ref{reflection}. Thus if $M$ is not
homeomorphic to an even dimensional sphere, then the orientation of
the cooriented wave front $W_x$ on $M$ is always the one such that
the lifted wave front with this orientation is homotopic to a
positively oriented fiber $S^{m-1}$ of $\pr:ST^*M\to M.$

As an example of the computation, consider a globally hyperbolic
$(\ss,g)$ such that its Cauchy surface $M$ is not homeomorphic to a sphere.
Thus in this case the orientation of the fronts does not have to be
included into their description and $\A(M)\neq 0,$ see
Proposition~\ref{usefulprop} and Theorem~\ref{Main}.

Let $(W_{x}, W_{y})$ be two wave fronts located in a chart
diffeomorphic to $\R^m.$ Assume that for some vector $\vec v\in
\R^m$ the straight line homotopy $h_{\tau, \vec v}=(W_{x}+\tau
\vec v, W_{y}), \tau\in [0, +\infty)$ separates the fronts to be
located in two different halfspaces of $\R^m$. Assume moreover that
this homotopy involves exactly one passage through a dangerous
tangency point and this tangency point is nondegenerate, see for
example Figure~\ref{example1causality.fig}. Then by
Theorem~\ref{hesse} and the discussion before it we have $\alk(\wt
W_{x}, \wt W_{y})=\alk(\mathfrak S_x, \mathfrak S_y)=\pm 1\neq 0\in
\A(M).$ Here the sign $\pm 1$ depends on the sign of the determinant
of the Hessian at the dangerous tangency point and on the
coorientations and the actual orientations of the fronts. Hence the
events $x$ and $y$ are causally related.

\begin{figure}[htbp]
\begin{center}
\epsfxsize 11cm \hepsffile{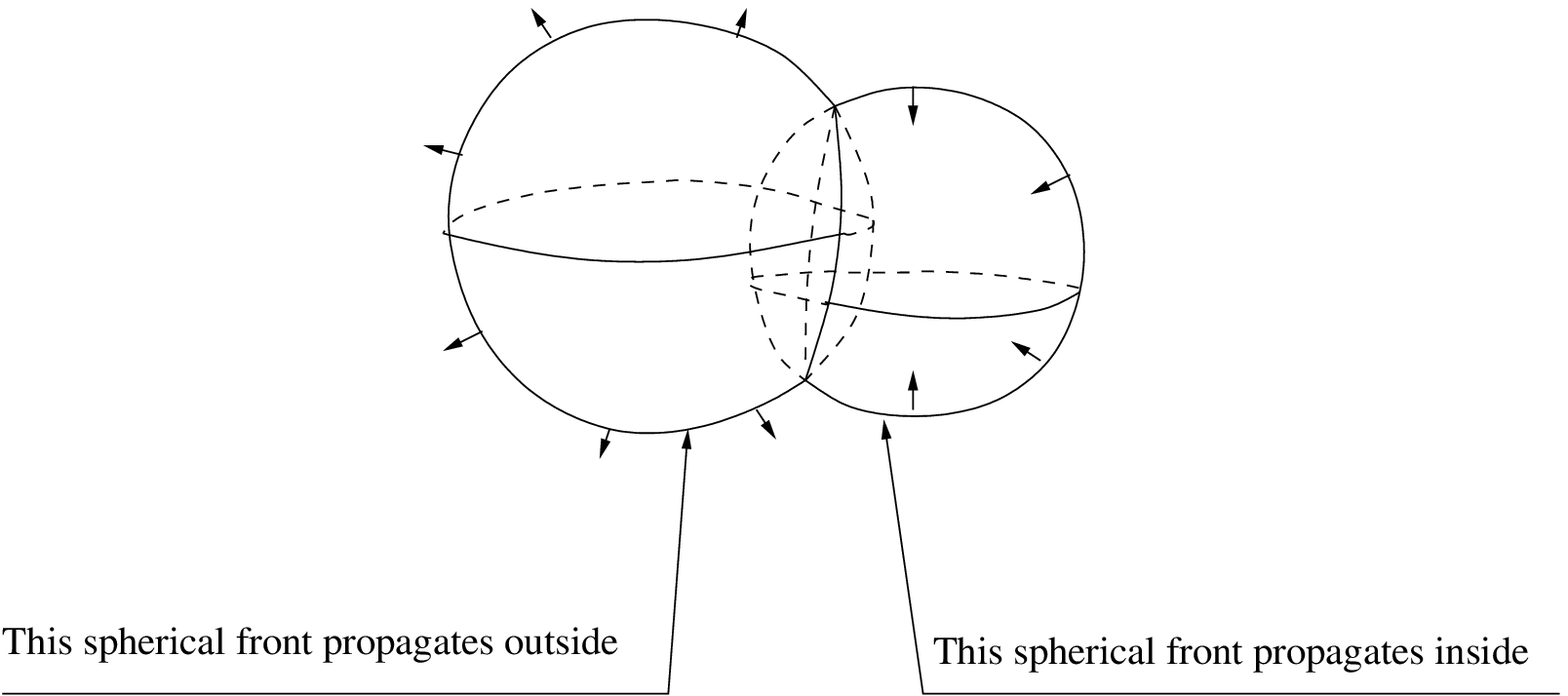}
\end{center}
\caption{}\label{example1causality.fig}
\end{figure}

\end{ex}

\begin{ex}\label{light}
Let us show how one can use $\alk$ to estimate the number of times
the exponent of the future directed null cone of a point $x$ crossed
a generic timelike curve joining two points $y, z$. {\em This number
can be interpreted as the number of times that an observer traveling
from $y$ to $z$ along a generic timelike curve sees the light from
the event $x$}.

Let $(\ss^{m+1}, g)$ be a globally hyperbolic space-time of
dimension $>2$ such that its Cauchy surface $M^m$ is not an
odd-dimensional rational homology sphere with finite $\pi_1(M).$
Theorem~\ref{Main} says that $\A(M)=\Z$ and $q:\Z\to \A(M)$ is the
identity map. Assume moreover that $M$ is not an even dimensional
homotopy sphere, so that as we discussed in~\ref{example1}, we do
not have to specify the orientations of the fronts $W_{x, M}$ when
depicting them.

Let $y, z\in \ss$ be two points that can be joined by a future
directed generic timelike curve from $y$ to $z.$ Let $L\ni y$ and
$N\ni z$ be two Cauchy surfaces.

Assume that $\Im W_{x,L}$ and $y\in L$ are in the same chart of $L$
and are shown in Figure~\ref{example2causality.fig}.a. Assume that
$\Im W_{x,N}$ and $z\in N$ are in the same chart of $N$ and are
shown in Figure~\ref {example2causality.fig}.b.
(Figure~\ref{example2causality.fig}.a depicts a trivially embedded
sphere with $y$ outside of it.
Figure~\ref{example2causality.fig}.b depicts a sphere
that can be obtained from the trivially embedded sphere located far
from $z$ by passing three times through a point $z$ and by creation
of some singularities far away from $z$.) The normal vector fields
to the fronts in Figure~\ref{example2causality.fig}.a and in
Figure~\ref{example2causality.fig}.b are such that the evaluations
of the covector fields defining the front lifts to $ST^*L$ and to
$ST^*N$ on the vector fields are positive. That is, these are the
vector fields defining the front lifts to $STL$ and to $STN$ that
are identified with $ST^*L$ and with $ST^*N$ via the Riemannian
metrics $g|_L$ and $g|_N.$

\begin{figure}[htbp]
\begin{center}
\epsfxsize 10cm \hepsffile{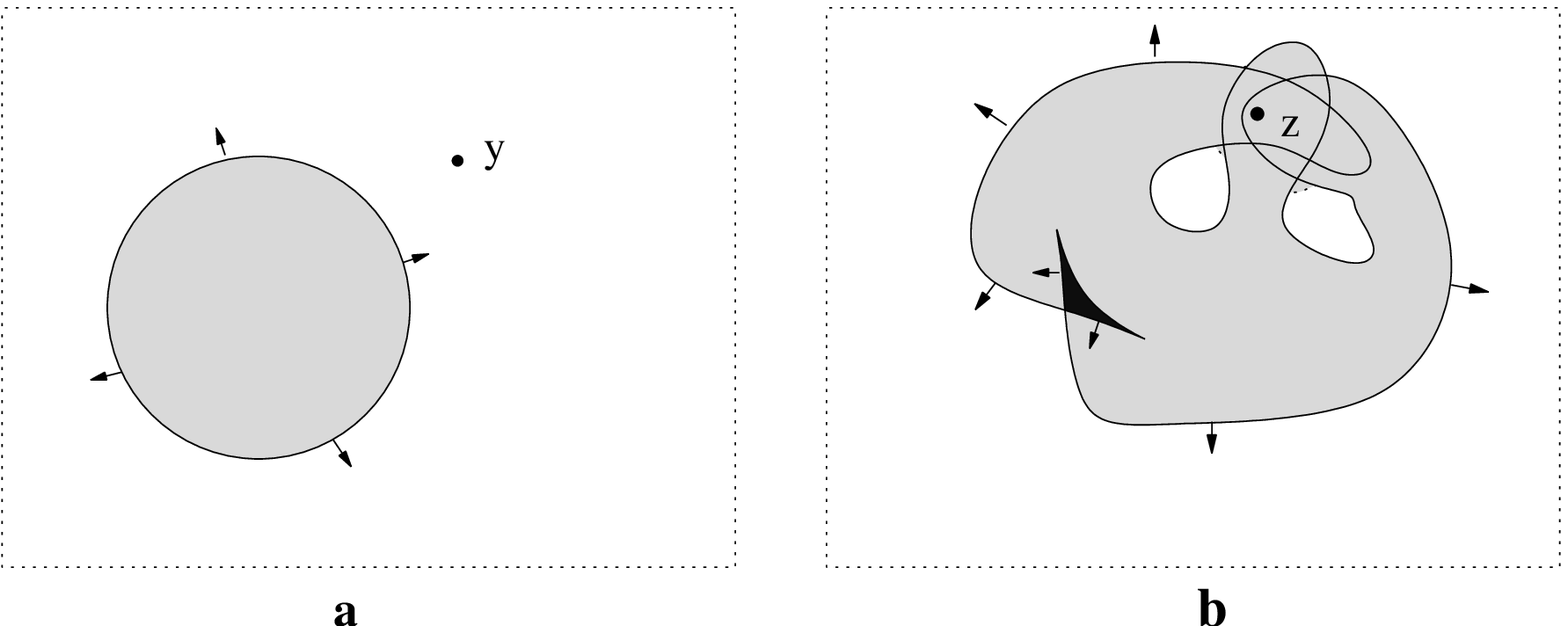}
\end{center}
\caption{}\label{example2causality.fig}
\end{figure}

Using Lemma~\ref{helpful} we get that $\alk(\mathfrak S_x, \mathfrak
S_y)=0$ and $\alk(\mathfrak S_x, \mathfrak S_z)=3.$ Let $\gamma$ be
a generic (as defined before \theoref{alkasintersection}) past
inextendible future directed curve ending at $y.$ Let $U\subset
T_x\ss$ be the part of the future pointing null hemicone where
$\exp_x$ is defined. Theorem~\ref{alkasintersection} says that
$\exp_{x}|_U\bullet \gamma=0\in \Z,$ where $\bullet$ is the
intersection number.

Let $\beta$ be a generic future directed timelike curve from $y$ to
$z.$ Then $\beta\cdot\gamma$ is a generic future directed past
inextendible curve ending at $z$. Theorem~\ref{alkasintersection}
says that $\exp_{x}|_U\bullet (\beta\cdot \gamma)=3$. Combining this
with equality $\exp_{x}|_U\bullet \gamma=0$, we conclude that
$\exp_{x}|_U\bullet \beta=3$. Thus an observer traveling from $y$
to $z$ along $\gb$ sees the light from the event $x$ at least $3$
times regardless of which generic timelike curve s/he chooses to
travel. (If $\beta$ is not generic and the points of
self-intersection of $\exp_{x}|_U$ belong to $\Im \gb$ then, at a
point of $\beta$, s/he might see the light coming from several
different directions, and the total number of times s/he sees light
may be less than $3.$)
\end{ex}

\section{Refocussing and nonrefocussing spaces\/}\label{s:refocussing}

In this section we discuss (non)refocussing space-times. We use
them in the next section.

Let $\SKY$ denote the set of all skies in $(\ss, g)$ with the following
topology. For $\wh {\mathfrak S}\in \SKY$ the topology base at $\wh {\mathfrak 
S}$ is given 
by $\{\mathfrak S| \mathfrak S\subset W\}$ for open $W\subset \mathcal N$ such 
that $\wh {\mathfrak S}\subset W.$

Consider the map
\begin{equation}\label{eq:mu}
\mu: \ss\to \SKY, \quad \mu(x)=\mathfrak S_x.
\end{equation}

One verifies that if $(\ss, g)$ is globally hyperbolic, then $\mu$ is 
continuous. 
Example~\ref{ex:sameskies} shows that $\mu$ is not even a bijection in
general. Is $\mu$ a homeomorphism provided that it is a bijection? In order for 
$\mu$ to be open it suffices to show that for every $x\in \ss$ and open $U\ni x$ 
there exists an open $V\ni x$ contained in $U$ such that 
$\mu(V)$ is open. This motivates the following definition, 
cf.~\cite{LowNullgeodesics}.

\begin{defin}\label{d:nonrefocussing}
A strongly causal space-time $(\ss,g)$ (that is not necessarily
globally hyperbolic) is called {\it refocussing at $x\in X$\/} if
there exists a neighborhood $O$ of $x$ with the following property:
For every open $U$ with $x\in U \subset O$ there exists $y\not \in
U$ such that all the null-geodesics through $y$ enter $U.$ A
space-time $(X,g)$ is called {\it refocussing\/} if it is
refocussing at some $x,$ and it is called {\it nonrefocussing\/} if
it is not refocussing at every $x\in X.$
\end{defin}

Low~\cite{LowNullgeodesics} introduced the concept of nonrefocussing space-times 
and observed that if a globally hyperbolic $(\ss, g)$ is nonrefocussing, then 
$\mu$ is bijective and open, i.e. $\mu: \ss \to \SKY$ is a homeomorphism. We 
note that in the original definition of Low~\cite{LowNullgeodesics} $U$ was 
allowed to be any open neighborhood containing $x$ that is not necessarily 
sufficiently small. This is clearly a typo, since for every $(\ss, g)$ and 
$U=\ss$ such a point $y\not \in U$ does not exist.

We need the following topological lemma.

\begin{lem}\label{homology}
Let $M^m$ be a non-compact manifold and $B$ an open ball in $M$.
Assume that the closure $\ov B$ of $B$ is a smoothly embedded ball.
Let $V$ be an open subset in $M$ such that its closure $\ov V$ is
compact and $\ov V \setminus V\subset B$. Then $\ov V \subset B$.
\end{lem}

\pp Without loss of generality we assume that $V$ is connected.
Since $\ov V\setminus V$ is compact, there exists an open disk
$B_0\subset B$ such that
$$
\ov V \setminus V\subset B_0\subset \ov B_0\subset B
$$
and the boundary $\ov B_0\setminus B_0$ of $B_0$ is a smoothly
embedded $(m-1)$-dimensional sphere $S$. We may assume that $\ov
V\cap S\neq\emptyset.$ Otherwise since $\ov V$ is connected, $\ov
V\subset B_0\subset B$ and the proof is finished. Furthermore, the
set $V\cap S= \ov V\cap S$ is open as well as closed in $S$. Since
$S$ is connected and $\ov V\cap S\ne\emptyset$, we conclude that $S
\subset V$.

Arguing by contradiction, suppose that $\ov V\setminus B \ne
\emptyset$. Since $\ov V\setminus V\subset B,$ we have that
$V\setminus B\neq \emptyset.$ Then $V\setminus \ov B_0\neq
\emptyset$ is an open subset of $M$ and $Y:=V\setminus B_0=\ov
V\setminus B_0$ is a compact connected smooth orientable manifold
with the interior $\Int Y=V\setminus \ov B_0.$ Hence $\partial
Y=\ov B_0 \setminus B_0=S$.

Take a point $a\in Y \setminus \partial Y$ and consider the
commutative diagram
$$
\CD
(Y, \partial Y) @>>> (M,B)\\
@VVV @VVV\\
(M, M\setminus a) @= (M, M\setminus a)
\endCD
$$
of inclusions. This diagram induces the commutative diagram
$$
\CD
\Z @= H_m(Y, \partial Y) @>>> H_m(M,B) @= 0\\
@. @VVV @VVV @.\\
\Z @= H_m(M, M\setminus a) @= H_m(M, M\setminus a) @=\Z
\endCD
$$
Here $H_m(M,B)=H_m(M)=0$ since $M$ is not compact. Since the left
map is an isomorphism, we conclude that such a diagram cannot exist.
Thus, $V\subset B$. \qed

\begin{defin}
A set $A$ in a (not necessarily globally hyperbolic) space-time
$(X,g)$ is {\it achronal\/} if no timelike curve intersects $A$ more
than once. In particular every subset of a Cauchy surface is
achronal.

For an achronal set $A$ its {\em future Cauchy
development\/} $D^+(A)$ is the set of all the points $x\in X$ such
that every past inextendible non-spacelike curve through $x$ meets
$A.$ Similarly, the {\it past Cauchy development\/} $D^-(A)$ is the
set of all $x\in X$ such that every future inextendible
non-spacelike curve through $x$ meets $A.$ In particular $A$ is a
subset of both $D^+(A)$ and $D^-(A).$ The {\it Cauchy development of
$A$\/} is $D(A)=D^+(A)\cup D^-(A).$

If $M$ is a Cauchy surface in a globally hyperbolic space-time
$(\ss, g),$ then $\ss=D^+(M)\cup D^-(M).$
\end{defin}

\begin{prop}[Low~\cite{LowNullgeodesics, 
LowRefocussing}]\label{p:nonrefocussing}
A globally hyperbolic space-time $(\ss, g)$ with a non-compact
Cauchy surface $M,$ is nonrefocussing.
\end{prop}

\pp A brief outline of the proof is contained in \cite[Theorem
5]{LowRefocussing}. We are grateful to Robert Low who explained us
the details of his proof.

Assume that $(\ss, g)$ is refocussing at a point $x.$ Take an open
neighborhood $O$ of $x$ such that for every open $V$ with $x\in
V\subset O$ there exists $y\not \in V$ such that all the null
geodesics through $y$ enter $V.$

Take a Cauchy surface $M$ through $x$ and an open ball $B$ in $M$
with $x\in B$ such that the closure $\ov B$ is a smoothly embedded
ball. Put $U=D(B).$ Then $U$ is open, globally hyperbolic and
contains $B$, see~\cite[Section 14, Lemma 42 and Lemma 43]{ONeill}.
Clearly $B$ is a Cauchy surface of $U.$ Assume moreover that $B$ is
sufficiently small so that $U\subset O.$

Take a point $y\in \ss$
with $y\not \in U=D(B)$ such that all the null-geodesics through $y$
cross $U$. Without loss of generality we assume that $y\in D^+(M).$
By~\cite[Proposition 3.16 and Lemma 3.5]{BeemEhrlichEasley}, the
set $J^-(y)$ is closed and the set $I^-(y)$ is open. Moreover,
$J^-(y)$ is the closure of $I^-(y)$ by~\cite[Section 14, Lemma
6]{ONeill}. Put $J^-=J^-(y)\cap M$ and $I^-=I^-(y)\cap M$. Because
of what we said above, $J^-$ is the closure (in $M$) of the open
subset $I^-$ of $M$. Since $J^-(y)\cap D^+(M)$ is compact
by~\cite[Section 14, Lemma 40]{ONeill}, we get that $J^-$ is compact
in $M.$

By~\cite[Corollary 4.14]{BeemEhrlichEasley} if $z\in J^-(y)\setminus
I^-(y),$ then there is a null-geodesic from $y$ to $z.$ Thus if
$z\in J^- \setminus I^-,$ then $z$ lies on a past directed
null-geodesic from $y.$ By our choice of $y$ this null geodesic has
to pass through $U.$ Since $B$ is a Cauchy surface of a globally
hyperbolic $U,$ this null geodesic crosses $M$ in some point of
$B\subset M.$ Thus all the points of $J^-\setminus I^-$ are in $B$,
i.e. $J^-\setminus I^-\subset B$.

By \lemref{homology} applied to the case $V=I^-$ we get the
inclusion $J^-(y)\cap M\subset B$. Thus $y\in D^+(B)\subset U$ and
we get a contradiction. \qed

\forget \m The proof of Proposition~\ref{p:nonrefocussing} implies
the following stronger statement.

\begin{prop}\label{Lowpropositionmodified}
Let $(\ss^{m+1}, g)$ be a globally hyperbolic space-time with a
noncompact smooth spacelike Cauchy surface $M.$ Let $B$ be an open
ball in $M$ such that the closure $\ov B$ is a smoothly embedded
closed ball and let $U=D(B)$ be the Cauchy development. Then there
are no points $y\not \in U$ such that all the null geodesics through
$y$ cross $U.$ We use Proposition~\ref{Lowpropositionmodified} to
get a contradiction
\end{prop}

We use Proposition~\ref{Lowpropositionmodified} to prove the
following result. \forgotten

\m Clearly if $p:(\ss_1, g_1)\to (\ss, g)$ is a Lorentz cover of a globally 
hyperbolic space-time and $(\ss_1, g_1)$ is refocussing, then $(\ss, g)$
is also refocussing. Below we prove the converse result.

\begin{thm}\label{refocussingtheorem}
Let $(\ss^{m+1}, g)$ be a globally hyperbolic space-time that is
refocussing, and let $p: \ss_1 \to \ss$ be a covering map. We equip
$\ss_1$ with the induced Lorentz metric $g_1$. Then $(\ss_1, g_1)$
is a refocussing globally hyperbolic space-time. In particular, if
$\ss$ has infinite fundamental group then $X$ is nonrefocussing,
see~Proposition~$\ref{p:nonrefocussing}$.
\end{thm}

\pp First, we prove that $(\ss_1, g_1)$ is globally hyperbolic. It
suffices to prove that $(\ss_1, g_1)$ admits a Cauchy surface.
Choose a Cauchy surface $M\subset X$ and put $M_1=p^{-1}(M)$. We
claim that $M_1$ is a Cauchy surface. Indeed, if $\gamma(t)$ is an
inextendible nonspacelike curve in $\ss_1,$ then $p\circ \gamma(t)$
is an inextendible nonspacelike curve in $\ss.$ Since $M$ is a
Cauchy surface, $p\circ \gamma(t)$ crosses $M$ at exactly one value
of $t.$ Hence $\gamma(t)$ also crosses $M_1$ at exactly one value of
$t,$ and thus $M_1$ is a Cauchy surface.

Now, suppose that $X$ is refocussing at some $x\in \ss$. Take a
Cauchy surface $M$ in $(\ss, g)$ with $x\in M$ and consider the
Cauchy surface $M_1=p^{-1}(M)$ in $(\ss_1, g_1)$. Choose $x_1\in
M_1$ such that $p(x_1)=x.$ Choose an open ball $B'_1 \subset M_1$
that is a normal neighborhood of $x_1$ with respect to the
exponential map $\exp^r_{x_1}:T_{x_1}M_1\to M_1$ constructed using
the Riemannian metric $g^r_{M_1}$ induced on $M_1$ from $g_1.$
Without loss of generality we can assume that $p|_{B'_1}:B'_1\to
M$ is an embedding. Choose an open ball $B_1\subset B'_1$ such that
the closure $\ov B_1$ is a smoothly embedded closed ball contained
in $B'_1.$ Put $U_1=D(B_1)$ to be the Cauchy development of $B_1.$
Then $U_1$ is open globally hyperbolic and contains $B_1,$
see~\cite[Section 14, Lemma 42 and Lemma 43]{ONeill}. Clearly $B_1$
is a Cauchy surface for $U_1.$

Put $B=p(B_1)\subset M$ to be the open ball containing $x.$ Put
$U=p(U_1)\ni x$. Since $p$ is a cover, $U$ is open. Clearly $B$ is a
Cauchy surface for $U$ and hence $U$ is globally hyperbolic.

Let $O$ be a neighborhood of $x$ described in
\defref{d:nonrefocussing}. It is not difficult to prove that the ball
$B_1$ can be chosen so that $U \subset O$. Hence, there exists
$y\not \in U$ such that all the null geodesics through $y$ cross
$U.$ Without loss of generality $y\in D^+(M).$

Choose an $M$-proper isometry $h:M\times \R\to \ss$ and put $(m_y,
t_y)\in M\times \R$ to be the point such that $h(m_y, t_y)=y.$
Define $F:S^{m-1}\times \R\to M\times \R$ via $F(s, t)=(W^t_{y,
M}(s), t).$ For $s\in S^{m-1}$ put $\gamma_s(t)=F(s,t).$ Clearly up
to reparameterization the curves $h\circ \gamma_s(t), s\in S^{m-1},$
are exactly all the null geodesics through $y.$ Also
$h(\gamma_s(0))\in B$ is exactly the intersection point of the
corresponding null geodesic with $B$ and $h(\gamma_s(t_y))=y$ for
all $s\in S^{m-1}.$

Put $B'=p(B'_1).$ For $s\in S^{m-1}$ put $\rho_s:[0,1]\to B'$ to be
the unique geodesic (with respect to the induced Riemannian metric
on $B'$) arc from $x\in B\subset B'$ to $h(\gamma_s(0))\in B\subset
B'.$ For $s\in S^{m-1}$ define the path $\delta_s:[0, 1+t_y]\to \ss$
from $x$ to $y$ via $\delta_s(t)=\rho_s(t)$ for $t\in [0,1]$ and
$\delta_s(t)=h(\gamma_s(t-1))$ for $t\in [1, 1+t_y].$

For every $s_0, s_1\in S^{m-1}$ the paths $\delta_{s_0}$ and
$\delta_{s_1}$ are homotopic relative boundary. The homotopy is
given by the family of paths $\delta_{\beta(\tau)}$ constructed from
a path $\beta:[0,1]\to S^{m-1}$ with $\beta(0)=s_0, \beta(1)=s_1.$

For $s\in S^{m-1}$ put $\delta_{1, s}:[0, 1+t_y]\to \ss_1$ to be the
lift of $\delta_s$ starting at $x_1.$ Since all the paths $\delta_s$
are homotopic relative boundary, we get that all the values
$\delta_{1, s}(1+t_y)\in \ss_1$ are equal and we put
$y_1=\delta_{1,s}(1+t_y).$

Since $y\not \in U$, we conclude that $y_1\not \in U_1.$ We claim
that all the null geodesics through $y_1$ pass through $U_1.$ Indeed
for every $s\in S^{m-1},$ the path $\delta_{1, s}|_{[1,1+t_y]}$ is (up
to reparameterization) an arc of the null geodesics through $y_1$
and $\delta_{1, s_1}(1)\in B_1\subset U_1$. Thus, $(\ss_1, g_1)$ is
refocussing. \qed

\begin{rem}[Refocussing space-times and the Blaschke conjecture type 
problems]\label{Blaschke} 
The following construction gives many examples of refocussing globally 
hyperbolic space-times.
Let $(M, \ov g)$ be a complete oriented Riemannian manifold, such that
for some $x\in M$ and positive $r\in \R$ the exponential $\exp_x:T_xM\to M$
maps the whole sphere of radius $r$ centered at ${\bf 0}\in T_xM$ to
one point. A static Lorentz manifold $(M\times \R, \ov g\oplus -dt^2)$ is 
globally hyperbolic, see~\cite[Theorem 3.66]{BeemEhrlichEasley}, and it is
clearly refocussing at $(x, r).$

One can show that $x$ is the end point of all the length $2r$ geodesic arcs in 
$M$ starting at $x,$~i.e. $(M, \ov g)$ is a $Y_{2r}^x$-manifold 
in terms of Besse~\cite[Chapter 7.B]{Besse}. The question on topology of such 
manifolds is closely related to the Blaschke conjecture type 
problems,~see~\cite{Besse}. A weak form of a Bott-Samelson Theorem says that 
every $Y_{2r}^x$ manifold is 
a closed manifold with finite $\pi_1$ whose rational cohomology ring is 
generated by one element, see~\cite{BerardBergery},~\cite[Theorem 7.37]{Besse}, 
cf.~\cite{Bott, Samelson}.

Clearly there are many examples of refocussing globally hyperbolic space-times 
that are not obtained by the above construction. However 
Theorem~\ref{refocussingtheorem} says that a Cauchy surface in all of them is a 
closed manifold with finite $\pi_1.$ It would be interesting to know if its 
rational cohomology ring is necessarily generated by one element, i.e. if the 
Bott-Samelson type result holds for a Cauchy surface of a refocussing globally 
hyperbolic space-time. Since the only oriented two-dimensional surface with 
finite $\pi_1$ is $S^2$, Theorem~\ref{refocussingtheorem} implies that this is 
indeed so for $(2+1)$-dimensional globally hyperbolic refocussing space-times.
\end{rem}

\section{A weakened Low conjecture is true.\/}\label{s:weakened}

We show that a certain weakened version of the Low conjecture holds
for a vast family of globally hyperbolic space-times $(\ss^{m+1},
g), m>1$. 

Natario and Tod~\cite[Figure 13, p. 18]{NatarioTod} considered 
$(2+1)$-dimensional space-times with a Cauchy surface diffeomorphic to $\R^2$ 
and presented several examples of causally related events whose skies are linked 
but have zero linking number. They also observed
that since the skies of events are
Legendrian submanifolds of $\mathcal N,$ it makes sense to
ask if the skies of two causally related events are always
nontrivially linked in the Legendrian sense. When a Cauchy surface $M$ is 
diffeomorphic to an open subset of $\R^m$,
this is the modified Low conjecture due to Natario and
Tod~\cite{NatarioTod}. 

However even for $(2+1)$-dimensional space-times not all of the
Legendrian embeddings $S^{m-1}\to ST^*M=\mathcal N$ that are
Legendrian isotopic to $\eps_v, v\in M,$ correspond to skies,
see~\cite[Theorem 4.5]{NatarioTod}. Thus one can weaken the Low
conjecture even further and ask if it is always true that the skies
of causally related events in $(\ss, g)$ can not be unlinked by an
isotopy through the skies of events in $(\ss, g).$

\begin{defin}[isotopy through skies]\label{isotopythroughskies}
Let $(\ss^{m+1}, g), m+1>2,$ be a globally hyperbolic space-time. We
say that two nonsingular links $(\mathfrak S_1, \mathfrak S_2)$ and
$(\mathfrak S_1', \mathfrak S_2')$ are {\it isotopic through
skies\/} if there exists a continuous map $\rho:[0,1]\to
\mathcal{SKY}\times \mathcal{SKY}, \rho(t)=(\rho_1(t), \rho_2(t))$
such that $\rho_i(0)=\mathfrak S_i$, $\rho_i(1)=\mathfrak S_i'$,
$i=1,2$ and for all $t\in [0,1]$ the intersection of the skies
$\rho_1(t)$ and $\rho_2(t)$ in $\mathcal N$ is empty.
\end{defin}

\begin{defin}[sky-isotopy]\label{skyisotopy}
Let $x_1, x_2, y_1, y_2\in \ss$ be such that neither $x_1, x_2$ nor
$y_1, y_2$ belong to a common null geodesic. We say that the pairs
$({x_1}, {x_2})$ and $({y_1}, {y_2})$ are {\it sky-isotopic\/} if
there exist paths $p_1, p_2:[0,1]\to \ss$ such that $p_i(0)=x_i,
p_i(1)=y_i, i=1,2,$ and the skies $\mathfrak S_{p_1(t)}$ and
$\mathfrak S_{p_2(t)}$ are disjoint, for all $t\in [0,1].$ (The last
condition is equivalent to requiring that for every $t\in [0,1]$ the
points $p_1(t), p_2(t)$ do not belong to a common null geodesic.)
\end{defin}

\begin{rem}[Comparison of the ``sky-isotopy'' and of
the ``isotopy through skies'' notions]\label{comparison} If $({x_1},
{x_2})$ and $({y_1}, {y_2})$ are sky-isotopic, then, clearly, the
links $(\mathfrak S_{x_1}, \mathfrak S_{x_2})$ and $(\mathfrak
S_{y_1}, \mathfrak S_{y_2})$ are isotopic through skies. Indeed
given the paths $p_1(t), p_2(t)$ as in the definition of
sky-isotopy, put $\rho_i(t)=\mathfrak S_{p_i(t)}, t\in [0,1], i=1,2.$

\m The converse is not true in general, the pairs $(x,y)$ and $(x',y)$ of events 
in Example~\ref{einstein} yield a counterexample.

\m For $(X,g)$ nonrefocussing, $({x_1}, {x_2})$ and $({y_1},
{y_2})$ are sky-isotopic if and only if $(\mathfrak S_{x_1},
\mathfrak S_{x_2})$ and $(\mathfrak S_{y_1}, \mathfrak S_{y_2})$ are
isotopic through skies. Indeed given $\rho_i:[0,1]\to \mathcal{SKY},
i=1,2,$ as in
\defref{skyisotopy}, put $p_i(t)=\mu^{-1}(\rho_i(t)),$ where $\mu:X\to
\mathcal{SKY}$ is the homeomorphism from \eqref{eq:mu}.
\end{rem}

The following Theorem~\ref{Low1true} says that any two
pairs of causally unrelated events in a globally hyperbolic $(\ss,
g)$ are sky-isotopic, and that no such pair is sky-isotopic to a pair
of causally related events.

\begin{thm}\label{Low1true}
Let $(\ss^{m+1}, g), m+1>2$ be a globally hyperbolic space-time. Let
$(x_1, x_2)$ be a pair of causally unrelated
events, and let $(y_1, y_2)$ be two events that do
not belong to a common null geodesic. Then the following two statements are 
equivalent:
\begin{description}
\item [1] The events $y_1$ and $y_2$ are not
causally related.
\item [2] The pairs $({x_1}, {x_2})$ and
$({y_1}, {y_2})$ are sky-isotopic.

\end{description}
\end{thm}

\pp Choose a Cauchy surface $M\subset X$ and an $M$-proper isometry
$h:M\times \R\to X.$

{\it The proof of the implication $1\implies 2$ \/} follows immediately from
the following three claims that are proved below:

{\bf Claim 1.} For any causally unrelated $v_1, v_2$ there exist
$t\in \R$ and $w_1, w_2\in M_t\subset X$ such that $({v_1}, {v_2})$
is sky-isotopic to $({w_1}, {w_2}).$

{\bf Claim 2.} If $(v_1, v_2)$ and $(w_1, w_2)$ are two pairs of
distinct events in the same Cauchy surface $M_{\tau}\subset X,$ then
$({v_1}, {v_2})$ is sky-isotopic to $({w_1}, {w_2}).$

{\bf Claim 3.} For $t_1\neq t_2\in \R, n_1\neq n_2\in M$ the pairs
of events $\bigl ({h(n_1, t_1)}, {h(n_2, t_1)}\bigr)$ and $\bigl
({h(n_1, t_2)}, {h(n_2, t_2)}\bigr)$ are sky-isotopic.

\m We prove Claim 1. Let $t_1, t_2\in \R$ be such that $v_i\in
M_{t_i}, i=1,2.$ Without loss of generality we assume that $t_1\leq
t_2.$ Let $\gamma$ be a future directed inextendible timelike curve
through $v_2.$ Reparameterize $\gamma$ so that $\gamma(t)\in
M_t\subset X$ for all $t\in \R.$
Since $v_1$ and $v_2$ are causally unrelated, we conclude that $\mathfrak
S_{v_1}\cap \mathfrak S_{\gamma(t)}= \emptyset$ for all $t\in [t_1,
t_2]$. Indeed, if $\mathfrak S_{v_1}\cap \mathfrak
S_{\gamma(\tau)}\neq \emptyset$ for some $\tau\in [t_1, t_2],$ then
the arc of a null geodesic $\nu\in \mathfrak S_{v_1}\cap \mathfrak
S_{\gamma(\tau)}$ from $v_1$ to $\gamma(\tau)$ followed by
$\gamma|_{[\tau, t_2]}$ is a future directed non-spacelike curve
from $v_1$ to $v_2.$ Put $w_1=v_1, w_2=\gamma(t_1)\in M_{t_1}.$ Now,
to see that $({w_1}, {w_2})=({v_1}, {\gamma(t_1)})$ is sky-isotopic
to $({v_1}, {\gamma(t_2)})=({v_1}, {v_2})$, put $p_1(t)=v_1,
p_2(t)=\gamma(t), t\in [t_1, t_2]$.

\m We prove Claim 2. Let $(v_1, v_2)$ and $(w_1, w_2)$ be two pairs
of distinct events in the same Cauchy surface $M_{\tau}.$ Since
$\dim (\ss)>2$ and hence $\dim M_{\tau}>1$, we can choose two paths
$p_1(t), p_2(t)$ in $M_ {\tau}, t\in [0,1]$ such that $p_i(0)=v_i,$
$p_i(1)=w_i, i=1,2$ and $p_1(t)\neq p_2(t)$ for all $t\in [0,1].$
Since any two distinct points in the same Cauchy surface are
causally unrelated, $\mathfrak S_{p_1(t)}\cap \mathfrak
S_{p_2(t)}=\emptyset,$ for all $t\in [0,1].$ Thus, $({v_1}, {v_2})$
is sky-isotopic to $({w_1}, {w_2})$.

\m We prove Claim 3. Assume without loss of generality that
$t_1<t_2.$ Put $p_1(t)=h(n_1, t),$ $p_2(t)=h(n_2, t), t\in [t_1,
t_2].$ Since $p_1(t), p_2(t)\in M_t,$ the events $p_1(t)$ and
$p_2(t)$ are causally unrelated, for $t\in [t_1, t_2].$ Hence
$\mathfrak S_{p_1(t)}\cap \mathfrak S_{p_2(t)}=\emptyset,$ for all
$t\in [t_1, t_2],$ and $({h(n_1, t_1)}, {h(n_2, t_1)})$ is
sky-isotopic to $({h(n_1, t_2)}, {h(n_2, t_2)}).$ This completes the
proof of Claim 3 and, hence, of the implication $1\implies 2$ of the Theorem.

\m {\it To prove the implication {\bf $2\implies 1$}\/}, recall the
notion of Lorentzian distance, see~\cite{BeemEhrlichEasley}. For
points $p,q$ in a (not necessarily globally hyperbolic) space-time
$(X,g)$ with $q\in J^+(p)$ put $\Omega_{p,q}$ to be the space of all
piecewise smooth future directed non-spacelike curves
$\delta:[0,1]\to X$ with $\gamma(0)=p, \gamma(1)=q.$ For $\delta\in
\Omega_{p,q}$ choose a partition $0=t_0<t_1<t_2\cdots
<t_{n-1}=t_n=1$ such that $\delta |_{(t_i, t_{i+1})}$ is smooth for
all $i\in\{0,1,\cdots(n-1)\},$ and define the {\em Lorentzian arc
length\/} $L(\delta)$ of $\delta$ by
$$
L(\delta)=L_g(\delta)=\sum_{i=0}^{n-1}\int_{t_i}^{t_{i+1}}\sqrt{-g(\dot
{\delta}(\tau), \dot {\delta}(\tau))}d\tau.
$$

For $p,q\in (X,g)$ define the {\em Lorentzian distance function \/}
$d=d_g:X\times X\to \R\sqcup\infty$ as follows: set $d(p,q)=0$ for
$q\not \in J^+(p);$ and set $d(p,q)=\sup\{L_g(\delta)|\delta\in
\Omega_{p,q}\},$ for $q\in J^+(p).$ By~\cite[Chapter 14, Corollary
1]{ONeill}, if $a<\!\!<b$ and $b\leq c,$ or if $a\leq b$ and
$b<\!\!<c,$ then $a<\!\!<c.$ Combining this with the definition of
the Lorentzian distance $d$ we get that $d(p,q)>0$ if and only if
$q\in I^+(p)$.

In general, the Lorentzian distance function is not continuous, and
$d(p,q)$ is not finite. (Also $d(p,q)\neq d(q,p)$ and $d(p,p)\neq 0$
in many cases.) However, for $(\ss, g)$ globally hyperbolic, $d$ satisfies 
finite distance condition and is a
continuous function on $\ss\times \ss$, see~\cite[Corollary
4.7]{BeemEhrlichEasley}.

{\it We argue by contradiction.\/}
Assume that $y_1, y_2$ are causally related, but the pairs $(x_1,
x_2)$ and $(y_1, y_2)$ are sky-isotopic. Since $y_1, y_2\in \ss$ are
causally related, either $y_1\in J^+(y_2)$ or $y_2\in J^+(y_1).$
Without loss of generality assume that $y_2\in J^+(y_1).$ If $y_2\in
J^+(y_1)\setminus I^+(y_1),$ then $y_1$ and $y_2$ lie on a common
null geodesic, see~\cite[Corollary 4.14]{BeemEhrlichEasley}. This
contradicts the Theorem assumptions. Hence $y_2\in I^+(y_1)$ and
$d(y_1, y_2)>0.$

Since $(x_1, x_2)$ and $(y_1, y_2)$ are sky-isotopic,
take $p_1, p_2:[0,1]\to \ss$ such that
$p_i(0)=y_i, p_i(1)=x_i, i=1,2,$ and such that $\mathfrak
S_{p_1(t)}\cap \mathfrak S_{p_2(t)}=\emptyset,$ for all $t\in
[0,1].$ Define a continuous function $\ov d:[0,1]\to \R$ by $\ov
d(t)=d(p_1(t), p_2(t)).$ We have $\ov d(0)=d(p_1(0), p_2(0))=d(y_1,
y_2)>0.$ Furthermore, $\ov d(1)=d(p_1(1), p_2(1))=d(x_1, x_2)=0,$
since $x_1, x_2$ are causally unrelated. Put $\tau=\inf \{t\in
[0,1]| \ov d(t)=0\}$, so that
\begin{equation}\label{positived}
\ov d (\tau)=0 \text{ and } \ov d(t)>0 \text{ for all } t<\tau.
\end{equation}
Below we show that $\mathfrak S_{p_1(\tau)}\cap \mathfrak
S_{p_2(\tau)}\neq \emptyset.$ This contradicts our assumptions about
$p_1, p_2.$

By~\cite[Proposition 6.6.1]{HawkingEllis} $(\ss, g)$ is {\it causally simple,\/} 
i.e.~the sets $J^{\pm}(K)=\cup_{k\in K}J^{\pm}(k)$ are
closed for every compact $K\subset \ss$.

By \eqref{positived}, $d(p_1(t), p_2(t))=\ov d(t)>0$ for all
$t<\tau.$ Hence $p_2(t)\in I^+(p_1(t))\subset J^+(p_1(t))$ for all
$t<\tau,$ and so $\Im (p_2|_{[t, \tau)})\subset J^+(\Im (p_1|_{[t,
\tau]}))$ for all $t<\tau.$ Since $\Im (p_1|_{[t, \tau]})$ is
compact and $(\ss, g)$ is causally simple, we conclude that $J^+(\Im
(p_1|_{[t, \tau]}))$ is closed, and hence $p_2(\tau)\in J^+(\Im
(p_1|_{[t, \tau]}))$ for all $t<\tau.$

Choose an increasing sequence $\{t_i\in [0,1]\}_{i\in \N}$ that
converges to $\tau.$ Then for each $i\in \N$ there exists $\wt
t_i\in [t_i, \tau]$ such that $p_2(\tau)\in J^+(p_1(\wt t_i)).$
Hence $p_1(\wt t_i)\in J^-(p_2(\tau))$ for all $i.$ Since $(\ss, g)$
is causally simple, $J^-(p_2(\tau))$ is closed and it contains the
point $\ds p_1(\tau)=\lim_{i\to \infty} p_1(\wt t_i).$ Since
$p_1(\tau)\in J^-(p_2(\tau)),$ we have $p_2(\tau)\in
J^+(p_1(\tau)).$

On the other hand, $p_2(\tau)\not \in I^+(p_1(\tau))$ since
$d(p_1(\tau), p_2(\tau))=\ov d(\tau)=0$ by \eqref{positived}. So,
$p_2(\tau)\in J^+(p_1(\tau))\setminus I^+(p_1(\tau))$, and therefore
the points $p_1(\tau), p_2(\tau)$ belong to a common null geodesic,
see~\cite[Corollary 4.14]{BeemEhrlichEasley}. Thus $\mathfrak
S_{p_1(\tau)}\cap \mathfrak S_{p_2(\tau)}\neq \emptyset.$
Contradiction. \qed

\begin{rem}
Looking carefully at the proof of the implication $2\implies 1$ of
Theorem~\ref{Low1true} one notices that in fact we proved the
following stronger statement. Let $(\ss^{m+1}, g)$ be a causally
simple space-time such that the Lorentzian distance on it is a
continuous function satisfying the finite distance condition. Let $(x_1, x_2)$ 
be a pair of causally unrelated events and let $(y_1, y_2)$ be a pair of 
causally related events. Then for every pair of continuous paths $p_i:[0,1]\to 
\ss$ such that $p_i(0)=x_i, p_i(1)=y_i, i=1,2,$ there exists $t\in [0,1]$
for which $p_1(t)$ and $p_2(t)$ belong to the common null geodesic.
\end{rem}

\m The following Corollary~$\ref{Low2true}$ can be viewed as the
proof of a weakened Low conjecture saying that two events $y_1, y_2$
n a nonrefocussing globally hyperbolic $(\ss^{m+1}, g), m+1>2,$ that
do not belong to a common null geodesic, are causally unrelated if
and only if the link $(\mathfrak S_{y_1}, \mathfrak S_{y_2})$ is
isotopic through skies to a trivial link. (Probably the best choice
for the trivial link consists of skies of two events on the same
Cauchy surface.)

\begin{cor}\label{Low2true}
Let $(x_1, x_2)$ be two causally unrelated events in a
nonrefocussing globally hyperbolic space-time $(\ss^{m+1}, g),
m+1>2$. Let $(y_1, y_2)$ be two events that do not belong to a
common null geodesic, then the following two statements are
equivalent:
\begin{description}
\item [1] The nonsingular links $(\mathfrak S_{x_1}, \mathfrak S_{x_2})$ and
$(\mathfrak S_{y_1}, \mathfrak S_{y_2})$ are isotopic through skies.
\item [2] The events $y_1, y_2$ are causally unrelated.
\end{description}
\end{cor}

\pp Remark~\ref{comparison} says that for nonrefocussing
globally hyperbolic $(\ss, g)$ two events are sky-isotopic if and
only if their skies are isotopic through skies. Now
Corollary~\ref{Low2true} follows from Theorem~\ref{Low1true}.\qed

\begin{rem}[Isotopies that consist of skies at each time moment]\label{consistingofskies}
Using Theorem~\ref{Low2true} and the proof of the implication $1\implies 2$ of 
Theorem~\ref{Low1true} one can show the following result. 

Let $(\ss^{m+1}, g), m>1$ be a nonrefocussing globally hyperbolic
space-time. Put $\Emb(S^{m-1}\sqcup S^{m-1}, \mathcal N)$ to be the space of 
smooth
embeddings $S^{m-1}\sqcup S^{m-1}\to \mathcal N.$
Let $x_1, x_2\in \ss$ be two causally unrelated points and let $y_1,
y_2$ be two points that do not lie on a common null geodesic. Then $y_1, y_2$ 
are causally unrelated if and only if 
there is an isotopy $r=r(t)=(r_1(t), r_2(t)):[0,1]\to \Emb(S^{m-1}\sqcup 
S^{m-1}, \mathcal N)$ such that $\Im r_i(t)$ is a sky for all $t\in [0,1]$ and 
$\Im r_i(0)=\mathfrak S_{x_i}, \Im r_i(1)=\mathfrak S_{y_i}, i=1,2.$
\end{rem}

{\bf Acknowledgments:} The first author was supported by the
free-term research money from the Dartmouth College. The second
author was supported by the by MCyT, projects BFM 2002-00788 and BFM
2003-02068/MATE, Spain, and by NSF grant 0406311. His visit to
Dartmouth College was supported by the funds donated by Edward
Shapiro to the Mathematics Department of Dartmouth College.

The authors are very thankful to Robert Caldwell, Paul Ehrlich,
Robert Low, Jose Natario, Jacobo Pejsachowicz, Miguel Sanchez and Sergey 
Shabanov for useful discussions. We are grateful to the anonymous referee for 
the valuable comments and especially for the suggestion to work with conformal 
classes of metrics in Theorem~\ref{curvature}.

\appendix

\section{A brief review of contact and Lorentz manifolds}\label{reviewcontact}
\begin{defin}[contact structures and Legendrian submanifolds] Let $Q^{2m-1}$ be
a smooth manifold equipped with a smooth hyperplane field
$\eta=\{\eta^{2m-2}_q\subset T_qQ^{2m-1}\bigm|q\in Q\}.$ This
hyperplane field is called {\em a contact structure,\/} if it can be
locally presented as the kernel of a $1$-form $\alpha$ with
$\alpha\wedge (d\alpha)^{m-1}\neq 0.$

An immersion (respectively an embedding) $f: Z^{m-1}\to Q$ of an
$(m-1)$-dimensional manifold $Z^{m-1}$ into a $(2m-1)$-dimensional
contact manifold $(Q^{2m-1},\eta)$ is called a {\em Legendrian
immersion {\rm (}respectively a Legendrian embedding{\rm )},\/} if
$(df)(T_zZ)\subset \eta_{f(z)},$ for all $z\in Z.$
\end{defin}

\begin{ex}[The contact structure on $ST^*M$]\label{contactSTM}
For a smooth manifold $M^m$ a point $p\in ST^*M$ can be regarded as
a linear functional $\wt p$ on $T_{\pr p}M$ that is defined up to a
multiplication by a positive number. Thus this point $p$ is
completely described by the hyperplane $\ell_p^{m-1}=\ker \wt
p\subset T_{\pr (p)}M$ and by the half-space $T_{\pr (p)}M\setminus
\ell^{m-1}_p$ where $\wt p$ is positive.

The natural contact structure
$$
\eta=\{\eta_p^{2m-2}\subset T_p(STM)^{2m-1}, p\in ST^*M\}
$$
is given by $\eta_p=(d \pr)^{-1}(\ell_p)$.

If $M$ is equipped with a Riemannian metric $\ov g,$ then we can
identify the tangent and the cotangent bundles of $M.$ Thus we can
also identify the spherical tangent bundle with the spherical
cotangent bundle. A smooth map $\gf:Z\to STM$ can be described as
the map $\psi:=\pr\circ \gf:Z\to M$ together with a smooth vector
field $\xi_{z}\in T_{\psi(z)} M, z\in Z,$ where $\xi_z$ points to
the direction $\gf(z).$ It is easy to see that for an
$(m-1)$-dimensional manifold $Z^{m-1}$ the mapping
$$
\CD f:Z @>\gf >> STM @>\cong >> ST^*M
\endCD
$$
is Legendrian exactly when $\xi_{z}$ is $\ov g$-orthogonal to
$d\psi(T_zZ),$ for all $z\in Z.$
\end{ex}

\begin{defin} [Levi-Civita connection on Lorentz manifolds, geodesic,
exponential map, curvature, etc.]\label{difgeomLorentz} Let $(\ss,
g)$ be a Lorentz manifold and let $\Xi(\ss)$ be the space of all
smooth vector fields $\ss\to T\ss$ on $\ss.$ A {\em Levi-Civita
connection\/} on $(\ss, g)$ is a connection $\nabla^g$ such that the
following metric compatibility and torsion free conditions hold for
every $\xi_1, \xi_2, \xi_3\in\Xi(\ss):$
\[ \xi_1g(\xi_2, \xi_3)=g(\nabla ^g_{\xi_1}\xi_2, \xi_3)+g(\xi_2,
\nabla^g_{\xi_1}\xi_3) \text{ and } [\xi_1, \xi_2]=\nabla
^g_{\xi_1}\xi_2-\nabla^g_{\xi_2}\xi_2.
\]
Every Lorentz manifold $(\ss,g)$ admits a unique Levi-Civita
connection, see for example~\cite[page 22]{BeemEhrlichEasley}. When
no confusion can arise we will often use $\nabla$ rather than
$\nabla^g.$ A {\em geodesic\/} $c:(a,b)\to (\ss, g)$ is a smooth
curve such that $\nabla_{ c'} c'=0$ for all of its points.

\m Similar to Riemannian manifolds one can use geodesics to define
the exponential map $\exp_p: T_p\ss\to \ss.$ The map $\exp_p$ is
defined not on the whole $T_p\ss$ but rather on a star-convex with
respect to $0\in T_p\ss$ set in it. There is an open neighborhood
$\wt U$ of $ 0\in T_p\ss$ such that $\exp_p|_{\wt U}$ is a
diffeomorphism onto a neighborhood of $p\in \ss.$ Such $\wt U$ is
called {\em a normal neighborhood.\/}

The {\em curvature $R$\/} of $\nabla$ is a function that assigns to
each pair $\xi_1, \xi_2\in\Xi(\ss)$ a map
$$
R(\xi_1, \xi_2)\!:\Xi(\ss)\to \Xi(\ss),\,R(\xi_1, \xi_2)
\xi_3=\nabla_{\xi_1}\nabla_{\xi_2}\xi_3-\nabla_{\xi_2}\nabla_{\xi_1}\xi_
3-\nabla_{[\xi_1 , \xi_2]}\xi_3.
$$
It is well-known that for
$p\in \ss,$ $R(\xi_1, \xi_2)\xi_3|_p$ depends only on $\nabla^g$ and
on $\xi_1(p), \xi_2(p), \xi_3(p),$ see for example~\cite[page
20]{BeemEhrlichEasley}. Moreover $R(\xi_1, \xi_2)\xi_3|_p$ linearly
depends on $\xi_1(p), \xi_2(p), \xi_3(p).$

A two-dimensional plane $E_p\subset T_p\ss$ is said to be {\em
spacelike\/} if $g|_{E_p}$ is positive definite, it is called {\em
timelike\/} if $g|_{E_p}$ is nondegenerate but it is not positive
definite, and $E_p$ is called {\em null or light-like\/} if
$g|_{E_p}$ is degenerate. Let $E_p$ be timelike or spacelike and let
$v, w$ be a basis of $E_p,$ then one defines the sectional curvature
$$
K(E_p)=\frac{g(R(w,v)v,w)}{g(v,v)g(w,w)-(g(v,w))^2},
$$
see~\cite[pages 29-30]{BeemEhrlichEasley}. (Note that for light-like
$E_p$ the expression in the denominator is zero.)
\end{defin}

\section{Manifolds for which the positively and the negatively oriented
$S^{m-1}$-fibers of $ST^*M\to M^m$ are homotopic}\label{backward}

\begin{defin}[good manifolds]
Let $M^m, m>2$ be a Cauchy surface in a globally hyperbolic $(\ss,
g).$ Let $r: S^{m-1} \to S^{m-1}$ be an autodiffeomorphism of degree
$-1.$ We call a manifold $M$ {\it ``good''\/} if the maps $\eps_v r,
\eps_v:S^{m-1}\to ST^*M $ are not free homotopic. Since for any
$v_1,v_2\in M$ the maps $\eps_{v_1}$ and $\eps_{v_2}$ are free
homotopic, this definition does not depend on the choice of $v\in
M.$
\end{defin}

For a generic cooriented wave front $W_{x, M}$ on $M$ we can reconstruct the
submanifold $\Im \wt W_{x, M}\subset ST^*M$ from the cooriented $\Im
W_{x, M}.$ The submanifold $\Im \wt W_{x, M}$ is diffeomorphic to
$S^{m-1}$ and the lifted wave front $\wt W_{x, M}:S^{m-1}\to ST^*M$
can be reconstructed up to an autodiffeomorphism of $S^{m-1}.$

If $M$ is good, then we can reconstruct $\wt W_{x, M}$ up to an
orientation preserving autodiffeomorphism of $S^{m-1}.$ Indeed,
choose a diffeomorphism $f:S^{m-1}\to \Im \wt W_{x, M}\subset M.$
Since $M$ is good, exactly one of the maps $f$ and $fr$ is homotopic
to $\eps_v$ and this map equals to $\wt W_{x, M}$ up to an
orientation preserving autodiffeomorphism of $S^{m-1}$.

Two links that are the same up to orientation preserving
autodiffeomorphisms of the linked spheres are link homotopic. Since
$\alk$ does not change under link homotopy, we see that for good
$M$ the methods of Examples \ref{example1} and \ref{light} work even
if the front orientations are not specified in the pictures of the
cooriented fronts.

The following theorem shows that almost all manifolds are good.

\begin{thm}\label{reflection} If a connected oriented manifold $M^m$
is not good, then $M$ is homeomorphic to an even-dimensional sphere
and $\Im\{\eps_*: \pi_{m-1}(S^{m-1})\to \pi_{m-1}(M)\}\cong \Z/2$.
\end{thm}

\pp Since $M$ is orientable, the $\pi_1(ST^*M)$-action on the class
in $\pi_{m-1}(ST^*M)$ of the positively oriented $S^{m-1}$-fiber of
$\pr$ is trivial, see the proof of Lemma~\ref{Bconnected}. So $\eps$
and $\eps r$ are homotopic if and only if the group $G:=\Im\{\eps_*:
\pi_{m-1}(S^{m-1})\to \pi_{m-1}(ST^*M)\}$ is $\Z/2$ or 0. Note that
if $M$ is not closed then it is good, since the bundle $ST^*M \to M$
has a section (the Euler class belongs to the trivial group), and
therefore $G=\Z$. So we assume that $M$ is a closed oriented
manifold and consider the following commutative diagram:
$$
\CD
@. \pi_m(M) @>\partial >> \pi_{m-1}(S^{m-1}) @>\eps_*>>\pi_{m-1}(ST^*M)\\
@. @VhVV @V\cong Vh'V @. \\
\Z @= H_m(M) @>\chi(M) >> H_{m-1}(S^{m-1}) @= \Z
\endCD
$$
Here $h$ and $h'$ are the Hurewicz homomorphism, the top sequence is
a segment of
the homotopy exact sequence of the spherical cotangent bundle $ST^*M \to M$,
and the bottom map is the multiplication by the Euler characteristic
of $M$. (The commutativity follows since in the Leray--Serre
spectral sequence of the spherical cotangent bundle the
transgression $\tau: H_m(M) \to H_{m-1}(S^{m-1})$ is the
multiplication by the Euler characteristic $\chi(M)$ of $M$.) Note
that $G\cong \pi_{m-1}(S^{m-1})/\Im \partial$.

If $G=\Z/2,$ then $\Im \partial =2\Z \subset \Z=\pi_{m-1}(S^{m-1}).$
Hence $h$ is a non-zero homomorphism, i.e. there exists a map $S^m
\to M^m$ of non-zero degree. Therefore $M$ is a rational homology
sphere, cf \lemref{oddrationalhomologysphere}. Hence $\chi(M)=0$ if
$m$ is odd and $\chi(M)=2$ if $m$ is even. The case $\chi(M)=0$ is
impossible, since $h'\partial \ne 0$. So $\chi(M)=2$ and therefore
$h$ must be surjective. Thus there exists a map $S^m \to M^m$ of
degree $1.$ Similarly to the proof of~\propref{usefulprop}.(ii), we
get that $M$ is homeomorphic to a sphere. Since $\chi(M)=2,$ $m$ is
even.

If $G=0,$ then $\partial$ is surjective. Hence $h$ must be
surjective and $\chi(M)=1.$ Similarly to the case considered before,
we get that $M$ is a homotopy sphere. However this contradicts to
$\chi(M)=1$. \qed

\forget
\begin{rem}[Cases when it suffices to specify the orientations of the depicted
fronts] Let $c:S^{m-1}\to S^{m-1}$ be the central symmetry
automorphism. When $m$ is odd, $c$ is an automorphism of degree
$-1.$ When $m$ is even, $c$ is free homotopic to the identity map
$\id:S^{m-1}\to S^{m-1}.$ Let $\wt c:ST^*M\to ST^*M$ be the
fiberwise automorphism that is the central symmetry on each
$S^{m-1}$-fiber of $\pr:ST^*M\to M.$

If a front $W_{x, M}$ is generic and $m=\dim M$ is odd, then it
suffices to specify the orientation of the branches of $\Im W_{x,
M}$ in order to be able to reconstruct $\wt W_{x, M}$ up to an
orientation preserving autodiffeomorphism of $S^{m-1}.$ (That is we
can forget about the front coorientation when depicting it.) Indeed
there are two possible ways to equip $W_{x, M}$ with a coorientation
that continuously depends on a point on a branch of $\Im W_{x, M}.$
Thus for a generic $W_{x, M}$ we can reconstruct the (unordered)
pair of oriented submanifolds $\Im \wt W_{x, M}, \Im \wt c\circ \wt
W_{x, M}$ from $\Im W_{x, M}$ equipped with the orientations of the
branches. However we can not immediately tell which one of the two
submanifolds is in fact $\wt W_{x, M}.$

For each of the two oriented submanifolds we pick an orientation
preserving diffeomorphism from $S^{m-1}$ to it. Such a
diffeomorphism is unique up to an orientation preserving
autodiffeomorphism of $S^{m-1}.$ Thus we get two inclusions
$S^{m-1}\to ST^*M$ whose images are $\Im \wt W_{x, M}, \Im \wt
c\circ \wt W_{x, M}.$ If one of these inclusions is homotopic to
$\eps_v,$ then the other one is homotopic to $\eps_v\circ c.$ When
$m$ is odd, $c$ is an automorphism of degree $-1$ and
Theorem~\ref{reflection} implies that $\eps_v$ and $\eps_v\circ c$
are not homotopic. Since $\wt W_{x, M}$ is homotopic to $\eps_v,$ we
can tell which one of the two oriented submanifolds of $ST^*M$ is in
fact $\Im \wt W_{x, M}.$

Similarly to the case considered above, we get that when $M$ is
odd-dimensional and the wave fronts are generic it is not necessary
to depict the front coorientation and it suffices to specify the
orientations of the front branches in order to be able to compute
$alk.$ The orientations of the immersed branches of a front can be
depicted as the orientations of the $1$-dimensional normal bundles
to the branches, i.e. as a vector field normal to the branches.
\end{rem}
\forgotten

\end{document}